\pdfoutput=1
\newif\ifpersonal

\personaltrue 

\ifpersonal
\newcommand*{\personal}[1]{\textcolor[rgb]{0,0,1}{(Personal: #1)}}
\newcommand*{\todo}[1]{\textcolor{red}{(Todo: #1)}}
\else
\newcommand*{\personal}[1]{\ignorespaces}
\newcommand*{\todo}[1]{\ignorespaces}
\fi

\documentclass[oneside, hidelinks, 12pt, a4paper,reqno]{amsart}
\usepackage[toc,page]{appendix}
\usepackage[foot]{amsaddr}
\usepackage{etoolbox}

\makeatletter
\patchcmd{\@setaddresses}{\indent}{\noindent}{}{}
\patchcmd{\@setaddresses}{\indent}{\noindent}{}{}
\patchcmd{\@setaddresses}{\indent}{\noindent}{}{}
\patchcmd{\@setaddresses}{\indent}{\noindent}{}{}
\makeatother

\usepackage{amsthm,mathtools,bm,tensor}
\newcommand*{\boldone}{\text{\usefont{U}{bbold}{m}{n}1}}
\usepackage[T1]{fontenc}\newcommand{\yo}{\text{\usefont{U}{min}{m}{n}\symbol{'110}}}
\DeclareFontFamily{U}{min}{}
\DeclareFontShape{U}{min}{m}{n}{<-> dmjhira}{}
\usepackage[bitstream-charter]{mathdesign}
\usepackage{XCharter}
\usepackage{enumerate,comment,braket,xspace,tikz-cd,csquotes}
\usepackage[driver=dvips,centering,vscale=0.7,hscale=0.8]{geometry}
\usepackage{url}
\usepackage{microtype}
\usepackage[greek.polutoniko,english]{babel}
\usepackage[style=authoryear,backend=biber,sorting=nyt,citestyle=ieee-alphabetic, bibstyle=alphabetic,texencoding=ascii,bibencoding=utf8,maxnames=5,minnames=3,maxbibnames=99]{biblatex}
\usepackage[nodisplayskipstretch]{setspace}
\usepackage{enumitem}
\usepackage{eucal}
\usepackage[scr=boondoxo]{mathalfa}
\usepackage{listings}
\usepackage{colordvi}
\usepackage{tikz,pgf}
\usepackage{pict2e}
\usetikzlibrary{arrows,decorations.pathmorphing,backgrounds,positioning,fit,matrix,calc}
\tikzset{curve/.style={settings={#1},to path={(\tikztostart)
			.. controls ($(\tikztostart)!\pv{pos}!(\tikztotarget)!\pv{height}!270:(\tikztotarget)$)
			and ($(\tikztostart)!1-\pv{pos}!(\tikztotarget)!\pv{height}!270:(\tikztotarget)$)
			.. (\tikztotarget)\tikztonodes}},
	settings/.code={\tikzset{quiver/.cd,#1}
		\def\pv##1{\pgfkeysvalueof{/tikz/quiver/##1}}},
	quiver/.cd,pos/.initial=0.35,height/.initial=0}
\tikzset{tail reversed/.code={\pgfsetarrowsstart{tikzcd to}}}
\tikzset{2tail/.code={\pgfsetarrowsstart{Implies[reversed]}}}
\tikzset{2tail reversed/.code={\pgfsetarrowsstart{Implies}}}
\tikzset{Rightarrow/.style={double equal sign distance,>={Implies},->},
triple/.style={-,preaction={draw,Rightarrow}},
quadruple/.style={preaction={draw,Rightarrow,shorten >=0pt},shorten >=1pt,-,double,double
distance=0.2pt}}
\usepackage{epsfig}
\usepackage{epstopdf}
\usepackage{tikz-cd}
\usepackage{fancyhdr}
\usepackage{graphicx}
\usepackage{fullpage}
\makeatletter
\newcommand{\mylabel}[2]{#2\def\@currentlabel{#2}\label{#1}}
\makeatother
\usepackage{hyperref}
\hypersetup{
	colorlinks=true,
	linkcolor=red,
        linktoc=page,
	anchorcolor=black,
	filecolor=blue,
	citecolor = blue, 
	menucolor=., 
	urlcolor=cyan,
	breaklinks=true,
    unicode=true,
}
\usepackage[capitalize]{cleveref}
\usepackage{nameref}
\usepackage{enumitem}
\crefalias{enumi}{lemman}
\usepackage{comment}
\usepackage{color}
\usepackage[totoc]{idxlayout}
\usepackage[lastexercise]{exercise}
\usepackage{xcolor}
\usepackage{breakcites}
\usetikzlibrary{positioning}
\newcommand{\bDelta}{\bm{\Delta}}
\makeatletter
\newcommand{\twoarrows}[3][0.2ex]{%
  \mathrel{\mathpalette\twoarrows@{{#1}{#2}{#3}}}%
}
\newcommand{\twoarrows@}[2]{\twoarrows@@#1#2}
\newcommand{\twoarrows@@}[4]{%
  \vcenter{\offinterlineskip\m@th
    \ialign{\hfil##\hfil\cr
      $#1#3$\cr
      \noalign{\vskip#2}
      $#1#4$\cr
    }%
  }%
}
\makeatother

\makeatletter
\newcommand*{\doublerightarrow}[2]{\mathrel{
		\settowidth{\@tempdima}{$\scriptstyle#1$}
		\settowidth{\@tempdimb}{$\scriptstyle#2$}
		\ifdim\@tempdimb>\@tempdima \@tempdima=\@tempdimb\fi
		\mathop{\vcenter{
				\offinterlineskip\ialign{\hbox to\dimexpr\@tempdima+1em{##}\cr
					\rightarrowfill\cr\noalign{\kern.5ex}
					\rightarrowfill\cr}}}\limits^{\!#1}_{\!#2}}}
\newcommand*{\triplerightarrow}[1]{\mathrel{
		\settowidth{\@tempdima}{$\scriptstyle#1$}
		\mathop{\vcenter{
				\offinterlineskip\ialign{\hbox to\dimexpr\@tempdima+1em{##}\cr
					\rightarrowfill\cr\noalign{\kern.5ex}
					\rightarrowfill\cr\noalign{\kern.5ex}
					\rightarrowfill\cr}}}\limits^{\!#1}}}
\makeatother

\makeatletter
\def\@tocline#1#2#3#4#5#6#7{\relax
	\ifnum #1>\c@tocdepth 
	\else
	\par \addpenalty\@secpenalty\addvspace{#2}%
	\begingroup \hyphenpenalty\@M
	\@ifempty{#4}{%
		\@tempdima\csname r@tocindent\number#1\endcsname\relax
	}{%
		\@tempdima#4\relax
	}%
	\parindent\z@ \leftskip#3\relax \advance\leftskip\@tempdima\relax
	\rightskip\@pnumwidth plus4em \parfillskip-\@pnumwidth
	#5\leavevmode\hskip-\@tempdima
	\ifcase #1
	\or\or \hskip 1em \or \hskip 2em \else \hskip 3em \fi%
	#6\nobreak\relax
	\dotfill\hbox to\@pnumwidth{\@tocpagenum{#7}}\par
	\nobreak
	\endgroup
	\fi}
\makeatother

\newlist{defenum}{enumerate}{1} 
\setlist[defenum]{label=\arabic*), ref=\thedefn.(\arabic*)}
\crefalias{defenumi}{definition} 
\newlist{lemenum}{enumerate}{1} 
\setlist[lemenum]{label=\arabic*), ref=\thelemma.(\arabic*)}
\crefalias{lemenumi}{lemma} 
\newlist{assumpenum}{enumerate}{1} 
\setlist[assumpenum]{label=\arabic*), ref=\theassumption.(\arabic*)}
\crefalias{assumpenumi}{assumption} 
\newlist{propenum}{enumerate}{1} 
\setlist[propenum]{label=\arabic*), ref=\theproposition.(\arabic*)}
\crefalias{propenumi}{proposition}
\newlist{exmpenum}{enumerate}{1} 
\setlist[exmpenum]{label=\arabic*), ref=\theexmp.(\arabic*)}
\crefalias{exmpenumi}{example}
\newlist{thmenum}{enumerate}{1} 
\setlist[thmenum]{label=\arabic*), ref=\thethm.(\arabic*)}
\crefalias{thmenumi}{theorem}
\newlist{factenum}{enumerate}{1} 
\setlist[factenum]{label=\arabic*), ref=\thefact.(\arabic*)}
\crefalias{factenumi}{fact}
\newlist{exercenum}{enumerate}{1} 
\setlist[exercenum]{label=\arabic*), ref=\theexercise.(\arabic*)}
\crefalias{exercenumi}{exercise}
\newlist{remarkenum}{enumerate}{1} 
\setlist[remarkenum]{label=\arabic*), ref=\theremark.(\arabic*)}
\crefalias{remarkenumi}{remark} 
\newlist{axiomenum}{enumerate}{1} 
\setlist[axiomenum]{label=\arabic*), ref=\theaxiom.(\arabic*)}
\crefalias{axiomenumi}{axiom} 
\newlist{questenum}{enumerate}{1} 
\setlist[questenum]{label=\arabic*), ref=\thequestion.(\arabic*)}
\crefalias{questenumi}{question}

   \makeatletter
\newcommand{\adjunction}{\@ifstar\named@adjunction\normal@adjunction}
\newcommand{\normal@adjunction}[4]{%
  #1\colon #2%
  \mathrel{\vcenter{%
    \offinterlineskip\m@th
    \ialign{%
      \hfil$##$\hfil\cr
      \longrightharpoonup\cr
      \noalign{\kern-.3ex}
      \smallbot\cr
      \longleftharpoondown\cr
    }%
  }}%
  #3 \noloc #4%
}
\newcommand{\named@adjunction}[4]{%
  #2%
  \mathrel{\vcenter{%
    \offinterlineskip\m@th
    \ialign{%
      \hfil$##$\hfil\cr
      \scriptstyle#1\cr
      \noalign{\kern.1ex}
      \longrightharpoonup\cr
      \noalign{\kern-.3ex}
      \smallbot\cr
      \longleftharpoondown\cr
      \scriptstyle#4\cr
    }%
  }}%
  #3%
}
\newcommand{\longrightharpoonup}{\relbar\joinrel\rightharpoonup}
\newcommand{\longleftharpoondown}{\leftharpoondown\joinrel\relbar}
\newcommand\noloc{%
  \nobreak
  \mspace{6mu plus 1mu}
  {:}
  \nonscript\mkern-\thinmuskip
  \mathpunct{}
  \mspace{2mu}
}
\newcommand{\smallbot}{%
  \begingroup\setlength\unitlength{.15em}%
  \begin{picture}(1,1)
  \roundcap
  \polyline(0,0)(1,0)
  \polyline(0.5,0)(0.5,1)
  \end{picture}%
  \endgroup
}
\makeatother

\newcommand{\infinity}{\mbox{\footnotesize $\infty$}}

\newcommand{\triv}{\operatorname{triv}}
\newcommand{\oblv}{\operatorname{oblv}}

\newcommand{\Einf}{\mathbb{E}_{{\scriptstyle\infty}}}

\newcommand{\Bbbbk}{\Bbbk}

\newcommand{\op}{\operatorname{op}}

\newcommand{\Shv}{\mathrm{Shv}}

\newcommand{{\Cn}}{\operatorname{Cn}}

\newcommand{\RR}{\mathbb{R}}

\DeclareMathOperator*{\colim}{colim}

\newcommand{\Homin}{{\underline{\smash{\Hom}}}}

\newcommand{\calA}{\mathcal{A}}

\newcommand{\scrC}{\mathscr{C}}
\newcommand{\scrD}{\mathscr{D}}
\newcommand{\scrE}{\mathscr{E}}

\newcommand{\scrG}{\mathscr{G}}
\newcommand{\scrO}{\mathscr{O}}

\newcommand{\scrP}{\mathscr{P}}

\newcommand{\Hom}{{\operatorname{Hom}}}

\DeclareMathOperator*{\limast}{lim^{\ast}}
\DeclareMathOperator*{\limsh}{lim^{!}}

\DeclareMathOperator*{\colimsh}{colim_{!}}
\DeclareMathOperator*{\colimastd}{colim^{\ast,\vee}}

\DeclareMathOperator{\Spec}{Spec}

\DeclareMathOperator{\Qcoh}{QCoh}
\DeclareMathOperator{\QCoh}{QCoh}

\newcommand{\Mod}{{\operatorname{Mod}}}

\newcommand{\lp}{\left(}
\newcommand{\rp}{\right)}

\newcommand{\scrV}{\mathscr{V}}

\newcommand*{\longhookrightarrow}{\ensuremath{\lhook\joinrel\relbar\joinrel\rightarrow}}

\makeatletter

\def\enddoc@text{}
\makeatother
\newcounter{savedchapter}
\preto\appendix{\setcounter{savedchapter}{\arabic{chapter}}}
\newcommand\resumechapters{
	\setcounter{chapter}{\arabic{savedchapter}}
	\setcounter{section}{0}
	\gdef\@chapapp{\chaptername}
	\gdef\thechapter{\@arabic\c@chapter}
}
\makeatother

\newcommand{\Ebb}{\mathbb{E}}

\newcommand{\PP}{\mathbb{P}}

\newcommand{\Map}{{\operatorname{Map}}}

\newcommand{\Fun}{{\operatorname{Fun}}}

\newcommand{\FunL}{\Fun^{\operatorname{L}}}

\newcommand{\Rani}{\mathrm{Ran}_{\iota}}

\newcommand{\CAlg}{{\operatorname{CAlg}}}

\newcommand{\Alg}{{\operatorname{Alg}}}

\newcommand{\Nuccat}{\Nuc^{\mathrm{cat}}}

\numberwithin{equation}{subsection}
\theoremstyle{plain}
\newtheorem{theorem}[equation]{Theorem}

\newtheorem{lemma}[equation]{Lemma}
\newtheorem{proposition}[equation]{Proposition}

\newtheorem{corollary}[equation]{Corollary}

\newtheorem{question}[equation]{Question}

\theoremstyle{definition}
\newtheorem{defn}[equation]{Definition}
\newtheorem{definition}[equation]{Definition}
\newtheorem{parag}[equation]{}

\newtheorem{remark}[equation]{Remark}
\newtheorem{construction}[equation]{Construction}
\newtheorem{porism}[equation]{Porism}

\newtheorem{warning}[equation]{Warning}
\newtheorem{notation}[equation]{Notation}
\newtheorem{assumption}[equation]{Assumption}
\newtheorem{stassumption}[equation]{Standing assumption}
\newtheorem{fact}[equation]{Fact}

\newtheorem{exmp}[equation]{Example}
\newtheorem*{warning*}{Warning}
\newtheorem*{assumption*}{Standing assumption}
\newcommand{\Sp}{\mathscr{S}\mathrm{p}}

\newcommand\restr[2]{{
		\left.\kern-\nulldelimiterspace
		#1
		\vphantom{\|}
		\right|_{#2} 
}}
\usepackage{lipsum}
\usepackage{fancyhdr}
\providecommand{\abstract}{}
\usepackage{abstract}
\DeclareGraphicsExtensions{.png,.jpg,.pdf,.mps}
\usepackage{fancyhdr}
\usepackage{lastpage}
\usepackage{multirow}
\usepackage{faktor}

\usepackage{array}
\usepackage{makecell}

\newcommand{\stackspace}{1.7}
\newcommand{\stack}[2][1cm]{\;\tikz[baseline, yshift=.65ex]%
	{\foreach \k [evaluate=\k as \r using (.5*#2+.5-\k)*\stackspace] in {1,...,#2}{%
			\ifodd\k{\draw[->](0,\r pt)--(#1,\r pt);}%
			\else{\draw[<-](0,\r pt)--(#1,\r pt);}\fi
	}}\;}
    \newcommand{\stackrev}[2][1cm]{\;\tikz[baseline, yshift=.65ex]%
	{\foreach \k [evaluate=\k as \r using (.5*#2+.5-\k)*\stackspace] in {1,...,#2}{%
			\ifodd\k{\draw[<-](#1,\r pt) -- (0,\r pt);}%
			\else{\draw[->](#1,\r pt) -- (0,\r pt);}\fi
	}}\;}

\usepackage{doc}

\newcommand{\CondAni}{\mathrm{Ani}^{\mathrm{cond}}}
\newcommand{\Prof}{\mathrm{ProFin}}
\newcommand{\Betti}{\mathrm{Betti}}
\newcommand{\cf}{\mathscr{f}}

\newcommand{\Anring}{\mathrm{AnRing}}
\newcommand{\AnStacks}{\mathrm{AnStacks}}

\newcommand{\Ansp}{\mathrm{AnSpec}}

\newcommand{\Afd}{\mathrm{Afd}}

\newcommand{\Spa}{\mathrm{Spa}}

\newcommand{\catrex}{\smash{\widehat{\mathrm{Cat}}}^{\mathrm{rex}}}

\newcommand{\LPr}{\mathrm{Pr}^{\mathrm{L}}}

\newcommand{\Ani}{\mathrm{Ani}}

\newcommand{\Nuc}{\mathrm{Nuc}}
\newcommand{\RFun}{\mathrm{Fun}^{\mathrm{R}}}

\newcommand{\cC}{\mathscr{C}}

\newcommand{\cV}{\mathscr{V}}

\newcommand{\Z}{\mathbb{Z}}

\newcommand{\N}{\mathbb{N}}

\usepackage[normalem]{ulem}

\title{An axiomatic approach to analytic $1$-affineness}
\author{Matteo Montagnani \and Emanuele Pavia}
\address{SISSA\\Via Bonomea 265\\34136 Trieste, TS, Italy}\email{\href{mailto:mmontagn@sissa.it}{mmontagn@sissa.it}}
\email{\href{mailto:epavia@sissa.it}{epavia@sissa.it}}
\begin{document}
\maketitle
\begin{abstract}
The notion of $1$-affineness was originally formulated by Gaitsgory in the context of derived algebraic geometry. Motivated by applications to rigid and analytic geometry, we introduce two very general and abstract frameworks where it makes sense to ask for objects to be $1$-affine with respect to some sheaf of categories. The first framework is suited for studying the problem of $1$-affineness when the sheaf of categories arises from an operation in a six-functor formalism over $\scrC$; we apply it to the setting of analytic stacks and condensed mathematics. The second one concerns $1$-affineness in the context of quasi-coherent sheaves of categorical modules over stable module categories: it simultaneously generalizes the algebro-geometric setting of Gaitsgory and makes it possible to formulate the problem also when dealing with rigid analytic varieties and categories of nuclear modules. 
\end{abstract}

\tableofcontents
{\section*{Introduction}
\addtocontents{toc}{\protect\setcounter{tocdepth}{0}}
\numberwithin{equation}{section}
\renewcommand\thesection{I}
\renewcommand\theequation{\thesection.\arabic{equation}}
The notion of sheaves of categories provides a categorification of the theory of quasicoherent sheaves from ordinary algebraic geometry. This is the suitable technical framework where to study higher categorical structures which vary "continuously" over some geometric object (with respect to some Grothendieck topology).

This notion has been first developed in the seminal work \cite{1affineness}: in \textit{loc. cit.}, Gaitsgory defines a functor
\[
\mathrm{ShvCat}\colon\mathrm{PSt}_{\Bbbk}\longrightarrow\catrex
\]
which assigns to any derived prestack $X$ (defined over a ground ring $\Bbbk$) a cocomplete $\infinity$-category consisting of “quasi-coherent sheaves of categories” on $X$. \\
Roughly speaking, an object of
$\mathrm{ShvCat}(X)$ is a rule that assigns to each affine test scheme $\mathrm{Spec}(A)\to X$ a presentable $\Qcoh(A)$-linear $\infinity$-category, together with a homotopy-coherent system of descent data that allows to glue such $\infinity$-categories along any fiber product of affine schemes over $X$. In this way, we are able to work with higher categorical data varying over a base $X$ in a way which is conceptually very close to what one is used to in ordinary algebraic geometry when one works with quasi-coherent sheaves of modules.\\
In particular, for a quasi-coherent sheaf of categories we can define operations which generalize familiar notions such as pullbacks, pushforwards, relative tensor products, which remain sensitive to the geometry of our base prestack.\\

The key concept of \cite{1affineness} is the notion of \textit{$1$-affineness}. Just like quasi-coherent sheaves on an affine scheme $\Spec(A)$ are completely determined by the $A$-module structure on their global sections -- this property actually characterizes affine schemes among all quasi-separated schemes -- a $1$-affine (pre)stack is defined to be a (pre)stack for which any quasi-coherent sheaf of categories is completely determined by some $\QCoh(X)$-linear $\infinity$-category. More precisely, we have a functor
\begin{equation*}
    \Gamma_{X}\colon \mathrm{ShvCat}(X) \longrightarrow \Mod_{\Qcoh(X)}(\LPr)
\end{equation*}
sending a sheaf of categories
$\scrC \in \mathrm{ShvCat}(X)$ to its $\infinity$-category of global sections
$\Gamma_{X}(\scrC)$
 which carries a natural action of the monoidal $\infinity$-category
$\Qcoh(X)$. Then, the prestack $X$ is said to be \emph{$1$-affine} if this functor is an equivalence of $\infinity$-categories. 

Establishing whether a given prestack is $1$-affine is a highly nontrivial problem and carries deep geometric and categorical implications. It allows to study categorified sheaf-theoretic phenomena by reducing to the study of a single $\infinity$-category of modules over a (hopefully well-understood) symmetric monoidal $\infinity$-category.\\
As such, the notion of $1$-affineness has been used for several geometric applications: for example, in \cite{stefanich2023tan}, Stefanich uses the notion of $1$-affineneness to improve on Lurie's Tannaka reconstruction theorem, discarding the tameness hypotheses of \cite{sag}. Along the way, he proves the $1$-affineness of a large class of stacks -- namely, the quasi-compact stacks with quasi-affine diagonal.
The notion of 1-affineness has also been studied in the context of Hecke categories and topological field theories, see for example \cite{benzvi2012morita} and \cite{Ben_Zvi_2018}.

However, as explained in \cite{1affineness} the main motivation for the notion of sheaf of categories and 1-affineness comes from the geometric Langlands program. Indeed, in the geometric Langlands program it is interesting to study $\infinity$-categories with an action of the loop group $G(\!(t)\!)$. Furthermore, the $\infinity$-category of all $\infinity$-categories with an action of $G(\!(t)\!)$ is closely related to the $\infinity$-category of quasi-coherent sheaves of categories $\mathrm{ShvCat}(\mathbf{B}G(\!(t)\!)_\mathrm{dR})$, where $\mathbf{B}G(\!(t)\!)_\mathrm{dR}$ is the de Rham stack associated to the classifing stack $\mathbf{B}G(\!(t)\!)$. The proof of the geometric Langlands conjecture \cite{GLC1, GLC2, GLC3, GLC4,GLC5} actually relies on the facts that $X_{\mathrm{dR}}$ is $1$-affine when $X$ is an ind-scheme locally of finite type, and that $\mathbf{B}G$ is 1-affine when $G$ is the formal completion of an affine algebraic group along a subgroup.\\

 The aim of this paper is to introduce the notions of sheaves of categories and $1$-affineness in the setting of rigid and analytic geometry. Our trivial observation is that the problem of $1$-affineness can be posed in a very abstract generality, and some affirmative results do not really depend on the fact that we are dealing with algebraic or analytic stacks but are formal consequences of the theory. For this reason, we study the notion of $1$-affineness from a somewhat axiomatic perspective: we allow $\scrC$ to be a general $\infinity$-site, with some reasonable assumptions on the $\infinity$-category and on the Grothendieck topology (see \cref{assumption:standing}). We then consider a sheaf of $\scrV$-linear presentable symmetric monoidal
$\infinity$-categories on $\scrC$
\begin{equation*}
    \scrD \colon  \scrC \longrightarrow \mathrm{CAlg}(\Mod_{\cV}(\LPr)),
\end{equation*}
which again satisfy some reasonable assumptions. We can thus define the "categorification” of $\scrD$ to be the functor
\begin{equation*}
    \begin{split}
\scrD^{\mathrm{cat}}\colon\scrC^{\op}&\longrightarrow\catrex\\
U&\mapsto\LPr_{\scrD(U)}.
\end{split}
\end{equation*}
The first problem we study is to find reasonable conditions on $\scrC$ and on $\scrD$ for which the functor $\scrD^{\mathrm{cat}}(-)$ is a sheaf for the Grothendieck topology on $\scrC$.
In particular, we show that this holds in the following cases.
\begin{enumerate}
\item The sheaf $\scrD$ is already the categorification of a sheaf of
commutative algebras
\begin{align*}
D\colon\scrC^{\op}\longrightarrow\CAlg(\cV),
\end{align*}
in the sense that we have a natural equivalence of functors
\[
\scrD\simeq\Mod_{D(-)}(\cV),
\]
and that moreover all the maps in $S$ are descendable, as in
\Cref{def: descendable}. (The problem of $1$-affineness studied by Gaitsgory falls in this setting.)
\item The sheaf of $\infinity$-categories $\scrD$ is obtained as one operation in a strongly
monoidal six-functor formalism
\[
(\scrD^\ast,\scrD_!)\colon
\mathrm{Corr}(\scrC)^{\otimes}_{S',\mathrm{all}}\longrightarrow
\mathrm{Pr}^{\mathrm{L,\otimes}}_{\cV}
\]
in the sense of \cite{mann2022p}, while the topology on $\scrC$ is the universal
!-able topology (see \Cref{def:!-able}). 
\end{enumerate}
When we have a fully faithful dense functor
$\iota \colon \scrC \to \widehat{\scrC}$, we have two different ways to extend $\scrD^{\mathrm{cat}}$ to an assignment on all of $\widehat{\scrC}$: we can either take the right Kan extension of $\scrD^{\mathrm{cat}}$ along $\iota$ ($\Rani\scrD^{\mathrm{cat}}(-)$), \textit{or} we can take categorical modules over the right Kan extension of $\scrD$ along $\iota$ ($\Mod_{\Rani\scrD(-)}(\LPr)$). These two operations are \textit{a priori} very different, but for any object
$Y$ of $\widehat{\scrC}$, we have a natural comparison functor
\begin{align}
\label{functor:1_affineness_intro}
\Mod_{\mathrm{Ran}_{\iota}\scrD(Y)}(\LPr)\longrightarrow\mathrm{Ran}_{\iota}\scrD^{\mathrm{cat}}(Y).
\end{align}
induced by the universal property of limits.
\begin{defn}
\label{def:1_affineness_intro}
Let $\scrC$, $\cV$ and $\scrD$ be as in \cref{assumption:standing}, and let $\iota\colon\scrC\to\widehat{\scrC}$ be a dense functor. We say that an object $Y$ in $\widehat{\scrC}$ is \textit{$1$-affine with respect to $\iota$} if the comparison functor \eqref{functor:1_affineness_intro} is an equivalence.
\end{defn}
We put ourselves this general setup motivated by two different problems, both concerning $1$-affineness in rigid or analytic geometry. The first one (studied in \cref{sec:modular_case}) is related to
rigid analytic geometry and $K$-theory. Indeed, thanks to the work of
Clausen and Scholze in \emph{condensed mathematics}, (\cite{clausen2019lectures}) and to the definition of Efimov of \emph{continuous $K$-theory} for dualizable $\infinity$-categories (\cite{efimov2025ktheorylocalizinginvariantslarge,efimov2025localizinginvariantsinverselimits}), nowadays we have a reasonable notion of algebraic
$K$-theory for “analytic spaces”. Indeed to every rigid analytic variety $X$ it is
possible to associate, functorially, the $\infinity$-category $\Nuc(X)$ of nuclear sheaves on $X$. The “algebraic $K$-theory” of $X$ is then defined to be the continuous $K$-theory of the $\infinity$-category $\Nuc(X)$. For this reason, at least from a $K$-theoretic point of view, we can consider the $\infinity$-category of nuclear sheaves as the analytic analogue of
the $\infinity$-category of quasi-coherent sheaves in algebraic geometry. Some of these ideas have already been studied in rigid analytic geometry, see for example \cite{andreychev2023ktheorieadischerraume, andreychev2021pseudocoherent}. In this setting, the site $\scrC$ is then defined to be the $\infinity$-category of affinoid
spaces, while the $\infinity$-category $\widehat{\scrC}$ is the $\infinity$-category of rigid analytic varieties,
and the sheaf $\scrD(-)$ is just the functor $\Nuc(-)$. In this way, we can prove the following theorem.
\begin{theorem}[{\cref{thm: 1-affine Nuc}}]
    Let $X$ be a quasi-compact and separated rigid analytic variety. Then $X$ is $1$-affine, i.e., there is an equivalence of categories
    \begin{equation*}
        \Mod_{\Nuc(X)}(\LPr_{\mathrm{st}}) \simeq \lim_{\Spa(A) \subset X} \Mod_{\Nuc(A)}(\LPr_{\mathrm{st}}),
    \end{equation*}
    where the limit ranges over an affinoid open cover of $X$. 
\end{theorem}
Following the ideas explained in \cite{stefanich2023tan} we were also able to deduce the following K\"unneth formula for nuclear categories in rigid analytic geometry.
\begin{corollary}[{\cref{cor: Kunneth nuclear rigid}}]
    Let $X \to Z$ and $Y \to Z$ be two maps of quasi-compact and separated rigid analytic varieties. Then the canonical functor
    \begin{equation*}
        \Nuc(X)\otimes_{\Nuc(Z)} \Nuc(Y) \longrightarrow \Nuc(X\times_{Z}Y) 
    \end{equation*}
    is an equivalence of $\infinity$-categories.
\end{corollary}

The second problem which we want to investigate is $1$-affineness for \emph{analytic stacks}, as defined in \cite{Clausen_Scholze_lectures}: this fits in the abstract setting we study in \cref{sec:six_functor}. Analytic stacks were recently introduced by Clausen and Scholze using the novel theoretic framework of \textit{condensed mathematics} (\cite{clausenlectures, analytic, clausen2019lectures}). Using this formalism they manage to capture a wide array of spaces appearing in complex geometry, non-archimedean geometry, algebraic geometry, and topology. They also come with a well-behaved notion of quasi-coherent sheaves, which enjoy all the functorial properties of a \textit{six-functor formalism} (\cite{six-functors-scholze, heyer6ff2024}).
We also remark that analytic geometry can be studied also following different approaches, which have been currently developed in \cite{benbassat2024perspectivefoundationsderivedanalytic, Bambozzi_2018, soor2024sixfunctorformalismquasicoherentsheaves, kelly2025}. This parallel line of work still exhibits notable structural affinities with the condensed approach.

To study the property of $1$-affineness for analytic stacks, we can consider $\scrC$ to be the $\infinity$-category of affine analytic stacks (modeled by analytic rings in the sense of \cref{def:analytic_ring}), and $\widehat{\scrC}$ to be the $\infinity$-category of all analytic stacks. The sheaf of $\infinity$-categories $\scrD(-)$ is precisely the sheaf sending an analytic stack to its derived $\infinity$-category, hence our choice of notation.\\
However, our assumptions are sufficiently general that we expect that this setting can be used to successfully explore the notion of $1$-affineness in the context of motivic six-functor formalisms for schemes; we leave this application for future work.\\

Our main result in this context is the following.
\begin{theorem}[{\cref{thm: main theorem}}]
\label{thm:main theorem intro}
    Let $X$ be an analytic stack admitting an affine universal !-cover $        \{\Ansp(A_{i}) \to X \}_{i \in I}$ such that all fiber products $\Ansp(A_{i_{1}})\times_{X}\dots \times_{X}\Ansp(A_{i_{n}})$ are again affine. Then $X$ is $1$-affine.
\end{theorem}

We need to stress that \cref{thm:main theorem intro}
can also be obtained using some results contained in \cite{kunneth6ff}, which was written and appeared on ArXiv at the final stages of preparation of the present work. While the stream of ideas is very similar, we obtain a slightly more general $1$-affineness results that holds under the hypotheses of \cref{assumption:modular_six_functor}, see \cref{thm:main 2-descent}.

The latter is one of the main motivations for this paper. In \cite{PPS_2} the authors studied $1$-affineness of a different incarnation of the Betti stack associated to topological cases -- that we call \textit{homotopical} Betti stack, since it captures only the underlying anima of a topological space. There, it is shown that the homotopical Betti stack often fails to be $1$-affine, and that a somewhat surprising topological
obstruction to $1$-affineness lies in any non-trivial element of the second homotopy group. 

On the other hand, in analytic geometry we have an analytic incarnation of the Betti stack of a topological space $X$, which is sensitive to the whole topology of $X$ (indeed, its derived $\infinity$-category agrees with the derived $\infinity$-category of sheaves over $X$ with values in derived condensed abelian groups). This was recently introduced in \cite{geometrizationreallanglands}.

We discover a fundamental difference between the “analytic” Betti stack and the "homotopical" Betti stack from algebraic geometry from the perspective of $1$-affineness, in the sense that the topological constraints are far less strict. 
\begin{theorem}[\cref{thm: Betti is 1-affine}]
    Let $X$ be a finite dimensional, metrizable, compact Hausdorff space. Then the analytic Betti stack associated to $X$ is $1$-affine.
\end{theorem}
Using the Riemann--Hilbert correspondence of \cite{geometrizationreallanglands}, we can relate the analytic Betti stack with the analytic De Rham stack over the complex numbers when $X$ is a compact complex manifold. In this way, we obtain an analytic analogue of the $1$-affineness for the de Rham stack in derived algebraic geometry from \cite{1affineness}.
\begin{theorem}[{\cref{cor:de_rham_1_affine}}]
    Let $X$ be a compact complex manifold. Then the analytic De Rham stack $X_{\mathrm{dR}}^{\mathrm{an}}$ is $1$-affine. 
\end{theorem}

}
\section*{Conventions and notations}\
\begin{itemize}
\item We will use throughout the language of $\infinity$-categories and higher homotopical algebra, as developed in \cite{htt,HA}, from which we borrow most of the notations and conventions. The only exception is provided by the $\infinity$-category of spaces: we prefer the terminology "anima", and denote it as $\Ani$.
\item Since our work heavily relies on intrinsically derived and homotopical concepts, we shall simply write ``limits'', ``colimits'', ``tensor products'', etc., suppressing adjectives such as ``homotopy'' or ``derived'' in the notation. 
\item We denote the subcategory of $\Einf$-algebras in a symmetric monoidal $\infinity$-category $\mathscr{C}$ as $\CAlg(\mathscr{C})$. When $\scrC$ is the $\infinity$-category of $\Bbbk$-modules over a commutative ring spectrum $\Bbbk$, we shall simply write $\CAlg_{\Bbbk}$ instead of $\CAlg(\Mod_{\Bbbk})$. 
\item We denote with $\catrex$ the (very large) $\infinity$-category of large cocomplete $\infinity$-categories, and we denote with $\LPr$ its large sub-$\infinity$-category spanned by presentable $\infinity$-categories. Both these $\infinity$-categories are symmetric monoidal under Lurie's tensor product (\cite[$\S4.8.1$]{HA}). So, for any $0\leqslant k\leqslant\infinity$, we call an object $\cV$ in $\Alg_{\Ebb_{k}}(\LPr)$ a \textit{presentably $\Ebb_k$-monoidal $\infinity$-category}: this is an $\Ebb_k$-monoidal $\infinity$-category $\cV$ which is presentable and such that the tensor product is compatible with colimits separately in each variable. When $k=\infinity$, we simply say that $\cV$ is a \textit{presentably symmetric monoidal $\infinity$-category}.
\item For any presentably  symmetric monoidal $\infinity$-category $\cV$, we write $\LPr_{\cV}$ for the $\infinity$-category $\Mod_{\cV}(\LPr)$ of \textit{presentably $\cV$-linear $\infinity$-categories}. When $\cV=\Sp$ is the $\infinity$-category of spectra, we shall simply write $\LPr_{\mathrm{st}}$ since this is the $\infinity$-category of stable presentable $\infinity$-categories. Again, the relative tensor product over $\cV$ yields a symmetric monoidal structure on $\LPr_{\cV}$, so we can consider $\Ebb_k$-algebra objects in it: they are presentably $\Ebb_k$-monoidal $\infinity$-categories which are tensored and enriched over $\cV$ (\cite{enrichment}). We call a commutative algebra in $\LPr_{\cV}$ a \textit{presentably $\cV$-linear symmetric monoidal $\infinity$-category}.
\item Most of the time we will consider categories enriched over some preferred symmetric monoidal $\infinity$-category. If $\scrC$ is enriched over a category $\cV$, we adopt the following notation. 
\begin{itemize}
  \item We let $\Map_{\scrC}(-,-)$ denote the underlying anima of maps between two objects in $\cC$;
  \item We let $\Homin_{\scrC}(-,-)$ denote the morphism object in $\cV$ providing the enrichment.
\end{itemize}
The latter applies in particular when $\scrC$ is closed symmetric monoidal and thus is enriched over itself. This convention allows us to distinguish whether we are viewing a morphism object 
as a mere anima or as a more structured object.


\item When we will need the theory of condensed mathematics, we will follow the notations and terminology from the video course \cite{Clausen_Scholze_lectures}. For the reader’s convenience, we will still recall all relevant definitions and results and introduce the necessary notations as they appear.
\end{itemize}
\section*{Acknowledgements}
The authors wish to express their gratitude to Mauro Porta and Nicolò Sibilla for many valuable discussions and for their constant interest and encouragement during the preparation of this work. Thanks are also due to Youshua Kesting for helpful conversations and for generously sharing a preliminary version of his paper \cite{kunneth6ff}. The first author also acknowledges Ken Lee and Qixiang Wang for insightful suggestions.

\addtocontents{toc}{\protect\setcounter{tocdepth}{2}}
\numberwithin{equation}{subsection}
\section{Categorified sheaves over \texorpdfstring{$\infinity$}{oo}-sites}
\label{sec:1affinenessgeneral}
\subsection{The general setup}
\begin{stassumption}
\label{assumption:standing}
Throughout this section, we fix the following.
\begin{assumpenum}
\item \label{assumption:site}A small $\infinity$-site $\scrC$ admitting all finite coproducts, binary products, and pullbacks. We call objects in $\scrC$ \textit{test spaces}, and denote them as $U\to X$: they have to be thought as open quasi-compact subobjects of some geometric object $X$ (which may, or may not, belong to the $\infinity$-site $\scrC$ itself). When dealing with morphisms of test spaces, we shall trim our notations and just write $\varphi\colon U\to V$.
\item A stable and presentably symmetric monoidal $\infinity$-category $\cV$.
\item\label{assump:stable} A sheaf of presentably $\cV$-linear symmetric monoidal $\infinity$-categories $\scrD\colon\scrC^{\op}\to\CAlg(\LPr_{\scrV})$ which sends finite limits existing in $\scrC$ to finite colimits in $\CAlg(\LPr_{\cV})$. 
\end{assumpenum}We moreover assume the following.
\begin{assumpenum}
\setcounter{assumpenumi}{3}
\item\label{assumption:disjoint_universal}Finite coproducts in $\scrC$ are \textit{disjoint} and \textit{universal}: let $\varnothing$ denote the initial object of $\scrC$. Then for any finite collection of objects $\left\{U_i\to X\right\}_{i\in I}$, we have that
\[
U_{i_1}\times_{U_{i_1}\coprod U_{i_2}}U_{i_2}\simeq \varnothing.
\]
Moreover, for any diagram of test spaces $U_i\to W\rightarrow V$, the natural map
\[
\coprod_{i\in I}U_i\times_WV\longrightarrow\lp\coprod_{i\in I}U_i\times_WV\rp
\]
is an equivalence.
\item\label{assumption:topology} The Grothendieck topology on $\scrC$ is generated by a collection of maps $S$ which contains all equivalences and is stable under compositions, pullbacks and finite coproducts, in the sense of \cite[Proposition A.$3.2.1$]{sag}. 
\end{assumpenum}
\end{stassumption}
\begin{remark}
Let us say something more on the set of conditions in \cref{assumption:standing}.
\begin{enumerate}
\item The $\infinity$-site admits all finite coproducts, but the condition of having all \textit{binary} products in \cref{assumption:site} means only that whenever $U\to X$ and $V\to X$ are test spaces, then there exists a test space $U\times_X V$ in $\scrC$. In particular, $\scrC$ can fail to have a terminal object: for example, this happens when $\scrC$ is the little site spanned by a base of quasi-compact open subsets of a topological space $X$ which is not quasi-compact itself.
\item The assumptions \ref{assumption:disjoint_universal} and \ref{assumption:topology} identify the Grothendieck topology on $\scrC$ as the one described as follows: a collection of maps $\left\{U_i \to V\right\}_{i\in I}$ is a covering if and only if there exists a \textit{finite} subset $J\subseteq I$ such that $\coprod_{j\in J}U_j\to V$ is a morphism in $S$. In particular, each test space is tautologically quasi-compact, in the sense that any $\tau$-covering can be refined to a finite one.
\item Since $\cV$ is assumed to be stable and presentable, it is a $\Sp$-module in $\LPr$ (\cite[Proposition $4.8.2.18$]{HA}). Being $\cV$ symmetric monoidal, \cite[Corollary $4.8.2.19$]{HA} implies that there exists a symmetric monoidal functor $\Sp\to\cV$, hence we obtain a forgetful functor
\[
\LPr_{\cV}\longrightarrow\LPr_{\Sp}\simeq\LPr_{\mathrm{st}},
\]
where $\LPr_{\mathrm{st}}$ is the $\infinity$-category of presentable and stable $\infinity$-categories. In particular, \cref{assump:stable} implies that for all test spaces $U\to X$ the $\infinity$-category $\scrD(U)$ is stable. 
\end{enumerate}
\end{remark}
\begin{defn}
\label{def:categorification}
Let $\scrC$, $\cV$ and $\scrD$ be as in \cref{assumption:standing}. The \textit{categorification of $\scrD$} is the functor
\begin{equation}
\label{functor:categorification}
\begin{split}
\scrD^{\mathrm{cat}}\colon\scrC^{\op}&\longrightarrow\catrex\\
U&\mapsto\LPr_{\scrD(U)}.
\end{split}
\end{equation}
\end{defn}
We are interested in extending the categorification of the sheaf $\scrD$ to some (possibly large) $\infinity$-category $\widehat{\scrC}$ containing $\scrC$, in such a way that objects in $\widehat{\scrC}$ are still governed by the ones in $\scrC$. In order to formalize this idea, we need to fix some notations.
\begin{defn}
\label{def:dense_cat}
Let $\scrC$ be a $\infinity$-category, and let $\iota\colon\scrC\to\widehat{\scrC}$ be a functor. We say that $\iota$ is \textit{dense} if the identity functor on $\widehat{\scrC}$ is a left Kan extension of $\iota$ along itself. Equivalently, this means that for all objects $D$ of $\widehat{\scrC}$ the functor
\[
\scrC\times_{\widehat{\scrC}}\widehat{\scrC}_{/D}\longrightarrow\scrC\overset{\iota}{\longrightarrow}\widehat{\scrC}
\]
admits a colimit, which is equivalent to $D$. If $\iota$ is moreover fully faithful, we shall say that $\iota$ is a \textit{dense embedding} and that $\scrC$ is a \textit{dense sub-$\infinity$-category}.
\end{defn}
\begin{remark}
\label{remark:dense_cat_localization}
As explained in \cite[ArXiv v4, Remark $5.2.9.4$]{htt}, if $\iota\colon\scrC\to\widehat{\scrC}$ is a functor and $\widehat{\scrC}$ admits all small colimits, then $\iota$ is dense if and only if the induced functor $\Fun(\scrC^{\op},\Ani)\to\widehat{\scrC}$ exhibits $\widehat{\scrC}$ as a localization of the $\infinity$-category of presheaves over $\scrC$. If we drop the assumption that $\scrC$ is small, the same holds if we accept to enlarge our universe and consider presheaves of \textit{large} anima $\Fun(\scrC^{\op},\widehat{\Ani})$.
\end{remark}
\begin{exmp}
\label{exmp:important}
Let $\scrC$ be a Grothendieck site. Then the Yoneda embedding $\yo\colon\scrC\to\Fun(\scrC^{\op},\Ani)$ is always a dense embedding, and composing with the natural sheafification functor $\Fun(\scrC^{\op},\Ani)\to\Shv_{\tau_{\scrC}}(\scrC)$ one obtains a dense functor $\yo\colon\scrC\to\Shv_{\tau_{\scrC}}(\scrC)$. This is a dense embedding if and only if the Grothendieck topology $\tau_{\scrC}$ is sub-canonical.
\end{exmp}
\begin{construction}
\label{constr:kan_extension}
Let $\scrC$, $\cV$ and $\scrD$ be as in \cref{assumption:standing}, let $\scrD^{\mathrm{cat}}\colon\scrC^{\op}\to\catrex$ be the categorification of $\scrD$ as in \cref{def:categorification}, and let $\iota\colon\scrC\to\widehat{\scrC}$ be a dense functor. We have two ways to extend $\scrD^{\mathrm{cat}}$ to an assignment to all $\widehat{\scrC}$.
\begin{enumerate}
 \item We can either consider the right Kan extension $\mathrm{Ran}_{\iota}\scrD(Y)$ of $\scrD$ along (the opposite functor to) $\iota\colon\scrC\subseteq\widehat{\scrC}$, which is given by the assignment
 \[
 \Rani\scrD(Y)\simeq \lim_{\substack{\iota(X)\to Y\\ X\in\scrC}}\scrD(X),
 \]and then consider its categorification
\begin{align*}
   \widehat{\scrC}&\longrightarrow\catrex\\
    Y&\mapsto \LPr_{\mathrm{Ran}_{\iota}\scrD(Y)}.
x\end{align*}
\item Or we can consider the right Kan extension of $\scrD^{\mathrm{cat}}$ along $\iota\colon\scrC\subseteq\widehat{\scrC}$, obtaining the functor $\mathrm{Ran}_{\iota}\scrD^{\mathrm{cat}}$. On an object $Y$ of $\widehat{\scrC}$, it is explicitly described by the assignment
\[
Y\mapsto\lim_{\substack{\iota(X)\to Y\\ X\in\scrC}}\LPr_{\scrD(X)}.
\]
\end{enumerate}
An object in $\LPr_{\Rani\scrD(Y)}$ is a \textit{categorical $\Rani\scrD(Y)$-module}, while an object in $\Rani\scrD^{\mathrm{cat}}(Y)$ can be loosely described as a \textit{sheaf of $\infinity$-categories} described by a rule
\begin{align*}
\scrC\times_{\widehat{\scrC}}\widehat{\scrC}_{/Y}&\longrightarrow\catrex\\
\left\{\alpha\colon \iota(U)\to Y\right\}&\mapsto\scrE_U,
\end{align*}
where $\scrE_U$ is a categorical $\scrD(U)$-module, together with suitable homotopy coherent cocycle conditions. For any object $Y$ of $\widehat{\scrC}$, we have a natural comparison functor
\begin{align}
\label{functor:1_affineness}
\LPr_{\mathrm{Ran}_{\iota}\scrD(Y)}\longrightarrow\mathrm{Ran}_{\iota}\scrD^{\mathrm{cat}}(Y).
\end{align}
induced by the universal property of limits.
\end{construction}
\begin{defn}
\label{def:1_affineness}
Let $\scrC$, $\cV$ and $\scrD$ be as in \cref{assumption:standing}, and let $\iota\colon\scrC\to\widehat{\scrC}$ be a dense functor. We say that an object $Y$ in $\widehat{\scrC}$ is \textit{$1$-affine with respect to $\iota$} if the comparison functor \eqref{functor:1_affineness} is an equivalence.
\end{defn}
\begin{remark}
If $\widehat{\scrC}$ admits all colimits, then \cref{remark:dense_cat_localization} implies that it is a full sub-$\infinity$-category of $\Fun(\scrC^{\op},\Ani)$. In this case, if $Y$ is an object of $\widehat{\scrC}\subseteq\Fun(\scrC^{\op},\Ani)$ we shall simply say that $Y$ is \textit{$1$-affine}. This choice of terminology agrees with the one used of \cite{1affineness}.

\end{remark}
The fundamental question of this paper is the following.
\begin{question}
\label{question:main}
Let $\scrC$, $\cV$, and $\scrD$ be as in \cref{assumption:standing}, and let $\iota\colon\scrC\to\widehat{\scrC}$ be a dense functor. Which objects of $\widehat{\scrC}$ can be proved to be $1$-affine?
\end{question}
\begin{parag}
\label{parag:affineness_explicit}
The quest for $1$-affineness with respect to a dense functor $\iota\colon \scrC\to\widehat{\scrC}$ can be rephrased in more explicit terms. Recall the functor \eqref{functor:1_affineness}: this is the left adjoint in the adjunction
\begin{equation}
    \label{adjunction:1_affineness}
    \adjunction{\widehat{(-)}}{\LPr_{\mathrm{Ran}_{\iota}\scrD(Y)}}{\Rani\scrD^{\mathrm{cat}}(Y)}{\Gamma(Y,-)}.
\end{equation}
The left adjoint is a categorified version of the sheafification of modules on a locally ringed space. Given a categorical $\mathrm{Ran}_{\iota}\scrD(Y)$-module $\scrE$, its sheafification is described by the rule
\[
\left\{\alpha\colon\iota(U)\to Y\right\}\mapsto \scrE\otimes_{\Rani\scrD(Y)}\scrD(U).
\]
This formula makes sense, since $\scrD(Y)$ is a limit of symmetric monoidal $\infinity$-categories (hence it is naturally symmetric monoidal), so $\scrD(U)$ becomes a categorical $\Rani\scrD(Y)$-module via pullback.

On the other hand, the right adjoint $\Gamma(Y,-)$ is a categorified version of the global sections functor. For an object in $\mathrm{Ran}_{\iota}\scrD^{\mathrm{cat}}(Y)$, described by the rule
\[
\left\{\alpha\colon \iota(U)\to Y\right\}\mapsto\scrE_U,
\]
we have an associated diagram
\begin{align*}
\scrC\times_{\widehat{\scrC}}\widehat{\scrC}_Y&\longrightarrow\LPr_{\mathrm{Ran}_{\iota_\ast\scrD(Y)}}\\
\left\{\alpha\colon\iota(U)\to Y\right\}&\mapsto \alpha_\ast\scrE_U.
\end{align*}Then we have
\[
\Gamma((\scrE_U)_{\alpha}, Y)\simeq\lim_{\substack{\alpha\colon \iota(U)\to Y\\ U \in \scrC}}\alpha_\ast \scrE_U.
\]
This discussion makes sense also in the case when $Y=\iota(U)$ lies in $\scrC$. Of course, since $\scrD$ is assumed to be a sheaf for the Grothendieck topology on $\scrC$, the discussion in this case greatly simplifies: for any test space $U$ in $\scrC$ we have equivalences
\[
\mathrm{Ran}_{\iota}\scrD(\iota(U))\simeq \scrD(U)\quad\text{ and }\Rani\scrD^{\mathrm{cat}}(\iota(U))\simeq\scrD^{\mathrm{cat}}(U),
\]
so the adjunction \eqref{adjunction:1_affineness} actually boils down to an adjunction
\begin{equation}
\label{adjunction:1_affineness_affine}
\adjunction{\widehat{(-)}}{\LPr_{\scrD(U)}}{\scrD^{\mathrm{cat}}(U)}{\Gamma(X,-)},
\end{equation}
which is trivially an adjoint equivalence. In general, for all objects $Y$ in $\widehat{\scrC}$ and for any $\scrE$ in $\LPr_{\mathrm{Ran}_{\iota}\scrD(Y)}$ we have a unit functor
\begin{equation}
\label{functor:unit}
\eta\colon \scrE\longrightarrow\Gamma(Y,\widehat{\scrE})
\end{equation}
and for any $(\scrG)_X$ in $\mathrm{Ran}_{\iota}\scrD^{\mathrm{cat}}(Y)$ we have counit functor
\begin{equation}
\label{functor:counit}
\epsilon\colon (\scrG_X)_{\alpha\colon \iota(X)\to Y}\longrightarrow\lp \lim_{\alpha\colon \iota(U)\to Y}\alpha_\ast\scrG\otimes_{\mathrm{Ran}_{\iota}\scrD(Y)}\scrD(X)\rp_X.
\end{equation}
An object is $1$-affine if and only if both functors \eqref{functor:unit} and \eqref{functor:counit} are equivalences.
\end{parag}
In general, it is difficult to characterize $1$-affine objects for an arbitrary choice of $\scrC$ and $\scrD$. For example, it is not clear at all that the categorification $\scrD^{\mathrm{cat}}$ is again a sheaf over $\scrC$, even if $\scrD$ was a sheaf in the first place. In the following, we propose two suitably general frameworks in which $\scrD^{\mathrm{cat}}$ can be proved to be still a sheaf.
\begin{assumption}
\label{assumption:modular_six_functor}
Let $\scrC$, $\scrD$ and $\cV$ be as in \cref{assumption:standing}. Suppose that one of the following two conditions hold.
\begin{assumpenum}
\item\label{assumption:six_functor_case}The sheaf of $\infinity$-categories $\scrD$ is obtained from a strongly monoidal six-functor formalism
\[
(\scrD^\ast,\scrD_!)\colon \mathrm{Corr}(\scrC)^{\otimes}_{E,\mathrm{all}}\longrightarrow\mathrm{Pr}^{\mathrm{L,\otimes}}_{\cV}
\]in the sense of \cite{mann2022p}, and the topology on $\scrC$ is the universal !-able topology (see \cref{def:!-able}).
\item\label{assumption:modular_case} The sheaf $\scrD$ already is the categorification of a sheaf of commutative algebras
\begin{align*}
D\colon\scrC^{\op}\longrightarrow\CAlg(\cV),
\end{align*}
in the sense that we have a natural equivalence of functors
\[
\scrD\simeq\Mod_{D(-)}(\cV).
\]
\end{assumpenum}
\end{assumption}
While the meaning of \cref{assumption:modular_case} is quite self-evident, \cref{assumption:six_functor_case} is a compact way to formulate many different sub-conditions: we shall describe in detail what we mean in \cref{sec:six_functor}. We finish this section by collecting some pieces of terminology that we will use later (essentially to describe non-trivial $1$-affine object).
\begin{defn}
\label{def:affine-object}
Let $\scrC$, $\cV$, and $\scrD$ be as in \cref{assumption:standing}. Let $\iota\colon\scrC\to\widehat{\scrC}$ be a dense functor.
\begin{defenum}
\item\label{def:affine} An object $X$ of $\widehat{\scrC}$ is \textit{$\scrC$-affine} if it belongs to the essential image of $\iota$.
\item\label{def:affine-mor} A morphism $\alpha\colon Y\to X$ is \textit{$\scrC$-affine} if for any $\scrC$-affine object $\iota(U)\to X$ the map $Y\times_X\iota(U)$ is $\scrC$-affine and the projection $Y\times_X\iota(U)\to\iota(U)$ is the image of a map of $\scrC$ under $\iota$.
\item \label{def:atlas}We say that a collection of maps $\left\{\varphi_i\colon W_i\to X\right\}$ is a \textit{$\scrC$-affine atlas} (or \textit{affine $\scrC$-atlas}) if each $W_i$ is $\scrC$-affine and for each $\scrC$-affine $\iota(U)\to X$ the collection of maps $\left\{W_i\times_X\iota(U)\to\iota(U)\right\}$ is the image under $\iota$ of a covering of $\iota(U)$ for the Grothendieck topology on $\scrC$.
\end{defenum}
\end{defn}
\begin{remark}
\label{remark:covering not affine}
Let $\scrC$ be as in \cref{assumption:standing}, and let $\iota\colon\scrC\to\widehat{\scrC}$ be a dense functor in a $\infinity$-category $\widehat{\scrC}$ where colimits are universal. Let $\left\{\iota(W_i)\to X\right\}$ be a $\scrC$-affine atlas of an object $X$ in $\widehat{\scrC}$. Then we can write
\[
X\simeq\colim_{\iota(U)\to X}\iota(U).
\]
Since $\left\{\iota(W_i)\to X\right\}$ is a $\scrC$-affine atlas, for any $\iota(U)\to X$ we let$$\pi_U\colon\coprod_{i=1}^n\iota(W_i)\times_X\iota(U)\to \iota(U)$$be the induced covering of $\iota(U)$, and let $\check{\mathrm{C}}(\pi_U)$ denote its \v{C}ech cover, so as to write
\[
\iota(U)\simeq\colim_{\bDelta^{\op}}\check{\mathrm{C}}(\pi_U).
\]
For each $[n]\in\bDelta$, we have
\[
\check{\mathrm{C}}(\pi_U)^n\simeq\coprod_{\left\{i_1,\cdots,i_n\right\}\subseteq\left\{1,\cdots,m\right\}}\iota(W_{i_1})\times_X\cdots\times_X\iota(W_{i_n})\times_X\iota(U),
\]
and so using the fact that colimits in $\widehat{\scrC}$ are universal we have
\[
\colim_{\iota(U)\to X}\check{\mathrm{C}}(\pi_U)^n\simeq \coprod_{\left\{i_1,\cdots,i_n\right\}\subseteq\left\{1,\cdots,m\right\}}\iota(W_{i_1})\times_X\cdots\times_X\iota(W_{i_n})\simeq\check{\mathrm{C}}(\varphi)^n.
\]
So we have
\[
X\simeq\colim_{\bDelta^{\op}}\colim_{\iota(U)\to X }\check{\mathrm{C}}(\pi_U)\simeq\colim_{\bDelta^{\op}}\check{\mathrm{C}}(\varphi).
\]
In particular, it is still true that if $X$ admits $\scrC$-affine atlas, then it is realized as the colimit of its associated \v{C}ech nerve.
\end{remark}
\subsection{Categorified sheaves under \cref{assumption:six_functor_case}}
\label{sec:six_functor}
Let $\scrC$, $\scrD$ and $\cV$ be as in \cref{assumption:standing}, We first spend some words on what the technical machinery of six-functor formalisms, established in \cite[$\S$A.5]{mann2022p}, specializes to in our setting, and explain clearly what we mean in \cref{assumption:six_functor_case}. 
\begin{fact}\label{fact:six-functor}\
    \begin{factenum}
         \item The pair $(\scrC,E)$ forms a geometric setup in the sense of \cite[Definition A.5.1]{mann2022p}. In particular, $E$ is a class of morphisms which contains all equivalences and is closed under pullbacks.
    \item\label{fact:monoida_structure_corr}Let $\mathrm{Corr}(\scrC)_{E,\mathrm{all}}$ denote the sub-$\infinity$-category of the $\infinity$-category of correspondences of $\scrC$, whose objects are the same as the objects of $\scrC$ but where we allow a span $V\leftarrow U\to W$ to be a morphism in $\mathrm{Corr}(\scrC)_{E,\mathrm{all}}$ only if $U\to W$ belongs to $E$. The $\infinity$-category $\mathrm{Corr}(\scrC)_{E,\mathrm{all}}$ can be regarded as the underlying $\infinity$-category of a $\infinity$-operad
    \[
    \mathrm{Corr}(\scrC)_{E,\mathrm{all}}^{\otimes}\coloneqq\mathrm{Corr}\lp((\scrC^{\op})^{\scriptscriptstyle\coprod})^{\op}\rp_{S^-,\mathrm{all}}
    \]
    using the formalism of coCartesian $\infinity$-operads (see \cite[Proposition $2.4.3.3$]{HA}).
    \item\label{fact:lax monoidal}We have a map of $\infinity$-operads
    \[
    (\scrD^\ast,\scrD_!)\colon \mathrm{Corr}(\scrC)_{E,\mathrm{all}}^{\otimes}\longrightarrow\mathrm{Pr}_{\cV}^{\mathrm{L},\otimes},
    \]
    such that for all morphisms $\varphi\colon U\to V$ we have a symmetric monoidal functor between presentably $\cV$-linear symmetric monoidal $\infinity$-categories $\varphi^\ast\colon\scrD(V)\to\scrD(U)$. If $\psi\colon W\to Y$ lies in $E$, then we have a functor $\psi_!\colon\scrD(W)\to\scrD(Y)$.
    \item Both $\varphi^\ast$ and $\psi_!$ admit right adjoints for abstract reasons. At the same time, the symmetric monoidal structure on each $\scrD(U)$ is closed.
    \end{factenum}
\end{fact}
\begin{remark}
If $\scrC$ admits a terminal object, and thus all limits, then $\scrC^{\op}$ admits all finite coproducts, hence $\smash{\lp(\scrC^{\op})^{\scriptscriptstyle\coprod}\rp^{\op}}$ can be regarded as a twist of the Cartesian monoidal structure on $\scrC$. Concretely, the fiber over a map of finite pointed sets $\alpha\colon J\to I$ of finite sets along the coCartesian fibration of $\infinity$-operads $\smash{\lp(\scrC^{\op})^{\scriptscriptstyle\coprod}\rp^{\op}}\to\mathrm{Fin}_\ast$ consists of maps
\[
X_i\longrightarrow\prod_{j\in\alpha^{-1}(i)}Y_j,
\]
rather than maps in the opposite direction. Still, this defines a symmetric monoidal $\infinity$-category (\cite[Remark 2.4.3.4]{HA}) and so $\mathrm{Corr}(\scrC)_{E,\mathrm{all}}^{\otimes}$ is a symmetric monoidal $\infinity$-category as well (\cite[Proposition 6.1.3]{liu2024enhancedoperationsbasechange}): its symmetric monoidal structure is inherited from such twisted Cartesian monoidal structure of $\scrC$. The $\infinity$-operad $\mathrm{Corr}(\scrC)_{E,\mathrm{all}}^{\otimes}$ in \cref{fact:monoida_structure_corr} is just the underlying $\infinity$-operad of such symmetric monoidal structure. Thus, the map of $\infinity$-operads in \cref{fact:lax monoidal} boils down to a lax monoidal structure on the functor $(\scrD^\ast,\scrD_!)$.
\end{remark}
Let us now explain what the universal !-able Grothendieck topology on $\scrC$ is.
\begin{defn}[{\cite[Definition 4.14]{six-functors-scholze}}]
\label{def:!-able}
Let $\scrC$ be any $\infinity$-category, and let $(\scrD^\ast,\scrD_!)$ a six-functor formalism on $(\scrC,E)$. We define a Grothendieck topology  $\tau^!_{\scrC}$ on $\scrC$ by saying that a collection of maps $\left\{\varphi_i\colon U_i\to V\right\}$ is a \textit{universal !-able cover} if all maps $\varphi_i$ lie in $E$ and the following conditions hold.
\begin{defenum}
    \item \label{def:subcanonical}(The topology is universally sub-canonical.) For any object $Z$ in $\scrC$ and for any map $W\to V$ from an object in $\scrC$, the functor $\Map_{\scrC}(-,Z)$ exhibits $\Map_{\scrC}(W,Z)$ as the limit for the diagram $\restr{\Map_{\scrC}(-,Z)}{U_i\times_VW\to W}.$
    \item\label{def:*-descent}(The topology satisfies universal \textasteriskcentered-descent.) For any map $W\to V$ from an object in $\scrC$, the functor $\scrD$ exhibits $\scrD(W)$ as the limit for the diagram $\restr{\scrD}{U_i\times_VW\to W}$ computed along the functors $\varphi^\ast$.
    \item\label{def:!-descent}(The topology satisfies universal !-descent.) For any map $X\to \yo(V)$ in $\Fun(\scrC^{\op},\Ani)$, the functor $\scrD$ exhibits $\scrD(X)$ as the limit for the diagram $\restr{\scrD}{\yo(U_i)\times_{\yo(V)}X\to X}$ computed along the functors $\varphi^!$.
\end{defenum}
We say that the Grothendieck topology $\tau^!_{\scrC}$ is the \textit{universal !-able topology} on $\scrC$. 
\end{defn}
All together, \cref{assumption:six_functor_case} means that the functor $\scrD$ is formally obtained as the restriction of the functor $(\scrD^\ast,\scrD_!)$ along the natural inclusion of $\infinity$-operads $(\scrC^{\op})^{\scriptscriptstyle\coprod}\subseteq\mathrm{Corr}(\scrC)_{E,\mathrm{all}}^{\otimes}$. Moreover, the functor $\scrD$ is a sheaf for the universal !-able topology on $\scrC.$
\begin{parag}
\label{parag:relative_tensor}
The assumption that $\scrD$ turns finite limits of $\scrC$ into finite colimits of $\CAlg(\LPr_{\cV})$ (\cref{assump:stable}) actually implies that the six-functor formalism
 \[
 (\scrD^\ast,\scrD_!)\colon\mathrm{Corr}(\scrC)^{\otimes}_{E,\mathrm{all}}\longrightarrow\mathrm{Pr}^{\mathrm{L},\otimes}_{\cV}
 \]
 is \textit{locally strongly monoidal}, in the following sense. For any test space $V$, the $\infinity$-category $\scrC_{/V}$ admits all finite limits (all finite limits but the terminal object are inherited by limits in $\scrC$, while the terminal object is realized by $V$), so the procedure described in \cref{fact:monoida_structure_corr} equips $\mathrm{Corr}(\scrC_{/V})_{E,\mathrm{all}}$ with a symmetric monoidal structure. Then the induced functor
 \[
\restr{(\scrD^\ast,\scrD_!)}{V}\colon\mathrm{Corr}(\scrC_{/V})^{\otimes}_{E,\mathrm{all}}\longrightarrow\mathrm{Pr}^{\mathrm{L},\otimes}_{\scrD(V)}
 \]
 is a symmetric monoidal functor of $\infinity$-categories. This is true because the product in $\scrC_{/V}$ of two objects $U\to V$ and $W\to V$ (which is the pullback $U\times_VW$ in $\scrC$) is sent to the pushout of the diagram $\scrD(U)\leftarrow\scrD(V)\to\scrD(W)$ in $\CAlg(\LPr_{\cV})$, which is precisely the relative tensor product $\scrD(U)\otimes_{\scrD(V)}\scrD(W)$.
\end{parag}
We now explain how \cref{assumption:six_functor_case} implies that we can deduce a strong dualizability condition at the categorified level, whenever a map of test spaes $\varphi\colon U\to V$ is a universal !-able cover.
\begin{construction}
Let $V$ be any test space in $\scrC$, let $\scrC_{/V}$ denote the slice $\infinity$-category over $V$. Let $\scrC_{/V}^S$ denote the sub-$\infinity$-category of $\scrC_{/V}$ described as follows.
\begin{enumerate}
\item Objects are maps $U\to V$ which belong to $E$.
\item A morphism between two arrows $\varphi\colon U\to V$ and $\psi\colon W\to V$ is a map in $\scrC_{/V}$.
\end{enumerate}
Notice that \cite[Lemma $2.1.5$.iii]{heyer6ff2024} implies that whenever $U\to V$ and $W\to V$ belong to $E$, then any morphism $U\to W$ commuting with the maps to $V$ is forced to lie in $E$ as well. This is easily seen to imply the following.
\begin{proposition}
\label{prop:finite_limits_local}
Let $\scrC$ and $E$ be as in \cref{assumption:six_functor_case}. Then $\scrC^{E}_{/V}\subseteq\scrC_{/V}$ is a full sub-$\infinity$-category closed under finite limits.
\end{proposition}
As a consequence of \cref{prop:finite_limits_local}, we have a strongly monoidal inclusion
\[
\mathrm{Corr}(\scrC^{E}_{/V})^{\otimes}\coloneq\mathrm{Corr}(\scrC^{E}_{/V})^{\otimes}_{\mathrm{all,all}}\subseteq\mathrm{Corr}(\scrC_{/V})^{\otimes}_{E,\mathrm{all}}.
\]
In particular, the restriction
\begin{align}
\label{functor:local_six_functor}
\restr{(\scrD^\ast,\scrD_!)^{E}}{V}\colon \mathrm{Corr}(\scrC^{E}_{/V})^{\otimes}\subseteq\mathrm{Corr}(\scrC_{/V})^{\otimes}_{E,\mathrm{all}}\longrightarrow\mathrm{Pr}^{\mathrm{L},\otimes}_{\scrD(V)}
\end{align}
remains strongly monoidal. 
\end{construction}

We want to show how \cref{prop:finite_limits_local} implies that  $\scrD(U)$ is $\scrD(V)$-dualizable, when $\varphi\colon U\to V$ lies in $E$. In order to do so, recall the following fundamental property of (higher) $\infinity$-categories of correspondences.
\begin{proposition}[{\cite[Theorem $1.4$]{iteratedspans}}]
\label{prop:iterated-spans}
If $\scrC$ is a $\infinity$-category admitting all finite limits, then for all $n\geqslant 1$ every object in the $(\infinity,n)$-category $\mathrm{Corr}_n(\scrC)$ of iterated spans is fully dualizable.
\end{proposition}
\begin{remark}
\label{remark:iterated-spans-2}
Actually, the proof of \cref{prop:iterated-spans} in \cite{iteratedspans} shows something more: each object is self dual. The evaluation and coevaluation correspondences
\[
\begin{tikzpicture}[scale=0.75,baseline=(current  bounding  box.center)]
\node (a) at (0,2){$X$};
\node (b) at (-2,0){$\left\{\ast\right\}$};
\node (c) at (2,0){$X\times X$};
\draw[->,font=\scriptsize](a) to[bend right] node[left]{}(b);
\draw[->,font=\scriptsize](a) to[bend left] node[right]{$\Delta$}(c);
\end{tikzpicture}\quad\text{ and }\quad\begin{tikzpicture}[scale=0.75,baseline=(current  bounding  box.center)]
\node (a) at (0,2){$X$};
\node (c) at (2,0){$\left\{\ast\right\}$};
\node (b) at (-2,0){$X\times X$};
\draw[->,font=\scriptsize](a) to[bend left] node[right]{}(c);
\draw[->,font=\scriptsize](a) to[bend right] node[left]{$\Delta$}(b);
\end{tikzpicture}
\]
are one the transpose of the other.
\end{remark}
\begin{corollary}[{\cite[Lecture 17, around 1:19:00]{Clausen_Scholze_lectures}}]
\label{cor:dualizability_six_functor}
Let $\scrC$, $\scrD$, $\cV$ be as in \cref{assumption:six_functor_case}. Then for any map of test spaces $\varphi\colon U\to V$ lying in $E$ the pullback $\varphi^\ast\colon\scrD(V)\to\scrD(U)$ turns $\scrD(U)$ into a dualizable $\scrD(V)$-module. In fact, this $\scrD(V)$-module is self-dual, and for any map of test spaces $\alpha\colon U\to W$ over $V$ which belongs to $E$ the dual of the functor $\alpha^\ast\colon\scrD(W)\to\scrD(U)$ corresponds to $\alpha_!\colon\scrD(U)\to\scrD(W)$ under the equivalences
\[
\scrD(U)\simeq\scrD(U)^{\vee}\quad\text{ and }\quad\scrD(W)\simeq\scrD(W)^{\vee}.
\]
\end{corollary}
\begin{proof}
Just notice that, being symmetric monoidal, the functor \eqref{functor:local_six_functor} preserves duals. Then everything is a formal consequence of \cref{prop:iterated-spans} and \cref{remark:iterated-spans-2}.
\end{proof}
We are now ready to prove that the functor \eqref{functor:categorification} is a sheaf for the Grothendieck topology on $\scrC$.
\begin{proposition}
\label{prop:six_functor_case}
Let $\scrC$, $\scrD$ and $\cV$ be as in \cref{assumption:standing}. Assume \cref{assumption:six_functor_case} to hold. Then the functor \eqref{functor:categorification} is a sheaf for the universal !-able Grothendieck topology on $\scrC$.
\end{proposition}
\begin{proof}
The assumption on our Grothendieck topology $\tau_{\scrC}$ allows us to reduce the above claim to the following (see also \cite[Proposition A.3.3.1]{sag}).
\begin{enumerate}
    \item The functor $\scrD^{\mathrm{cat}}$ preserves finite products.
    \item For any map of test spaces $\varphi\colon V\to U$ which is a universal !-able cover, the functor
    \[
    \scrD^{\mathrm{cat}}(U)\longrightarrow\lim_{[n]\in\bDelta}\lp\scrD^{\mathrm{cat}}(V)\stack{3}\scrD^{\mathrm{cat}}(V\times_UV)\stack{5}\scrD^{\mathrm{cat}}(V\times_UV\times_UV)\stack{7}\cdots\rp
    \]
    is an equivalence.
\end{enumerate}
The first claim follows from an easy check of the Barr--Beck--Lurie's monadicity theorem, using the fact that for \textit{any} collection of presentably symmetric monoidal $\infinity$-categories $\cV_i$ one has
\[
\prod_{i\in I}\LPr_{\cV_i}\simeq\LPr_{\prod_i\cV_i},
\]
see for example \cite[Proposition 3.2.32]{PPS_1}. In order to prove the second claim, we argue as follows: since all objects involved belong to $\scrC$, they are clearly $1$-affine with respect to $\iota$. So the claim reduces to proving that the functor
\[
\LPr_{\scrD(U)}\longrightarrow\lim_{[n]\in\bDelta}\LPr_{\scrD\lp V^{\times_U(n+1)}\rp}
\]
is an equivalence. Once again, this is a left adjoint functor which sends a categorical $\scrD(U)$-module $\scrE$ to the collection of $\scrD\lp V^{\times_U(n+1)}\rp$-modules $\scrE\otimes_{\scrD(U)}\scrD\lp V^{\times_U(n+1)}\rp$. Its right adjoint simply takes a compatible collection of $\scrD\lp V^{\times_U(n+1)}\rp$-modules $\scrE_n$ and computes the totalization of the diagram $\varphi_\ast\scrE_n$ in $\LPr_{\scrD(U)}$. The fact that the unit functor
\[
\scrE\longrightarrow\lim_{[n]\in\bDelta}\scrE\otimes_{\scrD(U)}\scrD\lp V^{\times_U(n+1)}\rp
\]
is an equivalence is essentially due to the fact that $\scrD$ is a sheaf for the universal !-able topology. Indeed, we argue as follows: we have an equivalence
\[
\scrD(U)\simeq{\limast_{[n]\in\bDelta}}\scrD\lp V^{\times_U(n+1)}\rp\simeq\limsh_{[n]\in\bDelta}\scrD\lp V^{\times_U(n+1)}\rp,
\]
where the superscripts $\ast$ and $!$ denote the functors with respect to which the limits are computed. Since each $\scrD\lp V^{\times_U(n+1)}\rp$ is self-dual over $\scrD(U)$ (\cref{cor:dualizability_six_functor}), we have that
\[
\scrE\otimes_{\scrD(U)}\scrD\lp V^{\times_U(n+1)}\rp \simeq\FunL_{\scrD(U)}\lp\scrD\lp V^{\times_U(n+1)}\rp, \scrE\rp
\]
so we have
\begin{align*}
{\limast_{[n]\in\bDelta}}\scrE\otimes_{\scrD(U)}\scrD\lp V^{\times_U(n+1)}\rp&\simeq{\limast_{[n]\in\bDelta}}\FunL_{\scrD(U)}\lp\scrD\lp V^{\times_U(n+1)}\rp,\scrE\rp\\&\simeq\FunL_{\scrD(U)}\lp{\colimastd_{[n]\in\bDelta^{\op}}}\scrD\lp V^{\times_U(n+1)}\rp,\scrE\rp\\
&\simeq\FunL_{\scrD(U)}{\lp\colimsh_{[n]\in\bDelta^{\op}}\scrD\lp V^{\times_U(n+1)}\rp,\scrE\rp}\\&\simeq\FunL_{\scrD(U)}\lp\limsh_{[n]\in\bDelta}\scrD\lp V^{\times_U(n+1)}\rp,\scrE\rp
\\&\simeq \FunL_{\scrD(U)}{\lp\scrD(U),\scrE\rp}\simeq\scrE.
\end{align*}

The fact that for an object $\left\{\scrE_n\right\}_{[n]}$ in $\lim_{[n]}\LPr_{\scrD\lp V^{\times_U(n+1)}\rp }$ the counit
\[
\left\{\lp\lim_{[m]\in\bDelta}\scrE_m\rp\otimes_{\scrD(U)}\scrD\lp V^{\times_U(n+1)}\rp \right\}_{[n]\in\bDelta}\longrightarrow\left\{\scrE_n\right\}_{[n]\in\bDelta}
\]
is an equivalence is a consequence of the fact that each $\scrD\lp V^{\times_U(n+1)}\rp $ is dualizable over $\scrD(U)$. Indeed, it is sufficient that the above functor is an equivalence at the $n=0$ stage, since for $n\geqslant1$ the functor at stage $n$ is obtained by tensoring with $\scrD\lp V^{\times_U(n+1)}\rp $. In this case, we have
\[
\lim_{[m]\in\bDelta}\scrE_0\otimes_{\scrD(U)}\scrD\lp V^{\times _U(m+1)}\rp \longrightarrow\scrE_0.
\]
But the cosimplicial object $\scrE_0\otimes_{\scrD(U)}\scrD\lp V^{\times_U(m+1)}\rp$ is cosplit via the diagonal, with coaugmentation given precisely by $\scrE_0\to\scrE_0\otimes_{\scrD(U)}\scrD(V)$. Therefore, $\scrE_0$ is the (universal) limit of such coaugmented cosimplicial diagram, and this finishes the proof.
\end{proof}
In particular from the above proposition we can deduce the following corollary 
\begin{corollary}
    Every object in $\scrC$ is 1-affine.
\end{corollary}

Let $\scrC$, $\cV$ and $\scrD$ be as in \cref{assumption:standing}, and suppose that \cref{assumption:six_functor_case} holds. In order to study the problem of $1$-affineness relatively to a dense functor $\iota\colon\scrC\to\widehat{\scrC}$, we shall assume the following.
\begin{assumption}
\label{assumption:reasonable}\
\begin{assumpenum}
\item \label{assumption:dense_embedding}The dense functor $\iota\colon\scrC\to\widehat{\scrC}$ is fully faithful.
\item \label{assumption:limit-preserving}The $\infinity$-category $\widehat{\scrC}$ admits all finite limits, and the dense functor $\iota\colon\scrC\to\widehat{\scrC}$ preserves pullbacks of $\scrC$ along all maps in $E$.
\item\label{assumption:reasonable-closure}If the pullback of a morphism $\alpha\colon X\to \iota(U)$ along a map $\iota(W)\to\iota(U)$ lies in $E$, then $X\simeq\iota(V)$ and $\alpha$ is the image under $\iota$ of a morphism of $E$.
\end{assumpenum}
\end{assumption}
We start by extending the six-functor formalism $(\scrD^\ast,\scrD_!)$ from $\scrC$ to $\widehat{\scrC}$ along $\iota$, following a recipe due to Heyer and Mann \cite{heyer6ff2024}.
\begin{proposition}[{\cite[Proposition 3.4.2]{heyer6ff2024}}]
\label{prop:extension}
Let $\scrC$, $\scrD$ and $\cV$ be as in \cref{assumption:six_functor_case}, and let $\iota\colon\scrC\to\widehat{\scrC}$ be a dense functor as in \cref{assumption:reasonable}. Let $\overline{E}$ be the class of morphisms $\varphi\colon X\to Y$ of $\widehat{\scrC}$ satisfying the following: for all test spaces $U$ in $\scrC$ and for all maps $\iota(U)\to Y$, the pullback $\iota(U)\times_YX$ is equivalent to $\iota(V)$ for some object $V$ in $\scrC$ and $\iota(U)\times_YX\to \iota(U)$ is the image of a morphism of $E$ under $\iota$. Then we have the following.
\begin{propenum}
    \item\label{prop:extension1} The pair $(\widehat{\scrC},\smash{\overline{S}}')$ is a geometric setup.
    \item\label{prop:extension2} The functor $\iota\colon\scrC\to\widehat{\scrC}$ is a morphism of geometric setups (\cite[$\S2.1$]{heyer6ff2024}).
    \item\label{prop:extension3} The six-functor formalism $(\scrD^\ast,\scrD_!)$ extends uniquely to a six-functor formalism $(\smash{\overline{\scrD}}^\ast,\smash{\overline{\scrD}}_!)$ on $\widehat{\scrC}$, via a left Kan extension
    \[
    \overline{\scrD}(X)\coloneqq\lim_{\iota(U)\to X}\scrD(U).
    \]
\end{propenum}
\end{proposition}

\begin{proof}
The first two claims can clearly be proved in the same way as in \cite[Proposition 3.4.2]{heyer6ff2024}, while the third claim is just \cite[Proposition A.5.16]{mann2022p}.
\end{proof}
\begin{proposition}[{\cite[Theorem 3.4.11]{heyer6ff2024}}]
\label{prop:extension strong}
Let $\scrC$, $\scrD$ and $\cV$ be as in \cref{assumption:six_functor_case}, and let $\iota\colon\scrC\to\widehat{\scrC}$ be a dense functor satisfying the assumptions of \cref{prop:extension}. Let $(\widehat{\scrC},\overline{E})$ denote the geometric setup obtained from the geometric $(\scrC,E)$ as in \cref{prop:extension}. Then there exists a collection of morphisms $\widehat{E}$ in $\widehat{\scrC}$ satisfying the following.
\begin{propenum}
\item\label{prop:extension strong 1}The functor $\iota$ defines a morphism of geometric setups $(\scrC,E)\to(\widehat{\scrC},\widehat{E})$, and the six-functor formalism $(\overline{\scrD}^\ast,\overline{\scrD}_!)$ on the source extends uniquely to a six-functor formalism $(\widehat{\scrD}^\ast,\widehat{\scrD}_!)$ on the target.
\item\label{prop:extension strong 2}The class $\widehat{E}$ is \emph{\textasteriskcentered-local on the target}: if a morphism $\alpha\colon X\to Y$ is such that for every morphism $\iota(U)\to Y$ the pullback $\iota(U)\times_YX\to \iota(U)$ lies in $\widehat{E}$, then $\alpha$ is in $\widehat{E}$.
\item\label{prop:extension strong 3}The class $\widehat{E}$ is \emph{!-local on both source and target}: if a morphism $\alpha\colon X\to Y$ is such that there exists either a small universal !-cover $\left\{\iota(U_i)\to Y\right\}$ such that all morphisms $\iota(U_i)\times_YX\to\iota(U_i)$ are in $\widehat{E}$, or a small universal !-cover $\left\{\iota(V_j)\to X\right\}$ such that all compositions $\iota(V_i)\to Y$ are in $\widehat{E}$, then $\alpha$ lies in $\widehat{E}$.
\item\label{prop:extension strong 4}The class $\widehat{E}$ is \emph{tame}: every morphism $Y\to \iota(U)$ which lies in $\widehat{E}$ is !-locally on the source in $E$.
\end{propenum}
Moreover, such choice can be chosen to be minimal with respect to the inclusion order.
\end{proposition}
\begin{remark}
The recipe in \cref{prop:extension strong} extends in a non-trivial way the class of morphisms $\alpha$ for which the functors $\alpha_!$ and $\alpha^!$ are defined, but it does not alter the stable and presentable $\infinity$-category $\overline{\scrD}(X)$ associated to an object $X$ in $\widehat{\scrC}$. In particular, the stable and presentable $\infinity$-category $\widehat{\scrD}(X)$ is once again $\overline{\scrD}(X)$, which in turn is just $\Rani\scrD(X)$. All in all, \cref{prop:extension strong} tells us that the left Kan extension of one of the operations the six-functor formalism is itself part of a six-functor formalism, and there exists a clever way to add maps to the collection $\overline{E}$ in order to keep all the good properties of the six-functor formalism. 
\end{remark}
\begin{proof}
While the statement in \cite[Theorem 3.4.11]{heyer6ff2024} only considers the case where $\widehat{\scrC}$ is the $\infinity$-category of sheaves for the universal !-able topology, the proof applies in our generality as well. Indeed, the only point where the proof explicitly relies on the fact that $\widehat{\scrC}=\mathrm{Shv}_{\tau^!_{\scrC}}(\scrC)$ is to ensure that the class $\mathcal{S}$ of all possible collections of morphisms $\widehat{E}$ which satisfy \cref{prop:extension strong 1} and \cref{prop:extension strong 4} is not empty. However, as already observed in \cref{prop:extension}, this is true also in our case since $\mathcal{S}$ contains at least $\overline{E}$.
\end{proof}
\begin{remark}
If $\iota$ satisfies all conditions in \cref{assumption:reasonable} except for the fully faithfulness, it should still be feasible to adapt the proofs of \cite[Proposition 3.4.2 and Theorem 3.4.11]{heyer6ff2024} obtaining a six-functor formalism $(\smash{\widehat{\scrD}}^\ast,\smash{\widehat{\scrD}}_!)$ on $\widehat{\scrC}$. Although this will not be an extension of the six-functor formalism $(\scrD^\ast,\scrD_!)$ anymore, it still is the universal six-functor formalism making the diagram of $\infinity$-operads
\[
\begin{tikzpicture}[scale=0.75]
\node (a) at (0,2){$\mathrm{Corr}(\scrC)^{\otimes}_{E,\mathrm{all}}$};
\node (b) at (0,0){$\mathrm{Corr}(\widehat{\scrC})^{\otimes}_{\overline{E},\mathrm{all}}$};
\node (c) at (4,2){$\mathrm{Pr}_{\cV}^{\mathrm{L},\otimes}$};
\draw[->,font=\scriptsize](a) to node[left]{$\iota$}(b);
\draw[->,font=\scriptsize] (a) to node[above]{$(\scrD^\ast,\scrD_!)$}(c);
\draw[->,dotted,font=\scriptsize](b) to[bend right] node[right]{$(\widehat{\scrD}^\ast,\widehat{\scrD}_!)$}(c);
\end{tikzpicture}
\]
commute. However, since the Grothendieck topology $\tau^!_{\scrC}$ is required to be sub-canonical, and since we want $\widehat{\scrD}(\iota(U))$ to be the same as $\scrD(U)$ for any test space $U$, we do not investigate the matter further.
\end{remark}
\begin{remark}
With the general recipe of \cref{prop:extension strong}, the six-functor formalism can be extended to the geometric setup $(\widehat{\scrC},\widehat{E})$. So it makes sense to define universal !-able covers also in $\widehat{\scrC}$.
\end{remark}
The following proposition is close in spirit to the results studied in \cite{kunneth6ff} concerning categorical K\"unneth formulas for analytic stacks, and especially to \cite[Theorem 3.26]{kunneth6ff}. The proof is almost the same as the one in \textit{loc. cit.}, but we still provide a proof for the reader's convenience.
\begin{proposition}\label{prop: Kunneth per !-covers}
Let $\scrC$, $\cV$, $\scrD$ be as in \cref{assumption:six_functor_case}, and let $\iota\colon\scrC\to\widehat{\scrC}$ be a dense embedding as in \cref{assumption:reasonable}. Let $Y$ be an object in $\widehat{\scrC}$, and suppose that $Y$ admits a universal !-able cover of the form $\left\{\iota(W_i)\to Y\right\}_{i\in I}$ such that all $n$-fold fiber products $\iota(W_{i_1})\times_Y\cdots\times_Y\iota(W_{i_n})$ are again $\scrC$-affine. Then for every couple of indices $\alpha,\beta\in I$ the natural functor
    \begin{equation*}
        \mathrm{Ran}_{\iota}\scrD(\iota(W_{\alpha})\times_{Y}\iota(W_{\beta}))\longrightarrow \mathrm{Ran}_{\iota}\scrD(W_{\alpha})\otimes_{\mathrm{Ran}_{\iota}\scrD(Y)}\scrD(W_{\beta})
    \end{equation*}
    is an equivalence.
\end{proposition}
\begin{proof}
Let us fix some notations first. For $\alpha$ and $\beta$ as in the statement, write
\begin{align*}
W_{i_1\cdots i_n}\quad&\text{for}\quad\coprod_{\left\{i_1,\dots,i_n\right\}\subseteq I} W_{i_1}\times_Y \cdots \times_Y W_{i_n},\\
W_{\alpha_1\cdots \alpha_n}\quad&\text{for}\quad\coprod_{\left\{i_1,\dots,i_n\right\}\subseteq I} W_{\alpha}\times_Y\lp W_{i_1}\times_Y\cdots\times_Y W_{i_n}\rp,\\
W_{\beta_1\cdots \beta_n}\quad&\text{for}\quad\coprod_{\left\{i_1,\dots,i_n\right\}\subseteq I} W_{\beta}\times_Y\lp W_{i_1}\times_Y\cdots\times_Y W_{i_n}\rp.
\end{align*}
We have
\begin{align*}
\mathrm{Ran}_{\iota}\scrD(\iota(W_{\alpha}))\otimes_{\mathrm{Ran}_{\iota}\scrD(Y)}\Rani\scrD(\iota(W_{\beta}))&\simeq\scrD(W_{\alpha})\otimes_{\Rani\scrD(Y)}\scrD(W_{\beta})\\
&\simeq\colim_{[n]\in\bDelta^{\op}}\scrD(W_{\alpha})\otimes_{\cV}\Rani\scrD(Y)^{\otimes n}\otimes_{\cV}\scrD(W_{\beta}).
\end{align*}
Using descent along the universal !-able covers $\left\{W_i\to Y\right\}$, $\left\{W_{\alpha}\times_YW_i\to W_{\alpha}\right\}$ and $\left\{W_{\beta}\times_YW_i\to W_{\beta}\right\}$ we deduce
\begin{align*}
\scrD(W_{\alpha})&\simeq\limast_{[m]\in\bDelta}\scrD(W_{\alpha_1\cdots\alpha_m})\simeq\limsh_{[m]\in\bDelta}\scrD(W_{\alpha_1\cdots\alpha_m})\simeq\colimsh_{[m]\in\bDelta^{\op}}\scrD(W_{\alpha_1\cdots\alpha_m}),\\
\Rani\scrD(Y)&\simeq\limast_{[m]\in\bDelta}\scrD(W_{i_1\cdots i_m})\simeq\limsh_{[m]\in\bDelta}\scrD(W_{i_1\cdots i_m})\simeq\colimsh_{[m]\in\bDelta^{\op}}\scrD(W_{i_1\cdots i_m}),\\
\scrD(W_{\beta})&\simeq\limast_{[m]\in\bDelta}\scrD(W_{\beta_1\cdots\beta_m})\simeq\limsh_{[m]\in\bDelta}\scrD(W_{\beta_1\cdots\beta_m})\simeq\colimsh_{[m]\in\bDelta^{\op}}\scrD(W_{\beta_1\cdots\beta_m}).
\end{align*}
Moreover, since $\bDelta^{\op}$ is sifted, $\bDelta^{\op}\to\bDelta^{\op,\times n}$ is a final functor for all integers $n$, thus we have
\[
\Rani\scrD(Y)^{\otimes n}\simeq \lp\colimsh_{[m]\in\bDelta^{\op}}\scrD(W_{i_1\cdots i_m})\rp^{\otimes n}\simeq\colimsh_{[m]\in\bDelta^{\op}}(\scrD(W_{i_1\cdots i_m})^{\otimes n}).
\]
All things considered, we can write
\begin{align*}
\colim_{[n]\in\bDelta^{\op}}\scrD(W_{\alpha})\otimes_{\cV}\Rani\scrD(Y)^{\otimes n}\otimes_{\cV}\scrD(W_{\beta})&\simeq \colimsh_{\bDelta^{\op,\times 4}}\scrD(W_{\alpha_1\cdots\alpha_{m_1}})\otimes_{\cV}\scrD(W_{i_1\cdots i_{m_2}})^{\otimes n}\otimes_{\cV}\scrD(W_{\beta_1\cdots\beta_{m_3}})\\
&\simeq\colimsh_{\bDelta^{\op,\times 2}}\scrD(W_{\alpha_1\cdots\alpha_{m}})\otimes_{\cV}\scrD(W_{i_1\cdots i_{m}})^{\otimes n}\otimes_{\cV}\scrD(W_{\beta_1\cdots\beta_{m}})\\
&\simeq\colimsh_{[m]\in\bDelta^{\op}}\scrD(W_{\alpha_1\cdots\alpha_m})\otimes_{\scrD(W_{i_1\cdots i_n})}\scrD(W_{\beta_1\cdots\beta_m})\\
&\simeq\colimsh_{[m]\in\bDelta^{\op}}\scrD\lp W_{\alpha_1\cdots\alpha_m}\times_{W_{i_1\cdots i_n}}W_{\beta_1\cdots\beta_m}\rp.
\end{align*}
Notice that in the second equivalence we used once again that $\bDelta^{\op}$ is sifted, while in the last equivalence we used the fact that the six-functor formalism $(\scrD^\ast,\scrD_!)$ is assumed to be strongly monoidal, and thus $\Rani\scrD$ sends fiber products of objects in $\scrC$ to relative tensor products of presentably symmetric monoidal $\infinity$-categories (see discussion in $\S$\ref{parag:relative_tensor}). We now observe that
\begin{align*}
W_{\alpha_1\cdots\alpha_m}\times_{W_{i_1\dots i_m}}W_{\beta_1\cdots\beta_m}&\coloneqq W_{\alpha}\times_Y W_{i_1\dots i_n}\times_{W_{i_1\dots i_m}}W_{\beta}\times_Y{W_{i_1\dots i_m}}\\
&\simeq W_{\alpha}\times_Y W_{\beta}\times_Y{W_{i_1\dots i_m}}
\end{align*}
is just the $m$-th step in the \v{C}ech nerve of the universal !-able cover
\[
\left\{W_{\alpha}\times_Y W_{\beta}\times_Y{W_{i_1\dots i_m}}\right\}\longrightarrow W_{\alpha}\times_Y W_{\beta}, 
\]
so using !-descent we conclude as desired.
\end{proof}
\cref{prop: Kunneth per !-covers} is the main stepping stone in proving the following theorem.
\begin{theorem}
\label{thm:main 2-descent}
Let $\scrC$, $\cV$, $\scrD$ be as in \cref{assumption:six_functor_case}, and let $\iota\colon\scrC\to\widehat{\scrC}$ be a dense embedding as in \cref{assumption:reasonable}. Let $Y$ be an object in $\widehat{\scrC}$, and suppose that $Y$ admits a universal !-able cover of the form $\left\{\iota(W_i)\to Y\right\}_{i\in I}$ such that all $n$-fold fiber products $\iota(W_{i_1})\times_Y\cdots\times_Y\iota(W_{i_n})$ are again $\scrC$-affine. Then $Y$ is $1$-affine relatively to $\iota$.
\end{theorem}
\begin{proof}
We endow the $\infinity$-category $\widehat{\scrC}$ of the universal !-able topology for the six-functor formalism established in \cref{prop:extension strong}. As explained in \cite[Theorem 3.4.11]{heyer6ff2024}, the functor $\Rani\scrD$ is again a sheaf for the universal !-able topology on $\widehat{\scrC}$, since $\scrD$ was a sheaf in the first place (\cref{prop:six_functor_case}).

In order to prove the $1$-affineness of $Y$, we simply need to prove that the functor $\LPr_{\Rani\scrD(-)}$ satisfies descent along the cover $\left\{\iota(W_i)\to Y\right\}$. Once we know this, since all $\scrC$-affine objects are trivially $1$-affine relatively to $\iota$, we would have a chain of equivalences
\[
\Rani\scrD^{\mathrm{cat}}(Y)\coloneqq\lim_{[n]\in\bDelta^{\op}}\prod_{\left\{i_1,\cdots,i_n\right\}\subseteq I}\scrD^{\mathrm{cat}}(W_{i_1}\times_Y\cdots\times_Y W_{i_n})\simeq\lim_{[n]\in\bDelta^{\op}}\prod_{\left\{i_1,\cdots,i_n\right\}\subseteq I}\LPr_{\scrD(W_{i_1}\times_Y\cdots\times_Y W_{i_n})}\simeq \LPr_{\Rani\scrD(Y)}.
\]
Notice that, while it is true that we have at our disposal a six-functor formalism $(\widehat{\scrD}^\ast,\widehat{\scrD}_!)$ on all $\widehat{\scrC}$ which extends the six-functor formalism $(\scrD^\ast,\scrD_!)$ over $\scrC$, it is not true that the six-functor formalism $(\widehat{\scrD}^\ast,\widehat{\scrD}_!)$ is locally strongly monoidal (in the sense of Paragraph \ref{parag:relative_tensor}) anymore. It follows that, in general, for a map $Z\to Y$ lying in $\widehat{E}$ the categorical $\Rani\scrD(Y)$-module $\Rani\scrD(Z)$ will not be self-dual anymore. But \cref{prop: Kunneth per !-covers} is precisely the ingredient that allows us to deduce that each $\scrD(W_i)$ is self-dual as a categorical $\Rani\scrD(Y)$-module. Indeed, let $\langle \iota(W_i)\rangle_{/Y}$ be the smallest sub-$\infinity$-category of $\widehat{\scrC}_{/Y}$ containing all maps $\iota(W_i)\to Y$ and admitting finite limits: this means that we are adding the identity $Y=Y$ and all maps from finitary fiber products $\iota(W_{i_1})\times_Y\cdots\times_Y\iota(W_{i_n})\to Y$. Then \cref{prop: Kunneth per !-covers} means that the composition
\[
\mathrm{Corr}(\langle \iota(W_i)\rangle_{/Y})^{\otimes}\longrightarrow\mathrm{Corr}(\widehat{\scrC}_{/V})^{\otimes}_{\widehat{E},\mathrm{all}}\overset{(\widehat{\scrD}^\ast,\widehat{\scrD}_!)}{\longrightarrow}\mathrm{Pr}^{\mathrm{L},\otimes}_{\Rani\scrD(Y)}
\]
is now strongly monoidal, hence it preserves duals, and so we can safely apply \cref{prop:iterated-spans} to establish that each $\scrD(W_i)$ is self-dual as a categorical $\Rani\scrD(Y)$-module. At this point, we can follow the same strategy of proof of \cref{prop:six_functor_case} to deduce that the natural functor
\[
\LPr_{\Rani\scrD(Y)}\longrightarrow\lim_{\iota(W_i)\to Y}\LPr_{\scrD(W_i)}
\]
is an equivalence, and so we conclude. 
\end{proof}
\subsection{Categorified sheaves under \cref{assumption:modular_case}}
\label{sec:modular_case}
Let $\scrC$, $\scrD$ and $\cV$ be as in \cref{assumption:standing}. Assume moreover that \cref{assumption:modular_case} holds: this means that for any morphism of test spaces $\varphi\colon U\to V$ the adjunction $\varphi^\ast\dashv\varphi_\ast$ can be interpreted as a free--forgetful adjunction of the form
\[
\adjunction{-\otimes_{D(V)}D(U)}{\Mod_{D(V)}(\cV)}{\Mod_{D(U)}(\cV)}{\oblv_{D(U)}}.
\]
In particular, the right adjoint commutes with \textit{all} limits and colimits and it is lax monoidal (\cite[Proposition A]{laxmonoidal}). It follows that $\oblv_{D(U)}D(U)$ is still a commutative algebra inside $\Mod_{D(V)}(\cV)$ (equivalently, it is a commutative $D(V)$-algebra). By abuse of notation, we shall still denote it as $D(U)$.
\begin{defn}[{\cite[Definition $3.18$]{mathew2016galois}}]
\label{def: descendable}
Let $\scrC$, $\scrD$ and $\cV$ be as in \cref{assumption:standing}. We say that a morphism of test spaces $\varphi\colon U\to V$ is \textit{descendable} if the thick $\otimes$-ideal generated by $f_\ast\boldone_{\scrD(U)}$ in $\Mod_{D(V)}$ coincides with the whole $\Mod_{D(V)}$.
\end{defn}
The following is a straight-forward generalization of \cite[Proposition 3.45]{mathew2016galois}.
\begin{proposition}
\label{prop:modular_case}
Let $\scrC$, $\scrD$ and $\cV$ be as in \cref{assumption:standing}, and assume that \cref{assumption:modular_case} holds. Assume moreover that any morphism of test spaces $\varphi\colon U\to V$ belonging to $S$ is descendable. Then the functor \eqref{functor:categorification} is a sheaf for the Grothendieck topology on $\scrC$.
\end{proposition}
\begin{proof}
This is a formal consequence of descent for $\infinity$-categories of modules, as showed in \cite[Proposition 3.45]{mathew2016galois}. There, this assertion is stated in the particular setting when $\cV=\Sp$ is the $\infinity$-category of spectra. However, its proof only relies on the notion of descendability (that makes sense also in our more general setting), on Barr--Beck--Lurie's (co)monadicity theorem (which holds for any comonad over any $\infinity$-category), and on the fact that given a map of test spaces $\varphi\colon U\to V$ which is descendable the functor $\scrD(-)$ satisfies descent along $\varphi$ (which is a piece of datum already in \cref{assumption:standing}).

In particular, that proof can be carried out \textit{verbatim} replacing $\Sp$ with $\cV$.
\end{proof}
\begin{remark}
Assuming that $\cV$ admits a left complete $Y$-structure, one could also expand the result of \cref{prop:modular_case} using the theory of \textit{faithfully flat monads}. Namely, for $\scrC$, $\cV$ and $\scrD$ as in \cref{assumption:standing}, with $\scrD$ further satisfying \cref{assumption:modular_case}, assume that for any morphism of test spaces $\varphi\colon U\to V$ the monadic adjunction $\varphi^\ast \dashv \varphi_\ast$ is faithfully flat in the sense of \cite[Definition D.6.4.1]{sag}. Then the functor \eqref{functor:categorification} is a sheaf. Already when $\cV$ is the $\infinity$-category of spectra, descendable morphisms and faithfully flat monads are not known to be compatible in general: for example, faithfully flat morphisms of commutative ring spectra always yield faithfully flat monads, but without some finiteness assumptions they are still not known to be descendable in the sense of \cref{def: descendable}. On the converse, descendable morphisms of commutative ring spectra need not to be faithfully flat (\cite[Remark D.3.3.3]{sag}).
\end{remark}
In this level of generality, studying the property of $1$-affineness for objects in $\widehat{\scrC}$ is very cumbersome. It turns out that both in the algebro-geometric and in the analytic settings the key property that we need is that, among all maps that form coverings for the topology on affine schemes, we can single out some particularly nice class of maps (the Zariski open immersions). Analytic open immersions are defined for analytic stacks generalizing some homological properties of pullback functors along finite open immersions (\cref{def:open immersion}). We use instead a different notion, already appeared in the literature as \textit{homotopy monormophism} (\cite{savage2024representabilitytheoremstacksderived, benbassat2024perspectivefoundationsderivedanalytic}), which asks instead for some good property with respect to self-intersections.
\begin{definition}
\label{def:top_open_immersion}\
Let $\scrC$ be an $\infinity$-category, and let $\alpha\colon U\to X$ be a morphism. We say that $\alpha$ is a \textit{Zariski morphism} if the diagonal map $U\to U\times_XU$ is an equivalence.
\end{definition}
\begin{remark}\
\begin{remarkenum}
\item \label{remark:topological_composition}
Let $\scrC$ be a $\infinity$-category. By a straightforward argument using the pasting law for pullbacks, it is easily seen that Zariski morphisms are stable under pullbacks and enjoy the following property: if $\varphi\colon Y\to X$ and $\psi\colon Z\to Y$ are morphisms of $\scrC$ such that $\varphi$ is a Zariski morphism, then $\varphi\circ\psi$ is Zariski if and only if $\psi$ is Zariski. In particular, Zariski morphisms are closed under composition, are right cancellative, and if $\varphi\colon Y\to X$ and $\psi\colon Z\to Y$ are Zariski morphisms then we have
\begin{align}
\label{eq:transitivity Zariski}
Z\simeq Z\times_YZ\simeq Z\times_XZ.
\end{align}
\item\label{remark:union}Let $\scrC$ be a $\infinity$-category, and let $\alpha\colon U\to X$ and $\beta\colon V\to X$ be two Zariski morphism. One can define the union of $U$ and $V$ as
\[
U\cup V\coloneqq U\coprod_{U\times_XV}V.
\]
It comes provided by a natural map $\gamma\colon U\cup V\to X$. If colimits in $\scrC$ are universal, we have that $\gamma$ is Zariski itself: indeed, in this case we have
\begin{align*}
(U \cup V)\times_X (U\cup V)&\coloneqq \lp U\coprod_{U\times_XV}V\rp \times_X\lp U\coprod_{U\times_XV}V\rp\\
&\simeq \lp (U\times_X U)\coprod_{U\times_XV}(V\times_XU)\rp\coprod_{U\times_XV}\lp (U\times_XV)\coprod_{U\times_XV}(V\times_XV)\rp
\\&\simeq U\coprod_{U\times_XV}V.
\end{align*}
It follows from \cref{remark:topological_composition} that also the inclusions $U\to U\cup V$ and $V\to U\cup V$ are Zariski morphisms.
\item \label{remark:zariski cover}Let $\scrC$, $\cV$ and $\scrD$ be as in \cref{assumption:modular_case}, and let $\iota\colon\scrC\to\widehat{\scrC}$ be a dense functor such that colimits in $\widehat{\scrC}$ are universal. If $\alpha\colon\iota(U)\to X$ and $\beta\colon\iota(V)\to X$ are Zariski morphisms which form a $\scrC$-affine atlas for $X$ (\cref{def:atlas}) then $\iota(U)\to \iota(U)\cup \iota(V)$ and $\iota(V)\to  \iota(U)\cup \iota(V)$ form a $\scrC$-affine atlas for $\iota(U)\cup\iota(V)$ as well: this is a clear consequence of the fact that $\iota(U)\to\iota(U)\cup \iota(V)$ and $\iota(V)\to\iota(U)\cup \iota(V)$ are Zariski morphisms, together with the fact that all squares in the following diagram
\[
\begin{tikzpicture}[scale=0.75]
\node (a) at (-6,0){$\iota(U)$};
\node (b) at (-2,0){$\iota(U)\cup \iota(V)$};
\node (c) at (2,0){$X$};
\node (a2) at (-6,2){$\iota(U)$};
\node (b2) at (-2,2){$\iota(U)\cup \iota(V)$};
\node (c2) at (2,2){$\iota(U)\cup \iota(V)$};
\node (a4) at (-6,4){$\iota(U)\times_X\iota(W)$};
\node (b4) at (-2,4){$\iota(W)$};
\node (c4) at (2,4){$\iota(W)$};
\draw[->,font=\scriptsize](a) to node[above]{}(b);
\draw[->,font=\scriptsize](b) to node[above]{}(c);
\draw[->,font=\scriptsize](a2) to node[above]{}(a);
\draw[->,font=\scriptsize](b2) to node[above]{}(b);
\draw[->,font=\scriptsize](c2) to node[above]{}(c);
\draw[->,font=\scriptsize](a2) to node[above]{}(b2);
\draw[->,font=\scriptsize](b2) to node[above]{}(c2);
\draw[->,font=\scriptsize](a4) to node[above]{}(a2);
\draw[->,font=\scriptsize](b4) to node[above]{}(b2);
\draw[->,font=\scriptsize](c4) to node[above]{}(c2);
\draw[->,font=\scriptsize](a4) to node[above]{}(b4);
\draw[->,font=\scriptsize](b4) to node[above]{}(c4);
\end{tikzpicture}
\]
are pullbacks. By an easy inductive argument (using the local property of covering sieves in $\infinity$-sites) it follows that if $\left\{\iota(U_i)\to X\right\}_{i=1}^k$ is a $\scrC$-affine atlas of Zariski morphisms, and $W_i$ denotes the union
\[
\iota(U_1)\cup\cdots\cup\widehat{\iota(U_i)}\cup\cdots\cup \iota(U_i),
\]
where the $i$-th element is omitted, then $\iota(U_i)\to X$ and $W_i\to X$ still form a $\scrC$-affine atlas of Zariski morphisms for $X$.
\item\label{remark:property_not_structure}If $\scrC$, $\cV$ and $\scrD$ are as in \cref{assumption:modular_case}, and if the functor
\[
D(-)\colon\scrC^{\op}\longrightarrow\CAlg(\cV)
\]
sends finite limits of $\scrC^{\op}$ to finite colimits of $\CAlg(\cV)$, then for any Zariski morphism of test spaces $\varphi\colon U\to V$ we have that the forgetful functor in the adjunction
\[
\adjunction{-\otimes_{D(V)}D(U)}{\Mod_{D(V)}(\cV)}{\Mod_{D(U)}(\cV)}{\oblv_{D(U)}}
\]
is fully faithful. Indeed, if $(X_i)_{i\in I}$ is a set of compact generators of $\cV$, the compact generators of $\Mod_{D(U)}(\cV)$ is given by $(X_i\otimes_{\cV}D(U))_{i\in I}$: in particular, on any generator of this form the counit of the adjunction $-\otimes_{D(V)}D(U)\dashv\oblv_{D(U)}$
\[
\oblv_{D(U)}(X_i\otimes_{\cV}D(U))\otimes_{D(V)}D(U)\longrightarrow X_i\otimes_{\cV}D(U)
\]
is readily seen to be an equivalence because $D(U)\otimes_{D(V)}D(U)\simeq D(U\times_VU)\simeq D(U)$. The equivalence on \textit{all} objects of $\Mod_{D(U)}(\cV)$ follows from the fact that both functors commute with colimits.
\end{remarkenum}
\end{remark}
The fundamental property of Zariski morphisms is that if an object $X$ in $\widehat{\scrC}$ admits a $\scrC$-affine atlas comprised of a single map $\iota(W)\to X$ which is a Zariski morphism, then $X$ is itself $\scrC$-affine. Moreover, if $X$ admits a finite $\scrC$-affine atlas $\left\{\alpha_i\colon\iota(W_i)\to X\right\}$ where each $\alpha_i$ is a Zariski morphism, then $X$ can be expressed as a \textit{finite} colimit of objects in the essential image of $\iota$.
\begin{proposition}
\label{prop:topological_open_immersion}
Let $\scrC$ be as in \cref{assumption:standing}, and let $\iota\colon\scrC\to\widehat{\scrC}$ be a dense functor. Assume that colimits in $\widehat{\scrC}$ are universal.\\
Consider an object $X$ in $\widehat{\scrC}$ admitting a $\scrC$-affine atlas $\left\{\alpha_i\colon \iota(W_i)\to X\right\}_{i=1}^k$ where each $\alpha_i$ is a Zariski morphism. Let $\varphi\colon \coprod_i \iota(W_i)\to X$ be the induced map from the disjoint union, and let $\check{\mathrm{C}}(\varphi)\to X$ be its \v{C}ech nerve. Then the maps in the natural composition
\[
\colim_{\bDelta^{\op}}\check{\mathrm{C}}(\varphi)\longrightarrow\colim_{\bDelta^{\op}_{\leqslant k-1}}\check{\mathrm{C}}(\varphi)\longrightarrow X
\]
are equivalences. In particular, if $k=1$ then $X\simeq\iota(W_1)$.
\end{proposition}
\begin{proof}
Since each $\alpha_i$ is a Zariski morphism, for all $h\geqslant k$ we have an equivalence
\[
\iota(W_{i_1})\times_X\cdots\times_X\iota(W_{i_h})\simeq\iota(W_1)\times_X\cdots\times_X\iota(W_k).
\]
Thus, the simplicial diagram corresponding to the \v{C}ech nerve $\check{\mathrm{C}}(\varphi)$ is constant for all $k\geqslant m-1$. In particular, it is $(m-1)$-skeletal, in the sense that the simplicial diagram $\bDelta^{\op}\to\widehat{\scrC}$ is a left Kan extension of its restriction to $\bDelta^{\op}_{\leqslant k}$. Since left Kan extensions along fully faithful functors do not alter colimits, the first map is naturally an equivalence.  So, it is sufficient to notice that this colimit is just $X$ because of \cref{remark:covering not affine}.
\end{proof}
\cref{prop:topological_open_immersion} offers a key technical advantage in working with objects of $\widehat{\scrC}$: when an object $X$ of $\widehat{\scrC}$ admits a \textit{finite} $\scrC$-affine atlas consisting only of Zariski morphisms, we are often allowed to run an inductive strategy on the cardinality of the $\scrC$-affine cover, starting by proving our statements when objects are $\scrC$-affine. We therefore introduce the following piece of terminology.
\begin{definition}
Let $\scrC$, $\cV$ and $\scrD$ be as in \cref{assumption:modular_case}, and let $\iota\colon\scrC\to\widehat{\scrC}$ be a dense functor.
\begin{defenum}
\item 
\label{def:quasi-compact}We say that an object $X$ in $\widehat{\scrC}$ is \textit{quasi-compact} if it admits a \textit{finite} affine $\scrC$-atlas $\left\{\alpha_i\colon\iota(W_i)\to X\right\}$ such that each $\alpha_i$ is a Zariski morphism.
\item\label{def:separated}We say that an object $X$ in $\widehat{\scrC}$ is \textit{separated} if all Zariski morphisms $\iota(U)\to X$ are $\scrC$-affine.
\end{defenum}
\end{definition}
\begin{remark}\
\begin{remarkenum}
    \item Given a dense functor $\iota\colon\scrC\to\widehat{\scrC}$, then any $\scrC$-affine object is always quasi-compact. If $\iota$ preserves pullbacks, then any $\scrC$-affine object is also separated. 
\item \label{transitivity separated}If $X$ is separated and $\alpha\colon Y\to X$ is a Zariski morphism, then $Y$ is separated as well. Indeed, if $\iota(U)\to Y$ is a Zariski morphism, then $\iota(U)\to X$ is a Zariski morphism and $\iota(U)\times_X Y\simeq \iota(U)$. So, for any map $\iota(W)\to Y$ we have an equivalence
\[
\iota(W)\times_Y\iota(U)\simeq\iota(W)\times_Y(Y\times_X\iota(U))\simeq \iota(W)\times_X\iota(U),
\]
which is $\scrC$-affine.
\item \label{transitivity quasi-compact}If $X$ is quasi-compact and $\alpha\colon Y\to X$ is a $\scrC$-affine Zariski morphism, then $Y$ is quasi-compact as well. Indeed, pulling back any finite $\scrC$-affine atlas of Zariski morphism for $X$ yields a finite $\scrC$-affine atlas of Zariski morphisms for $Y$.
\end{remarkenum}
\end{remark}
In order to obtain the best behavior of Zariski morphisms, for the rest of this section we add the following assumptions to our setting.
\begin{assumption}
 \label{assumption:modular_final}
Let $\scrC$, $\cV$ and $\scrD$ be as in \cref{assumption:modular_case}, and let $\iota\colon\scrC\to\widehat{\scrC}$ be a dense functor. We assume that $\iota$ preserves finite limits existing in $\scrC$, that colimits in $\widehat{\scrC}$ are universal, and that the sheaf of commutative algebras
\[
D(-)\colon\scrC^{\op}\longrightarrow\CAlg(\cV)
\]
sends finite limits of $\scrC$ to finite colimits in $\CAlg(\cV)$.
\end{assumption}
Our main $1$-affineness result of this section can then be stated as follows.
\begin{theorem}
\label{thm:modular_case}
Let $\iota\colon\scrC\to\widehat{\scrC}$, $\cV$ and $\scrD$ be as in \cref{assumption:modular_final}, and assume that $\cV$ is rigid. Let $X$ be an object in $\widehat{\scrC}$ admitting a finite affine $\scrC$-atlas $\left\{\varphi_i\colon\iota(W_i)\to X\right\}_{i=1}^k$ such that each $\varphi_i$ is a Zariski morphism, and such that for all subsets of indices $I=\left\{i_1,\ldots,i_n\right\}\subseteq \left\{1,\ldots k\right\}$ the functor
\[
\varphi_I^\ast\colon\scrD(\iota(W_{i_1})\times_X\cdots\times_X\iota(W_{i_n}))\longrightarrow\Rani\scrD(X) 
\]
admits a left adjoint $\varphi_{I,!}$. Then $X$ is $1$-affine.
\end{theorem}
\begin{remark}
The condition on the existence of left adjoints for pullback functors seems a bit unnatural. However, this is what happens for finite flat morphisms in algebraic geometry, for analytic open immersions in rigid analytic geometry (\cref{lem: left adjoint pullback opens-nuclear}), and in general for open immersions of affine analytic stacks (\cref{def:open immersion}).
\end{remark}
The proof of \cref{thm:modular_case} is quite complex and long, and relies on many auxiliary results that allow us to extend the good properties of Zariski open immersions in algebraic geometry to our abstract framework. We prove them following the same stream of ideas of \cite[Lemma 2.4.16]{mann2022p} and of \cite[Lemma 3.2]{soor2024sixfunctorformalismquasicoherentsheaves}.
\begin{proposition}
\label{prop:general_module}
Let $\iota\colon\scrC\to\widehat{\scrC}$, $\cV$ and $\scrD$ be as in \cref{assumption:modular_final}. Let $\alpha\colon Y\to X$ be a morphism between quasi-compact and separated objects of $\widehat{\scrC}$.
\begin{propenum}
\item\label{lemma:base change 1} (\emph{Base change}) Let $\beta\colon\iota(V) \to X$ be a morphism, and consider the diagram
               \[
        \begin{tikzpicture}[scale=0.75]
        \node (a) at (-3,2){$Y\times_{X}\iota(V)$};
        \node (b) at (3,2){$\iota(V)$};
        \node (c) at (-3,0){$Y$};
        \node (d) at (3,0){$X.$};
        \draw[->,font=\scriptsize](a) to node[above]{$\alpha'$}(b);
        \draw[->,font=\scriptsize](b) to node[right]{$\beta$}(d);
        \draw[->,font=\scriptsize](c) to node[above]{$\alpha$}(d);
        \draw[->,font=\scriptsize](a) to node[left]{$\beta'$}(c);
        \end{tikzpicture}
       \]
        Then for every object $M$ in $\Rani\scrD(Y)$ the canonical map 
        \begin{equation*}
            \beta^{\ast}\alpha_{\ast}M \longrightarrow \alpha'_{\ast}\beta'^{\ast}M
        \end{equation*}
        is an equivalence.
        \item\label{lemma:pushforward nuclear colimits}The functor $\alpha_\ast$ is conservative and commutes with colimits. In particular, it is a both monadic and comonadic functor.\\
        If $\alpha$ is Zariski, then $\alpha_\ast$ is moreover fully faithful.
    \end{propenum} 
\end{proposition}
\begin{proof}
Since $X$ is quasi-compact, we can consider a finite $\scrC$-affine atlas $\left\{\varphi_i\colon\iota(U_i)\to X\right\}_{i=1}^k$ consisting of Zariski morphisms. Since $\Rani\scrD(X)$ is a right Kan extension, using \cref{prop:topological_open_immersion} we have
\[
\Rani\scrD(X)\simeq\lim_{\bDelta_{\leqslant k}}\scrD(U_i),
\]
where the equivalence is given by the assignment
\[
N\mapsto\lim_{\bDelta_{\leqslant k}}\varphi_{i,\ast}\varphi^\ast_iN.
\]
Since we are in a stable setting (\cref{assump:stable}), all the functors and operations involved commute with finite limits. In particular, for most of the proof we can assume $X$ to be $\scrC$-affine.

We first prove \cref{lemma:base change 1}. Since $Y$ is quasi-compact we can find a finite $\scrC$-affine atlas of $Y$ consisting of Zariski morphisms, say $\left\{\iota(W_i)\to Y\right\}_{i=1}^k$. For every integer $n\geqslant0$, write
    \[
    Y_{i_{1}\cdots i_{n}}\quad\text{for}\quad\coprod_{\left\{i_1,\cdots i_n\right\}\subseteq I^{\times n}}\iota(W_{i_{1}})\times_{Y}\dots \times_{Y}\iota(W_{i_{n}}).
    \]
    Notice that each $\iota(W_{i_{1}})\times_{Y}\dots \times_{Y}\iota(W_{i_{n}})$ is still $\scrC$-affine. We denote with $j_{i_{1}\dots i_{n}}$ the induced map $Y_{i_{1}\dots i_{n}} \to Y$ and with $\alpha_{i_{1}\dots i_{n}}$ the composition $Y_{i_{1}\dots i_{n}} \to Y \to X$. In particular, the map $\alpha\colon Y\to X$ can be realized as a \textit{finite} colimit of the maps $\alpha_{i_1\cdots i_n}$. Thus, using once again the fact that $\Rani\scrD(-)$ sends colimits in $\widehat{\scrC}$ to limits and that we are working in the stable setting, we can write
   \begin{align*}
    \beta^\ast\alpha_\ast M&\simeq\beta^\ast\lp\lim_{\bDelta_{\leqslant k}}\alpha_{i_1\cdots i_n,\ast}\alpha^\ast_{i_1\cdots i_n}M\rp \\
    &\simeq\lim_{\bDelta_{\leqslant k}}\beta^\ast\alpha_{i_1\cdots i_n,\ast}\alpha^\ast_{i_1\cdots i_n}M\longrightarrow\alpha'_{\ast}\beta'^\ast\lp\lim_{\bDelta_{\leqslant k}}\alpha_{i_1\cdots i_n,\ast}j^\ast_{i_1\cdots i_n}M\rp\simeq\lim_{\bDelta_{\leqslant k}}\alpha'_\ast\beta'^\ast\alpha_{i_1\cdots i_n,\ast}k^\ast_{i_1\cdots i_n}M.
   \end{align*}
   So we are reduced to prove the statement when $\alpha\colon\iota(W)\to \iota(U)$ is a morphism between $\scrC$-affine objects; in this case, the claim is trivial because $Y\times_X\iota(V)$ is just the object $\iota(W\times_UV)$ and so the base change equivalence
   \[
   M\otimes_{D(W)}D(W\times_UV)\simeq M\otimes_{D(W)}D(W)\otimes_{D(U)}D(V)\longrightarrow M\otimes_{D(U)}D(V)
   \]
   reduces to the associativity of relative tensor products.
    
    
    To prove \cref{lemma:pushforward nuclear colimits}, we start by proving that $\alpha_\ast$ commutes with colimits. We observe that given a diagram $K\to \Rani\scrD(Y)$ selecting objects $M_k$, the colimit $\colim_K M_k$ can be written as
    \[
    \colim_{k\in K}M_k\simeq \colim_{k\in K}\lim_{[n]\in\bDelta_{\leqslant k}}j_{i_{1}\dots i_{n},\ast}j_{i_{1}\dots i_{n}}^\ast M_k\simeq\lim_{[n]\in\bDelta_{\leqslant k}}\colim_{k\in K}j_{i_{1}\dots i_{n},\ast}j_{i_{1}\dots i_{n}}^\ast M_k,
    \]
    where once again we used that the cosimplicial limit is actually finite. So, we are reduced to check the analogous statement when $Y\simeq\iota(W)$ is $\scrC$-affine: in this case, the claim is obvious since $\alpha_\ast$ is a forgetful functor between module $\infinity$-categories (\cref{remark:property_not_structure}). The same strategy implies that $\alpha_\ast$ is conservative.

    We finally prove the fully faithfulness of $\alpha_{\ast}$ in the case $\alpha\colon Y\to X$ is a Zariski morphism (which is the only part where we will not assume $X$ to be $\scrC$-affine). For any $M\in\Rani\scrD(Y)$ write
    \[
    M\simeq\lim_{[n]\in\bDelta_{\leqslant k}}j_{i_{1}\dots i_{n},\ast} j_{i_{1}\dots i_{n}}^{\ast}M.
    \]
So, in the counit morphism $\iota^\ast\iota_\ast M\to M$ for the adjunction $\alpha^\ast\dashv\alpha_\ast$, the source can be written as
\begin{align*}
\alpha^\ast\alpha_\ast M&\simeq\alpha^\ast\alpha_\ast\lp\lim_{[n]\in\bDelta}j_{i_{1}\dots i_{n},\ast} j_{i_{1}\dots i_{n}}^{\ast}M\rp\\&\simeq\lim_{[n]\in\bDelta}\lp\alpha^\ast\alpha_\ast j_{i_{1}\dots i_{n},\ast} j_{i_{1}\dots i_{n}}^{\ast}M\rp\\
&\simeq\lim_{[n]\in\bDelta_{\leqslant k}}\lp\alpha^\ast\alpha_{i_1\dots i_n,\ast}j_{i_{1}\dots i_{n}}^{\ast}M\rp\\
&\simeq\lim_{[n]\in\bDelta_{\leqslant k}}j_{i_1\dots i_n,\ast}j^\ast_{i_1\dots i_n}M\simeq M.
\end{align*}
Notice that in the equivalence in the last row we used both \cref{lemma:base change 1}, and the fact that $\iota(W_i)\times_XY\simeq \iota(W_i)$ (\cref{remark:topological_composition}). 
%
\end{proof}
\begin{proposition}\label{cor: projection formula of open}
Let $\iota\colon\scrC\to\widehat{\scrC}$, $\cV$ and $\scrD$ be as in \cref{assumption:modular_final}. Let $\alpha\colon Y \to X$ be a morphism between quasi-compact and separated objects of $\widehat{\scrC}$. Then for every $N$ in $\Rani\scrD(Y)$ and $M$ in $\Rani\scrD(X)$ the natural map 
    \begin{equation}\label{eq: projection formula}
        \alpha_{\ast}(N \otimes_{\Rani\scrD(Y)} \alpha^{\ast}M) \longrightarrow \alpha_{\ast}N \otimes_{\Rani\scrD(X)}M
    \end{equation} 
    is an equivalence.
\end{proposition}
\begin{proof}
 First, assume that $Y$ is itself $\scrC$-affine: write it $Y\simeq\iota(V)$. Since $X$ is quasi-compact, we can consider a finite $\scrC$-affine atlas $\{\varphi_i\colon \iota(W_i)\to X\}_{i=1}^k$ consisting of Zariski morphisms; let us denote with $j_{i}\colon \iota(W_i)\to X$ the inclusion. Arguing as in the previous proof, by descent along the $\scrC$-affine atlas we can check whether \eqref{eq: projection formula} is an equivalence 
    locally on the target. So, consider the  pullback diagram
    \[
    \begin{tikzpicture}[scale=0.75]
    \node (a) at (-3,2){$\iota(V)\times_X\iota(W_i)$};
    \node (b) at (3,2){$\iota(W_i)$};
    \node (c) at (-3,0){$\iota(V)$};
    \node (d) at (3,0){$X.$};
    \draw[->,font=\scriptsize] (a) to node[above]{$\alpha'$}(b);
    \draw[->,font=\scriptsize] (a) to node[left]{$j_i'$}(c);
    \draw[->,font=\scriptsize] (b) to node[right]{$j_i$}(d);
    \draw[->,font=\scriptsize] (c) to node[above]{$\alpha$}(d);
    \end{tikzpicture}
    \]
Since $X$ is separated, $\iota(W)\times_X\iota(W_i)$ is $\scrC$-affine: let us denote it by $\iota(U_i)$. We now have
\begin{align*}
j_i^{\ast}\alpha_\ast(N\otimes_{D(V)}\alpha^\ast M)&\simeq\alpha'_\ast j'^\ast_i(N\otimes_{D(V)}\alpha^\ast M)\\
&\simeq\alpha'_\ast(j'^\ast_i N\otimes_{D(V)}j'^\ast_i\alpha^\ast M)\\
&\simeq \alpha'_\ast j'^\ast_i M\otimes_{D(V)}j^\ast_i N\\
&\simeq j_i^\ast (\alpha_\ast M\otimes_{\Rani\scrD(X)}N).
\end{align*}
We highlight that we used the base change property of \cref{lemma:base change 1}, and the fact that for $\scrC$-affine objects the projection equivalence in \eqref{eq: projection formula} holds trivially (it boils down to the fact that $N\otimes_{D(V)}(D(V)\otimes_{D(U)}M)$ is equivalent to $N\otimes_{D(U)}M$). For a general $Y$, we consider a finite $\scrC$-affine atlas of $Y$ given by Zariski morphisms, and write each $N$ as a finite limit along such cover; thus, we reduce ourselves to the case when $Y$ is $\scrC$-affine, which we have just proved.
\end{proof}
Concatenating \cref{prop:general_module,cor: projection formula of open} we immediately obtain the following.
\begin{corollary}
\label{lemma:nuclear-modular-relative}
Let $\iota\colon\scrC\to\widehat{\scrC}$, $\cV$ and $\scrD$ be as in \cref{assumption:modular_final}. Let $ \alpha \colon\iota(U) \to X$ be a morphism to a quasi-compact and separated object of $\widehat{\scrC}$. We have an equivalence of $\infinity$-categories
    \[
    \scrD(U)\simeq\Mod_{\alpha_\ast D(U)}(\Rani\scrD(X)).
    \]
\end{corollary}
\begin{proof}
    The proof is a standard argument based on Barr--Beck--Lurie's monadicity theorem. Indeed, by \Cref{lemma:pushforward nuclear colimits} the functor $\iota_\ast$ is monadic, and one simply needs to do a routine computation to show that this monad corresponds to the monad $-\otimes_{\Rani\scrD(X)}\alpha_\ast D(U)$, using the projection formula.
\end{proof}
\begin{proposition}
\label{thm: Nuc is rigid}
Let $\iota\colon\scrC\to\widehat{\scrC}$, $\cV$ and $\scrD$ be as in \cref{assumption:modular_final}. If $X$ is a quasi-compact and separated object of $\widehat{\scrC}$ and $\cV$ is a rigid symmetric monoidal $\infinity$-category, then $\Rani\scrD(X)$ is rigid as well.
\end{proposition}
\begin{proof}
We argue inductively on the cardinality of $\scrC$-affine atlas $\left\{\varphi_i\colon\iota(W_i)\to X\right\}_{i=1}^k$, where each $\varphi_i$ is a Zariski morphism.\\
If $k=1$, then $X$ is $\scrC$-affine and in this case the assertion follows from the more general fact that if $\cV$ is rigid and $A$ is a commutative algebra object in $\cV$, then $\Mod_A(\cV)$ is rigid. Let now $k>1$, and assume that $\Rani\scrD(Y)$ is rigid whenever $Y$ admits a $\scrC$-affine atlas consisting of $k-1$ Zariski morphisms. Let us write
\[
Z\coloneqq\iota(W_1)\cup\cdots\cup\iota(W_{k-1}).
\]
By the inductive hypothesis (together with \cref{remark:union}) both $\Rani\scrD(Z)$ and $\Rani\scrD(Z\times_X\iota(W_k))$ are rigid. By descent along the $\scrC$-affine atlas, $\Rani\scrD(X)$ sits in a pullback diagram
 \begin{equation}\label{eq: pullback nuc cat}
     \begin{tikzpicture}[scale=0.75,baseline=(current  bounding  box.center)]
 \node (a) at(-3,2){$\Rani\scrD(X)$};
 \node (b) at (3,2){$\Rani\scrD(Z)$};
 \node (c) at (-3,0){$\scrD(W_k)$};
 \node (d) at (3,0){$\Rani\scrD(Z\times_X\iota(W_k)).$};
 \draw[->,font=\scriptsize](a) to node[above]{}(b);
 \draw[->,font=\scriptsize](a) to node[above]{}(c);
 \draw[->,font=\scriptsize](b) to node[above]{}(d);
 \draw[->,font=\scriptsize](c) to node[above]{}(d);
 \end{tikzpicture}
 \end{equation} 
\Cref{prop:general_module,cor: projection formula of open} imply that the pullback functors $\Rani\scrD(X)\to\Rani\scrD(Z)$ and $\Rani\scrD(X)\to\Rani\scrD(Z\times_X \iota(W_k))$ are internally left adjoints, in the sense of \cite[Definition 1.9]{ramzi2024dualizablepresentableinftycategories}. So we can simply apply \cite[Corollary 4.5]{ramzi2024dualizablepresentableinftycategories} and obtain that $\Rani\scrD(X)$ is a pullback also inside the full sub-$\infinity$-category $\mathrm{Pr}_{\mathrm{st}}^{\mathrm{L,dbl}}\subseteq\LPr_{\mathrm{st}}$ of \textit{dualizable} presentable stable $\infinity$-categories. Since the inclusion of rigid presentable symmetric monoidal $\infinity$-categories into dualizable presentable $\infinity$-categories preserves limits (\cite[Corollary 4.85]{ramzi2024locallyrigidinftycategories}), we can view the pullback in \eqref{eq: pullback nuc cat} as a limit of rigid $\infinity$-categories. In particular $\Rani\scrD(X)$ is rigid.
\end{proof}
The last result that we need in the proof of \cref{thm:modular_case} is the following lemma, due to Gaitsgory.
\begin{lemma}\label{lem: right adjointable modular case}
Let $\iota\colon\scrC\to\widehat{\scrC}$, $\cV$ and $\scrD$ be as in \cref{assumption:modular_final}. Let $\alpha\colon Y\to X$ be a morphism in $\widehat{\scrC}$. Assume that there is a $\scrC$-affine atlas $\{\iota(W_{i})\to X\}_{i\in I}$ of $X$ such that for every index $i$ the canonical map
\begin{equation*}
    \scrD(W_{i})\otimes_{\Rani\scrD(X)} \Rani\scrD(Y) \longrightarrow \Rani\scrD(\iota(W_{i})\times_{X}Y)
\end{equation*}
is an equivalence. Assume moreover that $\scrD(Y)$ is dualizable as a categorical $\scrD(X)$-module. Then the diagram
\[
\begin{tikzpicture}[scale=0.75]
\node(a) at (3,2.5){$\mathrm{Ran}_{\iota}\scrD^{\mathrm{cat}}(X)$};
\node (c) at (3,0){$\mathrm{Ran}_{\iota}\scrD^{\mathrm{cat}}(Y)$};
\node (b) at (-3,2.5){$\LPr_{\mathrm{Ran}_{\iota}\scrD(X)}$};
\node (d) at (-3,0){$\LPr_{\mathrm{Ran}_{\iota}\scrD(Y)}$};
\draw[->,font=\scriptsize] (b) to node[above]{$\widehat{(-)}_X$}(a);
\draw[->,font=\scriptsize] (a) to node[right]{$\widetilde{\alpha}^\ast$}(c);
\draw[->,font=\scriptsize] (b) to node[left]{$\alpha^\ast$}(d);
\draw[->,font=\scriptsize] (d) to node[below]{$\widehat{(-)}_Y$}(c);
\end{tikzpicture}
\]
is vertically right adjointable.
\end{lemma}
\begin{proof}
Completely analogous to the proof of \cite[Proposition 3.2.6 (ii)]{1affineness}. 
\end{proof}
However, in practice, often one would wish for some concrete way to check whether the K\"unneth formula in the hypotheses of \cref{lem: right adjointable modular case} is satisfied. The following offers a nice criterion when $\cV$ is a rigid symmetric monoidal $\infinity$-category.
\begin{proposition}\label{prop: kunneth for Nuc using opens}
    Let $\iota\colon\scrC\to\widehat{\scrC}$, $\cV$ and $\scrD$ be as in \cref{assumption:modular_final}, and assume $\cV$ to be rigid as a symmetric monoidal $\infinity$-category. Let $\alpha\colon Y\to X$ be a morphism between quasi-compact and separated objects of $\widehat{\scrC}$. For every morphism $\varphi\colon \iota(V) \to X$, the canonical functor
    \begin{equation*}
\Rani\scrD(Y)\otimes_{\Rani\scrD(X)}\scrD(V)\longrightarrow\Rani\scrD(Y\times_X\iota(V))
    \end{equation*}
    is an equivalence.
\end{proposition}
\begin{proof}
First assume that $Y\simeq\iota(U)$ is itself $\scrC$-affine. Now $\iota(U)\times_X\iota(V)$ is $\scrC$-affine because $X$ is separated: write $\iota(W)$ for such pullback, and let $\alpha'\colon\iota(W)\to \iota(V)$ and $\varphi'\colon\iota(W)\to\iota(U)$ denote the corresponding morphisms. Using \cref{lemma:nuclear-modular-relative}, we can write
\[
\scrD(U)\simeq\Mod_{\alpha_\ast D(U)}(\Rani\scrD(X))
\]
and so the tensor product $\scrD(U)\otimes_{\Rani\scrD(X)}\scrD(V)$ is equivalent to $\Mod_{\alpha_\ast D(U)}(\scrD(V))$ via \cite[Theorem 4.8.4.6]{HA}. This in turn is by definition $\Mod_{\varphi^\ast\alpha_\ast D(U)}(\scrD(V))$, since $\scrD(V)$ is tensored over $\Rani\scrD(X)$ via the symmetric monoidal pullback $\varphi^\ast$. 

On the other hand, using the base change property (\cref{lemma:base change 1}), we can write
\begin{align*}
\scrD(\iota(U)\times_X\iota(V))&\simeq\Mod_{\alpha'_\ast D(W)}(\scrD(V))\\
&\simeq\Mod_{\alpha'_\ast\varphi'^\ast D(U)}(\scrD(V))\\
&\simeq\Mod_{\varphi^\ast\alpha_\ast D(U)}(\scrD(V)).
\end{align*}
 For the general case, consider a $\scrC$-affine atlas $\left\{\iota(U_i)\to Y\right\}_{i\in I}$ consisting of Zariski morphisms. Since $\scrD(V)$ and $\Rani\scrD(X)$ are rigid (\cref{thm: Nuc is rigid}) it follows that $\scrD(V)$ is dualizable as a categorical $\Rani\scrD(X)$-module, so we have
\begin{align*}
 \Rani\scrD(Y)\otimes_{\Rani\scrD(X)}\scrD(V)&\simeq\lp\lim_{i\in I}\scrD(U_i)\rp\otimes_{\Rani\scrD(X)}\scrD(V)\\
 &\simeq\lim_{i\in I}\lp\scrD(U_i)\otimes_{\Rani\scrD(X)}\scrD(V)\rp\\
 &\simeq\lim_{i\in I}\Rani\scrD(\iota(U_i)\times_X\iota(V))
 \\&\simeq\Rani\scrD(Y\times_X\iota(V)).
\end{align*}
\end{proof}
\begin{proof}[Proof of \cref{thm:modular_case}]
    Let $X$ be quasi-compact and separated, and let $\left\{\varphi_i\colon\iota(W_i)\to X\right\}_{i=1}^k$ be a finite $\scrC$-affine atlas as in the statement. We will prove the theorem by induction on the cardinality $k$ of the affine $\scrC$-atlas.

    If $k=1$, $X$ is $\scrC$-affine so there is nothing to prove. For $k>1$, let$$Y\coloneqq\iota(W_2)\cup\cdots\cup\iota(W_{k})$$denote the union of the last $k-1$ elements of the atlas; let $\psi\colon Y\to X$, $\alpha_1\colon Y\times_X\iota(W_1)\to \iota(W_1)$ and $\beta\colon Y\times_X\iota(W_1)\to Y$ denote the natural maps (which are all Zariski, in virtue of \cref{remark:union}). Consider the following diagram.
    \begin{equation}
    \label{diagram:importante}
    \begin{tikzpicture}[scale=0.75,baseline=(current  bounding  box.center)]
\node (a) at (-4,7){$\LPr_{\Rani\scrD(X)}$};
\node (b) at (6,7){$\Rani\scrD^{\mathrm{cat}}(X)$};
\node (c) at (-4,3){$\LPr_{\scrD(W_1)}$};
\node (d) at (6,3){$\scrD^{\mathrm{cat}}(W_1)$};
\node(a1) at (-8,4.5){$\LPr_{\Rani\scrD(Y)}$};
\node (b1) at (2,4.5){$\Rani\scrD^{\mathrm{cat}}(Y)$};
\node (c1) at (-8,0.5){$\LPr_{\Rani\scrD(Y\times_X\iota(W_1))}$};
\node (d1) at (2,0.5){$\Rani\scrD^{\mathrm{cat}}(Y\times_X\iota(W_1))$};
\draw[->,font=\scriptsize] (a) to node[above]{$\widehat{(-)}_X$}(b);
\draw[->,font=\scriptsize] (a) to node[above left]{$\varphi_1^\ast$}(c);
\draw[->,font=\scriptsize] (b) to node[right]{$\widetilde{\varphi}_1^\ast$}(d);
\draw[->,font=\scriptsize] (a1) to node[above]{$\simeq$}(b1);
\draw[->,font=\scriptsize] (a1) to node[left]{$\beta^\ast$}(c1);
\draw[->,font=\scriptsize] (b1) to node[below left]{$\widetilde{\beta}^\ast$}(d1);
\draw[->,font=\scriptsize] (c) to node[above]{$\simeq$}(d);
\draw[->,font=\scriptsize] (c1) to node[above]{$\simeq$}(d1);
\draw[->,font=\scriptsize] (a) to node[above left]{$\psi^\ast$}(a1);
\draw[->,font=\scriptsize] (b) to node[above left]{$\widetilde{\psi}^\ast$}(b1);
\draw[->,font=\scriptsize] (c) to node[above left]{$\alpha_1^\ast$}(c1);
\draw[->,font=\scriptsize] (d) to node[above left]{$\widetilde{\alpha}_1^\ast$}(d1);
    \end{tikzpicture}
 \end{equation}
 Everything is naturally commutative, and the equivalences are all determined by the fact that the adjunction \eqref{adjunction:1_affineness} is an equivalence for $\iota(W_1)$ (being $\scrC$-affine) and for $Y$ and $Y\times_X\iota(W_1)$ (by the induction hypothesis). We want to prove that $\widehat{(-)}_X$ is an equivalence: for this purpose, we will prove that for all objects $\scrE$ in $\LPr_{\Rani\scrD(X)}$ and for all objects $(\scrE_U)_{\gamma\colon \iota(U)\to X}$ the unit \eqref{functor:unit} and the counit \eqref{functor:counit} are equivalences.

 Notice that both the top and the bottom faces of the diagram \eqref{diagram:importante} are right adjointable, because of \cref{lem: right adjointable modular case} (the hypotheses are satisfied thanks to \cref{prop: kunneth for Nuc using opens}). In particular we have the following.
 \begin{enumerate}
     \item The functor $\psi^\ast$ sends the unit functor for the adjunction $\widehat{(-)}_X\dashv \Gamma(X,-)$ to the unit functor for the adjunction  $\widehat{(-)}_Y\dashv \Gamma(Y,-)$, which is an equivalence.
     \item The functor $\varphi_1^\ast$ sends the unit functor for the adjunction $\widehat{(-)}_X\dashv \Gamma(X,-)$ to the unit functor for the adjunction  $\widehat{(-)}_{\iota(W_1)}\dashv \Gamma(\iota(W_1),-)$, which is an equivalence.
     \item The functor $\widetilde{\psi}^\ast$ sends the counit functor for the adjunction $\widehat{(-)}_X\dashv \Gamma(X,-)$ to the counit functor for the adjunction  $\widehat{(-)}_Y\dashv \Gamma(Y,-)$, which is an equivalence.
     \item The functor $\widetilde{\varphi}_1^\ast$ sends the counit functor for the adjunction $\widehat{(-)}_X\dashv \Gamma(X,-)$ to the unit functor for the adjunction  $\widehat{(-)}_{\iota(W_1)}\dashv \Gamma(\iota(W_1),-)$, which is an equivalence.
 \end{enumerate}
 So it is sufficient to prove that the functors $\psi^\ast$, and $\varphi_1^\ast$ (respectively $\widetilde{\psi}^\ast$ and $\widetilde{\varphi}_1^\ast$) are jointly conservative functors. For $\widetilde{\psi}^\ast$ and $\widetilde{\varphi}_1^\ast$, this is a obvious consequence of the fact that $Y$ and $\iota(U_1)$ form a $\scrC$-affine atlas for $X$, with respect to which which $\Rani\scrD^{\mathrm{cat}}$ satisfies descent. For $\psi^\ast$ and $\varphi_1^\ast$, it is enough to show that for each categorical $\Rani\scrD(X)$-module $\scrE$ the diagram
\begin{equation}
\label{diagram:cartesian}
\begin{tikzpicture}[scale=0.75,baseline=(current  bounding  box.center)]
 \node (a) at (-4,2.5){$\scrE$};
\node (b) at (4,2.5){$\scrE\otimes_{\Rani\scrD(X)}\Rani\scrD(Y)$};
\node (c) at (-4,0){$\scrE\otimes_{\Rani\scrD(X)}\scrD(W_1)$};
\node (d) at (4,0){$\scrE\otimes_{\Rani\scrD(X)}\Rani\scrD(Y\times_X\iota(W_1))$};
\draw[->,font=\scriptsize] (a) to node[above]{$\mathrm{id}_{\scrE}\otimes \psi^\ast$} (b);
\draw[->,font=\scriptsize] (a) to node[left]{$\mathrm{id}_{\scrE}\otimes \varphi_1^\ast$} (c);
\draw[->,font=\scriptsize] (b) to node[right]{$\mathrm{id}_{\scrE}\otimes \alpha_1^\ast$} (d);
\draw[->,font=\scriptsize] (c) to node[below]{$\mathrm{id}_{\scrE}\otimes \beta^\ast$}(d);
 \end{tikzpicture}
\end{equation}
 is Cartesian. For this, we will need the following lemma.
 \begin{lemma}
 \label{lemma:auxiliary finale}
 Let $\iota\colon\scrC\to\widehat{\scrC}$, $\cV$, $\scrD$ and $X$ be as in the statement of \cref{thm:modular_case}. Then, for every index $i$, for every quasi-compact and separated object $Y$ and for every composition
 \[
 \varphi_i\colon\iota(W_i)\overset{\psi_i}{\longrightarrow} Y\overset{\alpha}{\longrightarrow} X,
 \]
where $\alpha$ is Zariski, the functor $\psi_i^\ast\colon\Rani\scrD(Y)\to\scrD(W_i)$ admits a left adjoint $\psi_{i,!}$. If $Y$ is of the form
 \[
 Y\coloneqq\bigcup_{i_1,\ldots,i_n} \iota(W_i)
 \]
 for some subset of indices $\left\{i_1,\ldots, i_n\right\}\subseteq\left\{1,\ldots, k\right\}$, the same holds for $\alpha^\ast$.
 \end{lemma}
\begin{proof}
We start by proving the existence of the left adjoint $\psi_{i,!}\colon\scrD(W_i)\to \Rani\scrD(Y)$. This is characterized by the property that for all $M$ in $\scrD(W_i)$ and for all $N$ in $\Rani\scrD(Y)$, we have
\begin{align*}
\Map_{\Rani\scrD(Y)}(\psi_{i,!}M N)&\simeq\Map_{\scrD(W_i)}(M,\psi_i^\ast N).
\end{align*}
But since $\alpha_\ast\colon \Rani\scrD(Y)\to \Rani\scrD(X)$ is fully faithful (\cref{lemma:pushforward nuclear colimits}) we can write
\begin{align*}
\Map_{\scrD(W_i)}(M,\psi_i^\ast N)&\simeq\Map_{\scrD(W_i)}(M,\psi_i^\ast\alpha^\ast\alpha_\ast N)\\
&\simeq\Map_{\scrD(W_i)}(M,\varphi^\ast_i\alpha_\ast N)\\
&\simeq \Map_{\Rani\scrD(X)}(\varphi_{i,!}M,\alpha_\ast N)\\
&\simeq \Map_{\scrD(W_i)}(\alpha^\ast\varphi_{i,!}M,N).
\end{align*}
So we just define $\psi_{i,!}$ to be $\alpha^\ast\circ\varphi_{i,!}$.

For the existence of $\alpha_!$ in the case $Y$ is a union of objects in the atlas, we proceed by induction on $n$, where the case $n=1$ is trivial (indeed, in this case, we have an equivalence $Y\simeq \iota(W_n)$). For $n\geqslant 2$, we consider the affine $\scrC$-atlas of Zariski morphisms given by all inclusions $\iota(W_{i_j})\to Y$, and write
\[
\Rani\scrD(Y)\simeq\lim_{j=1,\ldots n}\scrD(W_{i_j}).
\]
By hypothesis, each functor in the above diagram admits a left adjoint; so in particular we can write
\[
\Rani\scrD(Y)\simeq\colim_{j=1,\ldots, n}\scrD(W_{i_j}).
\]
So we just let $\alpha_!$ be the colimit of all the functors $\varphi_{I,!}$.
\end{proof}
Using \cref{lemma:auxiliary finale}, we see that both $\beta^\ast$ and $\alpha^\ast_1$ admit left adjoints $\alpha_{1,!}$ and $\beta_!$. Passing to the diagram of left adjoints (and using the compatibility of colimtis and tensor products), from the diagram \eqref{diagram:cartesian} we obtain a pushout square
\[
\begin{tikzpicture}[scale=0.75]
 \node (a) at (-6,2.5){$\scrE\otimes_{\Rani\scrD(X)}\lp\scrD(W_1)\coprod_{\Rani\scrD(Y\times_X\iota(W_1))}\Rani\scrD(Y)\rp$};
\node (b) at (6,2.5){$\scrE\otimes_{\Rani\scrD(X)}\Rani\scrD(Y)$};
\node (c) at (-6,0){$\scrE\otimes_{\Rani\scrD(X)}\scrD(W_1)$};
\node (d) at (6,0){$\scrE\otimes_{\Rani\scrD(X)}\Rani\scrD(Y\times_X\iota(W_1)),$};
\draw[->,font=\scriptsize] (b) to node[above]{} (a);
\draw[->,font=\scriptsize] (c) to node[left]{} (a);
\draw[->,font=\scriptsize] (d) to node[right]{$\mathrm{id}_{\scrE}\otimes \alpha_{1,!}$} (b);
\draw[->,font=\scriptsize] (d) to node[below]{$\mathrm{id}_{\scrE}\otimes \beta_!$}(c);
 \end{tikzpicture}
\]
so our claim follows if
\[
\scrD(W_1)\coprod_{\Rani\scrD(Y\times_X\iota(W_1))}\Rani\scrD(Y)\simeq\Rani\scrD(X).
\]
Again, by passing to the left adjoints, this is equivalent to saying that
\[
\Rani\scrD(X)\simeq\scrD(W_1)\times_{\Rani\scrD(Y\times_X\iota(W_1))}\Rani\scrD(Y),
\]
which is true because the $\iota(W_i)$'s are a $\scrC$-affine atlas of $X$, and they form a $\scrC$-affine atlas for $Y$ as well.

\end{proof}

\section{Recollections on analytic geometry after Clausen--Scholze}
Our main goal is to apply the general formalism established in \cref{sec:1affinenessgeneral} to the setting of rigid and analytic geometry. More precisely:
\begin{enumerate}
    \item The site $\scrC$ will be either the opposite of the $\infinity$-category $\Anring$ equipped with the universal !-able topology, or the $\infinity$-category $\mathrm{Afd}$ equipped with the analytic topology. The former is the opposite of the $\infinity$-category of condensed analytic rings (equivalently: the $\infinity$-category of affine analytic stacks), while the latter is the $\infinity$-category of affinoid spaces over a complete non-archimedean field $\Bbbk$.
     \item The sheaf of presentably symmetric monoidal and stable $\infinity$-categories $\scrD$ will be either the functor $\scrD(-)$ or its sub-functor $\Nuc(-)$. The former assigns to an analytic ring $A$ its derived $\infinity$-category $\scrD(A)$ of complete modules, while the latter is defined for analytic rings coming from affinoid algebras and selects the sub-$\infinity$-category $\Nuc(A)\subseteq\scrD(A)$ spanned by \textit{nuclear} modules.
     
    \item Accordingly, the stable and presentably symmetric monoidal $\infinity$-category $\cV$ will be either the derived $\infinity$-category $\scrD(\Z)$, where $\Z$ is equipped with the trivial analytic structure,  or the $\infinity$-category $\Nuc(\Bbbk)$.
    \item The dense functor $\iota\colon\scrC\to\widetilde{\scrC}$ will be either the Yoneda embedding (in the analytic case) or the inclusion of affinoid spaces in the $\infinity$-category of all rigid varieties.
\end{enumerate}
While the theory of rigid varieties is classical (see for example \cite{Bosch_rigid_geometry} for an introduction), the theory of analytic stacks has been only recently developed in the framework of condensed mathematics. Thus, before proceeding in studying the problem of $1$-affineness in the context of analytic stacks, we start by reviewing some fundamental definitions and fixing the notations that will be used in \cref{sec:1_affineness_condensed}. We stress that this section does not bear any addition or improvement to the existing literature, but we include it in the paper for the reader's convenience: if they are already acquainted with the language of condensed mathematics, they may safely skip it.  

We mainly adapt the conventions and terminology from \cite{Clausen_Scholze_lectures}; further references are provided in the text when needed.
\subsection{Condensed objects in \texorpdfstring{$\infinity$}{oo}-categories}
Recall that a \textit{pro-finite set} is by definition a pro-finite object in the category $\mathrm{Fin}$ of finite sets. More explicitly, a pro-finite set is a set $X$ which can be realized as a projective limit of finite sets, i.e., it is a set admitting a presentation
\[
X\cong\lim_{i\in I}X_i
\]
where $I$ is a cofiltered diagram and each $X_i$ is finite. We denote the category of pro-finite sets as $\Prof$.
\begin{definition}
A pro-finite set $X$ is \textit{light} if it can be presented as a limit of finite sets over a \textit{countable} cofiltered diagram. 
\end{definition}
\begin{remark}
\label{remark:stone_duality}
Equipping finite sets with the discrete topology, every pro-finite comes equipped with the inverse limit topology. Indeed, by Stone duality, pro-finite sets are precisely totally disconnected compact Hausdorff topological spaces. In \cref{sec:main_applications}, we shall often blur this distinction and consider pro-finite sets as objects inside the category of topological spaces.

Under this equivalence, light pro-finite sets correspond to those totally disconnected compact Hausdorff topological spaces which are moreover metrizable \cite[Lecture 2]{Clausen_Scholze_lectures}.
\end{remark}
Light pro-finite sets define a full sub-category of $\Prof$, that we denote as $\mathrm{ProFin}^{\mathrm{light}}$. This category naturally admits a Grothendieck topology, for which a cover is a finite collection of jointly surjective maps.
\begin{definition}
Let $\scrC$ be any $\infinity$-category. A \textit{condensed object} in $\scrC$ is a sheaf from the category of light pro-finite sets to $\scrC$. We denote the $\infinity$-category of condensed objects in $\scrC$ as $\scrC^{\mathrm{cond}}$.
\end{definition}
\begin{remark}
\label{remark:condensed_stable}
If $\scrC$ is presentable, then $\scrC^{\mathrm{cond}}$ is presentable as well. This follows from the fact that $\Prof^{\mathrm{light}}$ is small enough to let $\smash{\Fun{\lp(\Prof^{\mathrm{light}})^{\op},\scrC\rp}}$ be presentable itself -- light pro-finite sets are singled out precisely for this reason -- and so $\scrC^{\mathrm{cond}}$ defines a topological localization of $\smash{\Fun{\lp(\Prof^{\mathrm{light}})^{\op},\scrC\rp}}$.

At the same time, if $\scrC$ is stable then so is $\scrC^{\mathrm{cond}}$. Again, this follows from the fact that $\scrC^{\mathrm{cond}}$ is a topological localization of the stable $\infinity$-category $\smash{\Fun{\lp(\Prof^{\mathrm{light}})^{\op},\scrC\rp}}$.
\end{remark}
We will be mainly interested in the $\infinity$-categories of condensed anima, of derived condensed modules over a ground derived commutative ring $\Bbbk$, and of derived condensed $\Einf$-$\Bbbk$-algebras. The latter in particular (almost) play the role of rings of functions of affine analytic stacks in condensed mathematics (see \cref{def:analytic_ring}).
\begin{remark}
In the foundational work \cite{toen2008homotopical} To\"en-Vezzosi developed derived geometry using the model of simplicial commutative algebras: from a more model-independent perspective, these correspond to \textit{derived commutative rings} in the sense of \cite[$\S4$]{raksit2020hochschildhomologyderivedrham}. On the other hand, the approach of Lurie is based on the notion of $\Einf$-algebras and leads to spectral algebraic geometry.

Along the same lines, in condensed mathematics one can define a notion of what should morally be a "simplicial condensed algebra": in the $\infinity$-categorical setting, this theory has been developed and used in \cite{analytic,clausenlectures,mann2022p,camargo2024analytic}, and such objects have been called \textit{condensed animated algebras}. However, it is possible to develop, in a similar way, an analogous theory using the notion of condensed $\Einf$-algebras, as shown in \cite{fedeli2023topological}. Our point of view follows the latter approach.

Notice that in characteristic zero the two points of view coincide, essentially because in characteristic zero $\Einf$-algebras and derived commutative rings are algebras for the same monad over $\Ani$.
\end{remark}

Let $\scrC$ be a presentably symmetric monoidal $\infinity$-category. Then the $\infinity$-category of $\scrC$-valued presheaves over $\Prof^{\mathrm{light}}$ admits a point-wise tensor product turning it into a presentably symmetric monoidal $\infinity$-category as well. Hence, one can define a presentably symmetric monoidal structure also on $\scrC^{\mathrm{cond}}$, given by sheafifying the point-wise monoidal structure on $\Fun{\lp(\Prof)^{\mathrm{light}},\scrC\rp}$. 

Thus, it is natural to ask what is the relationship between condensed commutative algebras in $\scrC$, and commutative algebras in $\scrC^{\mathrm{cond}}$. The following proposition shows that there is no ambiguity whatsoever. 
\begin{proposition}
\label{prop:commutative_algebras_condensed}
Let $\scrC$ be a presentably symmetric monoidal $\infinity$-category. Then there is an equivalence
    \begin{equation*}
        \CAlg(\scrC)^{\mathrm{cond}}\simeq\CAlg(\scrC^{\mathrm{cond}})
    \end{equation*}
    between condensed $\Einf$-algebras in $\scrC$, and $\Einf$-algebras in the $\infinity$-category of condensed objects of $\scrC$. 
\end{proposition}
\begin{remark}
\cref{prop:commutative_algebras_condensed} applies in particular to the case when $\scrC$ is either the $\infinity$-category of derived $\Bbbk$-modules, or the $\infinity$-category of \textit{connective} derived $\Bbbk$-modules.
\end{remark}
\begin{proof}[{Proof of \cref{prop:commutative_algebras_condensed}}]Since both $\scrC$ and $\CAlg(\scrC)$ admit all limits, \cite[Proposition 1.3.1.7]{sag} yields equivalences
\[
\scrC^{\mathrm{cond}}\coloneqq\mathrm{Sh}(\mathrm{ProFin}^{\mathrm{light}},\scrC)\simeq\Fun^{\mathrm{R}}{\lp(\CondAni)^{\op},\scrC\rp}
\]
and
\[\CAlg(\scrC)^{\mathrm{cond}}\coloneqq \mathrm{Sh}(\mathrm{ProFin}^{\mathrm{light}},\CAlg(\scrC))\simeq\Fun^{\mathrm{R}}{\lp(\CondAni)^{\op},\CAlg(\scrC)\rp},
\]
where the decoration R denotes those functors which preserve limits. So we can define a natural functor
\begin{equation*}
    \smash{\widetilde{\phi}} \colon \CAlg(\scrC)^{\mathrm{cond}} \longrightarrow \ \CAlg(\Fun(\Prof^{\mathrm{light}},\scrC))
\end{equation*}
as the composition
\[
\begin{tikzpicture}[scale=0.75]
\node (a) at (-4,2){$\RFun{\lp(\CondAni)^{\op},\CAlg(\scrC)\rp}$};
\node (b) at (4,2){$\Fun{\lp(\CondAni)^{\op},\CAlg(\scrC)\rp}$};
\node (c) at (-4,0){$\CAlg(\scrC)^{\mathrm{cond}}$};
\node (d) at (4,0){$\CAlg\lp\Fun{\lp(\CondAni)^{\op},\scrC\rp}\rp.$};
\draw[right hook->,font=\scriptsize] (a) to node [above]{}(b);
\draw[->,font=\scriptsize] (c) to node [rotate=90,above]{$\simeq$}(a);
\draw[->,font=\scriptsize] (b) to node [rotate=-90,above]{$\simeq$}(d);
\draw[->,font=\scriptsize] (c) to node [above]{$\smash{\widetilde{\phi}}$}(d);
\end{tikzpicture}
\]
In the above diagram, the vertical map on the right is an equivalence by \cite[Remark 2.1.3.4]{HA}, and the top horizontal arrow is a fully faithful embedding. In particular it follows that $\smash{\widetilde{\phi}}$ is a fully faithful embedding as well.

To prove that its essential image is the whole $\CAlg(\scrC^{\mathrm{cond}})$, it is enough to show the following: given a commutative algebra object $F$ in $\Fun{\lp(\CondAni)^{\op},\scrC\rp}$ such that the underlying functor $(\CondAni)^{\op}\to\scrC$ commutes with limits, then $F$ is in the image of $\smash{\widetilde{\phi}}$. But under the equivalence
\[
\Fun{\lp(\CondAni)^{\op},\CAlg(\scrC)\rp}\simeq\CAlg\lp\Fun{\lp(\CondAni)^{\op},\scrC\rp}\rp,
\]
the forgetful functor $\CAlg\lp\Fun{\lp(\CondAni)^{\op},\scrC\rp}\rp\to\Fun{\lp(\CondAni)^{\op},\scrC\rp}$ corresponds to the functor
\[
\Fun{\lp(\CondAni)^{\op},\CAlg(\scrC)\rp}\xrightarrow{\oblv_{\CAlg}\circ-}\Fun{\lp(\CondAni)^{\op},\scrC\rp}.
\]Since the forgetful functor $\CAlg(\scrC)\to\scrC$ is conservative and commutes with limits (\cite[Corollary 3.2.2.4 and Lemma 3.2.2.6]{HA},) our claim follows tautologically.
\end{proof}

\begin{remark}
Let $\scrC$ be a presentably symmetric monoidal $\infinity$-category, and let $\calA$ be a condensed commutative algebra in $\scrC$. In virtue of \cref{prop:commutative_algebras_condensed} $A$ corresponds to a commutative algebra in the presentably symmetric monoidal $\infinity$-category $\scrC^{\mathrm{cond}}$. In particular, \cite[$\S4.8$]{HA} yields a well-defined functor
\[
\Mod_{(-)}(\scrC^{\mathrm{cond}})\colon\CAlg(\scrC)^{\mathrm{cond}}\simeq\CAlg(\scrC^{\mathrm{cond}})\longrightarrow\CAlg(\mathrm{LPr}_{\scrC^{\mathrm{cond}}})
\]
sending a condensed $\Einf$-algebra in $\scrC$ to its (symmetric monoidal) derived $\infinity$-category of modules inside the $\infinity$-category of condensed objects of $\scrC$. When $\scrC=\Mod_{\mathbb{Z}}$, this procedure yields a stable derived $\infinity$-category. In this case, given a condensed $\Einf$-ring $\calA$, we shall denote the $\infinity$-category $\Mod_{\calA}(\Mod^{\mathrm{cond}}_{\mathbb{Z}})$ as $\Mod_{\calA}$.
\end{remark}

\subsection{Analytic rings and stacks}
In this subsection we review the six-functor formalism for analytic rings as defined by Clausen-Scholze in \cite{Clausen_Scholze_lectures}, and we show how to extend it to the context of analytic stacks. A nice overview of the following definitions and results can also be found in \cite{kunneth6ff}.

First, let us recall the following.
\begin{definition}
\label{def:analytic_ring}
    An analytic ring is a pair $A=(\calA,\scrD(A))$, where $\calA$ is a condensed $\Einf$-algebra, and $\scrD(A)$ is a full subcategory of $\Mod^{\mathrm{cond}}_{\calA}$ enjoying the following properties. 
    \begin{defenum}
    \item The $\infinity$-category $\scrD(A)$ contains the object $\calA$.
        \item \label{def:completion}Let $\iota_{A}\colon \scrD(A) \hookrightarrow \Mod_{ \calA}$ denote the natural inclusion functor. Then $\iota_A$ commutes with all limits and colimits, and it has a left adjoint (called the \emph{completion along $A$} functor) that we denote with $-\otimes_{\calA}A$. The composition $\iota_A\circ(-\otimes_{\calA}A)$ sends connective objects to connective objects.
        \item For every object $M$ in $\scrD(A)$ and every $N$ in $\Mod^{\mathrm{cond}}_{\Z}$, the internal mapping object $\Homin_{\Z}(N,M)$ lies in $\scrD(A)$. 
    \end{defenum}
A morphism between analytic rings $f\colon A\to B$ is a morphism of the underlying condensed $\Einf$-rings $\cf\colon\calA\to\mathcal{B}$ such that the induced forgetful functor $\Mod_{\mathcal{B}}\to\Mod_{\calA}$ sends $\scrD(B)$ to $\scrD(A)$. Analytic rings are naturally gathered in a $\infinity$-category $\Anring$. 
\end{definition}
\begin{remark}
   At a heuristic level, objects in $\Mod^{\mathrm{cond}}_{\calA}$
  play the role of topological $\calA$-modules. Following this analogy, objects in $\scrD(A)$ correspond to those topological $\calA$-modules which are complete. For this reason, this $\infinity$-category was originally denoted as $\scrD(\calA,\mathcal{M})$ in \cite{analytic}, where $\mathcal{M}$ should suggest a space of measures. This justifies the name of the left adjoint appearing in \cref{def:completion}.
\end{remark}
\begin{remark}
    By definition, $\iota_{A}$ is a localization functor. In particular, since $\Mod_{\calA}$ is presentable and stable (\cref{remark:condensed_stable}), it follows that $\scrD(A)$ is presentable and stable as well. Moreover, it admits a symmetric monoidal structure given by completing the natural tensor product of $\Mod_{\calA}$ along $A$. We denote such tensor product as $\otimes_{A}$.
\end{remark}
\begin{exmp}\
\begin{exmpenum}
\item\label{exmp:trivial_ring} As explained in \cite[Lecture 19]{Clausen_Scholze_lectures}, to every $\Einf$-$\Z$-algebra $R$ one can associate an analytic ring $R^{\mathrm{triv}}$, in the following way. Via the functor
    \[
    \CAlg_{\Z}\overset{\mathrm{const}}{\longrightarrow}\Fun{\lp(\Prof^{\mathrm{light}})^{\op},\CAlg_{\Z}\rp}\overset{(-)^+}{\longrightarrow}\CAlg^{\mathrm{cond}}_{\Z}
    \]
    we can associate a discrete condensed ring $\underline{R}$ to any $\Einf$-$\Z$-algebra $R$. So the analytic ring $R^{\triv}$ is explicitly given by the pair $(\underline{R},\Mod_{\underline{R}})$. 
    \item\label{exmp:rigid} Another source of analytic rings comes from rigid analytic geometry. As explained in \cite[Proposition~3.34]{andreychev2021pseudocoherent}, there is a fully faithful functor from the category of affinoid rings (and more generally from the category of complete Huber pairs) to the category of analytic rings. By abuse of notation we continue to denote by $A$ the analytic ring associated to an affinoid ring $A$, and we denote with $\scrD(A)$ the associated $\infinity$-category of "complete” modules. 
\end{exmpenum}
\end{exmp}

We are interested in (derived) analytic stacks, which are locally modeled by analytic rings. As customary in algebraic geometry, we shall denote the opposite $\infinity$-category of $\Anring$ as $\mathrm{Aff}$, and call it the $\infinity$-category of affine analytic stacks. Given an analytic ring $A$, its associated affine stack will be denoted as $\Ansp(A)$.
\begin{parag}
Let $f\colon \Ansp(A)\to\Ansp(B)$ be a morphism of affine analytic rings. Let us denote by $f_\ast$ the restriction of the forgetful functor $\Mod_{\mathcal{B}}\to\Mod_{\calA}$ to the sub-$\infinity$-category $\scrD(A)$: by the definition of morphism of analytic rings, this is a functor $f_\ast\colon\scrD(B)\to\scrD(A)$, which we can write as a composition of right adjoints
\[
\scrD(B)\overset{\iota_B}{\longhookrightarrow}\Mod_{\mathcal{B}}\longrightarrow\Mod_{\calA}\overset{f^{\mathrm{R}}}{\longrightarrow}\scrD(A).
\]
The last functor is the right adjoint of $\iota_A$, which exists because $\scrD(A)$ is a presentable $\infinity$-category closed under all colimits inside $\Mod_{\calA}$. The fact that the functor acts as the identity on the image of $\scrD(B)$ follows from the facts that the forgetful functor preserves complete modules and $f^{\mathrm{R}}$ is a colocalization functor.

It follows that the functor $f_\ast$ admits a left adjoint $f^{\ast}\colon \scrD(B) \to \scrD(A)$ given by the composition 
\[
\scrD(B)\overset{\iota_B}{\longhookrightarrow}\Mod_{\mathcal{B}}\overset{-\otimes_{\mathcal{B}}\calA}{\longrightarrow}\Mod_{\calA}\overset{-\otimes_{\calA}A}{\longrightarrow}\scrD(A).
\]
\end{parag}
In virtue of \cite[Proposition 4.11]{kunneth6ff} the assignment
\[
A\mapsto\scrD(A),\quad f\mapsto f^\ast
\]
can be promoted to a symmetric monoidal functor
\begin{equation}\label{eq: 6 functor formalism}
    \scrD\colon \mathrm{AnRing}\simeq\mathrm{Aff}^{\op} \longrightarrow \CAlg(\LPr_{\mathrm{st}}),
\end{equation}
which generalizes the usual functor sending a commutative ring to its stable $\infinity$-category of derived modules. 
\begin{definition}
    Let $f\colon\Ansp(A)\to\Ansp(B)$ be a map of affine analytic stacks.
    \begin{defenum}
        \item\label{def:proper} We say that $f$ is \emph{proper} if $f^\ast$ induces an equivalence $\scrD(A)\simeq \Mod_{A}(\scrD(B))$. This means that $B$ has the analytic ring structure induced by the one of $A$.
        \item \label{def:open immersion}We say that $f$ is an \emph{open immersion} if $f^{\ast}$ has a left adjoint $f_{!}$ that satisfies the projection formula. This means that there is a canonical equivalence of functors
        \[
        f_!(-)\otimes_B(-)\simeq f_!((-)\otimes_Af^\ast(-))
        \]
        from $\scrD(A)\otimes\scrD(B)$ to $\scrD(B).$
        \item We say that $f$ is \emph{!-able} if can be written as a composition $f \simeq p \circ j$, where $j$ is an open immersion and $p$ is proper. 
    \end{defenum}
\end{definition}
\begin{remark}
If $f$ is a !-able map with given decomposition $f \simeq p \circ j$, then we can define a functor $f_{!}\colon\scrD(A)\to\scrD(B)$ to be $p_{\ast} \circ j_{!}$. It is indeed a left adjoint: this follows from the fact that $j_!$ is a left adjoint by definition, while $p_\ast$ (which is normally only a right adjoint) is also a \textit{left} adjoint when $p$ is proper. Indeed, in this case, $p_\ast$ identifies with the forgetful functor $\Mod_A(\scrD(B))\to\scrD(A)$, which preserves all limits and colimits. In particular, $f_!$ admits a right adjoint $f^!$.
\end{remark}
The following theorem appears in the video series on analytic stacks by Clausen and Scholze. For completeness, we include a proof here, drawing on results from the existing literature.
\begin{theorem}[\cite{Clausen_Scholze_lectures}]
\label{thm:clausen_scholze_6ff}
Let $E$ denote the class of !-able morphisms, let $I$ denote the class of open immersions, and let $P$ denote the class of proper maps in $\mathrm{Aff}$.
\begin{thmenum}
\item\label{thm:geom_setup}The pair $(\mathrm{Aff},E)$ defines a geometric setup, in the sense of \cite[Definition 2.2.1]{heyer6ff2024}.
\item \label{thm:suitable_decomp}The pair $(P,I)$ is a suitable decomposition of $E$, in the sense of \cite[Definition 3.3.2]{heyer6ff2024}
\item \label{thm:6ff_strongly_monoidal}The functor \eqref{eq: 6 functor formalism} is a sheaf for the universal !-topology, and it is one of the operations in a symmetric monoidal six-functor formalism over $(\mathrm{Aff},E)$.
\end{thmenum}
\end{theorem}
\begin{proof}
    \cref{thm:geom_setup} and \cref{thm:suitable_decomp} are proved in \cite[Lemmas 4.6 and 4.7]{kunneth6ff}. For \cref{thm:6ff_strongly_monoidal}, we note that the fact that $\scrD(-)$ is a sheaf for the universal !-able topology follows directly from the definition of universal !-able topology (\cref{def:!-able}), while the fact that $\scrD(-)$ can be promoted to a six-functor formalism follows from \cite[Proposition 3.3.3]{heyer6ff2024}. The fact that $\scrD(-)$ is a symmetric monoidal functor follows from \cite[Proposition 4.11]{kunneth6ff}.
\end{proof}
\begin{exmp}
Let $R$ be an affinoid ring, seen as an analytic ring as in \cref{exmp:rigid}. Then, \cite{Clausen_Scholze_lectures} shows that the functor from affinoid rings to analytic rings sends analytic open immersions of affinoid rings to universal !-able maps of analytic rings. This makes it possible to regard rigid analytic varieties as analytic stacks. Again, when $X$ is a rigid analytic variety we shall abuse our notations and still denote its associated analytic stack by $X$.
\end{exmp}

\subsection{Proper objects and descent}
In this subsection we will briefly describe the notion of $\scrD$-proper maps, as defined in \cite{heyer6ff2024} (where they are called prim maps). For the reader's convenience, we also recall some properties and results borrowed from \cite{six-functors-scholze} (where they are called "cohomologically proper”) and \cite{camargo2024analytic} (where they are called "co-smooth”) that we will need in the rest of the paper.

We refer the reader to the above references for a complete treatment of this topic and all the relevant proofs.
\begin{notation}
Let $\scrC$ be a $\infinity$-category with finite limits. Fix a class of morphisms $E$ such that the pair $(\scrC, E)$ is a geometric setup, and let $(\scrD^\ast,\scrD_!)$ be a six-functor formalism over $(\scrC,E)$. Given a map $f \colon X \to S$, we denote with $p_{1,2}\colon X\times_{S}X \to X$ the two projection maps, and with $\Delta\colon X \to X\times_{S}X$ the diagonal.

If $P$ is an object in $\scrD(X)$ we denote with $\scrP_{f}(P)$ the object of $\scrD(X)$ given by $p_{2,\ast}\Homin_{\scrD(X\times_{S}X)}(p_{1}^{\ast}(P), \Delta_{!}(\boldone_{\scrD(X)}))$.
\end{notation}
\begin{definition}[{\cite[Proposition 3.1.11]{camargo2024analytic},\cite[Lemma 4.4.6]{heyer6ff2024}}]
     Let $f \colon X \to S$ be a map lying in $E$.
     \begin{defenum}
        \item We say that an object $P \in \scrD(X)$ is \textit{$f$-proper} if the natural map
     \[
     f_!(P\otimes\scrP_f(P))\overset{f_\ast}{\longrightarrow}f_\ast\Homin_{\scrD(X)}(P,P)
     \]
is an equivalence. 
\item We say that $f$ is \textit{$\scrD$-proper} if $\boldone_{\scrD(X)}$ is $f$-proper and $\scrP_{f}(\boldone_{\scrD(X)})$ is invertible with respect to the symmetric monoidal structure of $\scrD(X)$.
     \end{defenum}
\end{definition}
\begin{remark}
Note that in \cite{heyer6ff2024} the terminology $f$-prim is used for what we call $f$-proper. In that reference, the term $f$-proper is reserved for a stronger condition.
\end{remark}

The following result is the key tool that allows us to establish whether a morphism constitutes a universal !-able cover of the target.
\begin{theorem}[{\cite[Theorem 6.19]{six-functors-scholze}}]\label{thm: proper descent}
Let $f\colon X \to S$ be a morphism in $E$. If $f$ is descendable (\cref{def: descendable}) and $\boldone_{\scrD(X)}$ is $f$-proper, then $\scrD(-)$ satisfies universal \textasteriskcentered- and $!$-descent along $f$.
\end{theorem}
It follows that in order to understand whether a map is a universal !-able cover, we have to be able to detect whether an object is proper. The next proposition allows us to check this condition \textasteriskcentered-locally on the target.
\begin{proposition}[{\cite[Lemma 4.4.8]{heyer6ff2024}}]
    Let $f \colon X \to S$ be a morphism in $E$. Suppose that $S$ admits a universal \textasteriskcentered-cover $\{g_i\colon T_{i} \to S\}$, and denote with $f_{i}\colon X \times_{S} T_{i} \to T_{i}$ and $g_{i,X}\colon X \times_{S} T_{i} \to X$ the obvious projection maps. An object $P$ of $\scrD(X)$ is $f$-proper if the objects $g_{i,X}^{\ast}(P)$ are $f_{i}$-proper for all $i$. 
\end{proposition}
In order to show that $\boldone_{\scrD(X)}$ is $f$-proper and apply \Cref{thm: proper descent}, in many cases we can use the following weaker definition of properness.
\begin{definition}[{\cite[Definition 3.1.19]{camargo2024analytic}}]
\label{def:weakly_cohom_proper}
    Let $f\colon X\to S$ be a morphism of $\scrC$. We say that $f\colon X \to S$ is \emph{weakly cohomologically proper} if the following conditions hold.
    \begin{defenum}
        \item The object $\boldone_{\scrD(X)}$ is $\Delta_{f}$-proper.
        \item There is an equivalence $\scrP_{\Delta_{f}}(\boldone_{\scrD(X)}) \overset{\simeq}{\to} \boldone_{\scrD(X)}$.
        \item The object $\boldone_{\scrD(X)}$ is $f$-proper.
    \end{defenum}
\end{definition}
\begin{remark}\label{rmk: cohomologically proper maps affine}
\cref{def:weakly_cohom_proper} bears an important advantage over \cref{def:proper}. Indeed, when the six-functor formalism we are considering admits a suitable decomposition $(I,P)$ then all proper maps are also weakly cohomologically proper, as proved in \cite[Lemma 3.1.21]{camargo2024analytic}; in particular, in this case, given any proper map $f\colon X\to S$ the object $\boldone_{\scrD(X)}$ is $f$-proper. Thanks to \cref{thm:suitable_decomp}, this discussion applies to the case of the six-functor formalism over $\mathrm{AnRing}$ as well.
\end{remark}

\section{1-affineness in analytic geometry}
\label{sec:main_applications}
Finally, we are ready to prove our main $1$-affineness results in the context of rigid and analytic geometry, highlighting a class of analytic stacks and rigid varieties which are $1$-affine. While the final results that we obtain are similar in both the analytic and rigid cases, the strategy of the proofs and the reason why certain objects are $1$-affine are different, so we split our presentation in two parts.
\subsection{\texorpdfstring{$1$}{1}-affineness for derived categories of analytic stacks}
\label{sec:1_affineness_condensed}
Consider the $\infinity$-category $\Anring$ of analytic rings (in the sense of \cref{def:analytic_ring}), and consider the sheaf of presentably $\scrD(\Z)$-linear symmetric monoial $\infinity$-categories $\scrD(-)\colon\Anring\to\CAlg(\LPr_{\scrD(\Z)})$. Endow $\Anring$ with the universal !-able topology. Then, \cref{thm:6ff_strongly_monoidal} guarantees that this piece of data fits in the setting of \cref{assumption:six_functor_case}. In particular, \cref{prop:six_functor_case} tells us the following.
\begin{proposition}\label{}
    The functor
\[
\Mod_{\scrD(-)}\colon\Anring_{\Z}\longrightarrow\LPr_{\scrD(\Z)}
\]
is a sheaf for the universal !-able topology.
\end{proposition}
 Our goal is to single out a class of analytic stacks that are $1$-affine. This problem can be seen as a particular instance of \Cref{question:main}. Indeed, in this case \cref{thm:main 2-descent} immediately specializes to the following theorem.
\begin{theorem}\label{thm: main theorem}
    Let $X$ be an analytic stack admitting an affine !-cover $\{\Ansp(A_{i}) \to X \}_{i \in I}$ such that all the fiber products $\Ansp(A_{i_{1}})\times_{X}\dots \times_{X}\Ansp(A_{i_{n}})$ are affine. Then $X$ is $1$-affine. 
\end{theorem}

\begin{corollary}\label{cor: 2 descent}
    Let $X$ be an analytic stack admitting an affine universal !-cover $\{\Ansp(A_{i}) \to X \}_{i \in I}$ such that all the fiber products $\Ansp(A_{i_{1}})\times_{X}\dots \times_{X}\Ansp(A_{i_{n}})$ are affine. Then the functor
    \begin{equation}
        \Mod_{\scrD(-)} \colon \AnStacks \to \mathrm{Cat}^{\mathrm{rex}}
    \end{equation}
 satisfies descent along the cover $\{\Ansp(A_{i}) \to X \}_{i \in I}$.
\end{corollary}
\begin{proof}
    This follows using \Cref{thm: main theorem} since by definition $\scrD^{\mathrm{cat}}(-)\coloneqq\Mod_{\scrD(-)}$ satisfies descent. 
\end{proof}

\cref{thm: main theorem} in particular subsumes the case when $X$ is an analytic stack with affine diagonal, which has also been independently proved in \cite{kunneth6ff}. However, our formulation allows to deduce the $1$-affineness of more analytic stacks, for which we can prove that fiber products of affines are affine only in the case of a specific affine universal !-able cover. This is the case of the analytic Betti stack. Similarly \Cref{cor: 2 descent} can also be deduced by \cite[Lemma 3.32]{kunneth6ff}
\begin{parag}
Let $\CondAni$ be the $\infinity$-category of condensed anima (that is, condensed objects in the $\infinity$-category of anima). Inside this $\infinity$-category, there are two classes of objects coming from topology.
\begin{enumerate}
    \item \textit{Homotopy types} or \textit{anima}: they are the objects in the essential image of the constant functor $\Ani\to\CondAni$.
    \item \textit{Topological spaces} (\textit{not} up to homotopy): all topological spaces can be understood as condensed sets (\cite[Example 1.5]{clausen2019lectures}). The functor $\mathrm{Top}\to\mathrm{Set}^{\mathrm{cond}}$ is always faithful, and moreover fully faithful on the subcategory of $\kappa$-compactly generated topological spaces (\cite[Proposition 1.7]{clausen2019lectures}). So, post-composing with the limit-preserving inclusion $\iota_0\colon \mathrm{Set}\subseteq\Ani$ we obtain a faithful functor $\mathrm{Top}\to\CondAni$.
\end{enumerate}
If we denote with $|-|$ the functor sending a topological space to its underlying anima, and with $\underline{(-)}$ the functor from topological spaces to condensed sets as described above, we obtain the diagram
\begin{equation}\label{eq: Top into CondAni}
\begin{tikzpicture}[scale=0.75,baseline=(current  bounding  box.center)]
\node (a) at (0,2) {$\CondAni$};
\node (b) at (-3,0){$\mathrm{Set}^{\mathrm{cond}}$};
\node (c) at (3,0){$\Ani$};
\node (d) at (0,-2){$\mathrm{Top}.$};
\draw[right hook->,font=\scriptsize] (b) to node[above left]{$\iota_0$} (a);
\draw[->,font=\scriptsize] (c) to node[above right]{$\mathrm{const}$} (a);
\draw[->,font=\scriptsize] (d) to node[below left]{$\underline{(-)}$} (b);
\draw[->,font=\scriptsize] (d) to node[below right]{$|-|$} (c);
\end{tikzpicture}
\end{equation}
which however is \textit{not} commutative.
\end{parag}
In particular we can see two different realizations of a topological space in the category $\CondAni$ . Using the right hand side of the diagram \eqref{eq: Top into CondAni} we can assign to a topological space the sheafification of the constant presheaf with values in its underlying homotopy type. On the other side, we can assign to a topological space a $0$-truncated object in the $\infinity$-topos of condensed anima (in the sense of \cite[Definition 5.5.6.1]{htt}). This means that we have a collection of homotopy \textit{condensed} groups, and all the higher homtopy condensed groups are trivial. 

In this section we will view topological spaces as objects in condensed anima using the composition on the left hand side of the diagram \eqref{eq: Top into CondAni}. This functor is actually very well behaved. 
\begin{lemma}\label{lem: top space as condnesed anima}
    The functor
    \[
    \mathrm{Top}\overset{\underline{(-)}}{\longrightarrow}\mathrm{Set}^{\mathrm{cond}}\overset{\iota_0}{\longhookrightarrow}\CondAni
    \]
appearing in the left hand side of the diagram \eqref{eq: Top into CondAni} commutes with all limits. Moreover, it is fully faithful when restricted to $\kappa$-compactly generated topological spaces.
\end{lemma}
\begin{remark}
Two nice classes of topological spaces to which \cref{lem: top space as condnesed anima} applies are the class of all first-countable topological spaces and the class of all metrizable topological spaces (\cite[Remark 1.6]{clausen2019lectures}).
\end{remark}
\begin{proof}[{Proof of \cref{lem: top space as condnesed anima}}]
    The first assertions follows from the fact that both the functor $\mathrm{Top} \to \mathrm{Set}^{\mathrm{cond}}$ and the functor $\mathrm{Set}^{\mathrm{cond}} \to \CondAni$ are right adjoints. For the first functor, this is proved in \cite[Proposition 1.7]{clausen2019lectures}; for the latter, just notice that limits in $\infinity$-categories of sheaves are computed as in the $\infinity$-category of presheaves (hence, they are computed point-wise), and $\mathrm{Set}\simeq\Ani_{\leqslant0}\subseteq\Ani$ is closed under all limits (\cite[Proposition 5.5.6.5]{htt}).
    
    The second statement follows again from the fact that the first functor is fully faithful when restricted to $\kappa$-compactly generated topological spaces, and the second is always fully faithful.
\end{proof}
The above observations allow us to define the \textit{condensed Betti stack of a topological space}, following \cite{geometrizationreallanglands}. We start by defining this functor on light pro-finite sets, and then we extend it to every condensed anima (hence, to every topological space in virtue of \cref{lem: top space as condnesed anima}).
\begin{construction}
Let $(-)_{\Betti}$ denote the functor
\begin{equation}\label{eq: betti light pro-finite}
\Prof^{\mathrm{light}} \longrightarrow \AnStacks 
\end{equation}
defined by sending a pro-finite set $T$ to the affine analytic stack $\Ansp(\mathcal{C}^{0}(T,\Z)^{\mathrm{triv}})$. Here, $\mathcal{C}^{0}(T,\Z)^{\triv}$ denotes the trivial analytic ring on the $\Einf$-ring of continuous functions from $T$ (equipped with its inverse limit topology) to $\Z$, as described in \cref{exmp:trivial_ring}. In other words, the class of complete condensed modules over $\scrC^{0}(T,\Z)$ coincides with the whole $\Mod_{\underline{\mathcal{C}^{0}(T,\Z)}}.$
\end{construction}
\begin{proposition}[{\cite[Lecture 20]{Clausen_Scholze_lectures}}]\label{prop: affine betti preserves colimits and finite limits}
    The functor \eqref{eq: betti light pro-finite} preserves finite limits and sends hypercovers of pro-finite sets to universal !-covers. 
\end{proposition}
\begin{proof}[Sketch of proof]
Since each pro-finite set is a inverse limit of finite sets one reduces to prove the claim for finite sets, for which the proofs are easy.
\end{proof}
\begin{corollary}\label{prop: betti esteso}
The functor \eqref{eq: betti light pro-finite} admits a unique left exact and colimit-preserving extension to a functor $\CondAni\to\AnStacks$. We still denote such extension as $(-)_{\Betti}$, and call this functor the \emph{extended Betti stack functor}.
\end{corollary}
\begin{proof}
    This is a direct consequence of \Cref{prop: affine betti preserves colimits and finite limits} applying \cite[Proposition 6.2.3.20]{htt}. 
\end{proof}
Let now $X$ be a topological space and let $T$ be  a light pro-finite set such that there exists a surjection $f \colon T \to X$ (here we are viewing light pro-finite sets inside $\mathrm{Top}$ in virtue of \cref{remark:stone_duality}). For any $i\geqslant 0$, let us denote with $T_{i}$ the $i$-th object in the simplicial hypercover associated to \v{C}ech nerve of $f$. Since $(-)_{\Betti}$ commutes with colimits and finite limits, we can describe $X_{\Betti}$ as the geometric realization of the simplicial diagram
\[
(T_1)_{\Betti}\stack{3}(T_2)_{\Betti}\stack{5}(T_3)_{\Betti}\stack{7}\cdots.
\]
In general, the $(T_i)_{\Betti}$'s are not affine analytic stacks. The goal of the following lemmas is to show that, at least when $X$ is a finite dimensional, metrizable and compact topological space, then the above simplicial diagram arises as the \v{C}ech nerve of a universal !-cover of analytic stacks, and at each step such simplicial diagram consists of affine analytic stacks. In particular, this applies to the case when $X$ is a compact topological manifold.
\begin{lemma}\label{lem: pullback of light pro-finite}
    Let $X$ be a Hausdorff topological space and let $T_{1}$ and $T_{2}$ two light pro-finite sets mapping into $X$. Then the pullback $T_{1}\times_{X}T_{2}$ is still a light pro-finite set.
\end{lemma}
\begin{proof}
Let $f_i\colon T_i\to X$ be the two maps in the statement. We have a natural map $\langle f_1,f_2\rangle\colon T_1\times T_2\to X\times X$. Taking the pullback of this map along the diagonal $\Delta\colon X\to X\times X$ yields a natural map
\[
\alpha\colon (T_1\times T_2)\times_{X\times X}X\longrightarrow T_1\times T_2.
\]
Since $X$ is Hausdorff, the diagonal $X \to X \times X$ is a closed embedding, and so is $\alpha$. However, the source of $\alpha$ is homeomorphic to $T_1\times_XT_2$. In particular, it is a closed subset of the compact metrizable totally disconnected Hausdorff topological space $T_1\times T_2$, so it is compact, metrizable, totally disconnected and Hausdorff itself. 
\end{proof}
\begin{lemma}
    Let $X$ be a metrizable and compact topological space. Then there exist a light pro-finite set $K$ and a surjection $f \colon K \to X$. Its \v{C}ech yields a hypercover in $\CondAni$, which at each step consists only of light pro-finite sets.\\    
    Furthermore $f_{\Betti}\colon K_{\Betti} \to X_{\Betti}$ is a universal \textasteriskcentered-cover. 
\end{lemma}
\begin{proof}
    Since $X$ is a metrizable compact topological space, we have a continuous surjection form the Cantor set $K \coloneq\{0,1\}^{\N}$ to $X$. Let $f$ denote such cover map.
    
    Using \Cref{lem: pullback of light pro-finite}, we have that every pullback $K\times_{X} \dots \times_{X}K$ is again a light pro-finite set. Now the \v{C}ech nerve of $f$ is the required hypercover of $X$ consisting of light pro-finite sets. To prove that $f_{\Betti}$ is a \textasteriskcentered-cover we need to show that the functor $\scrD(-)$ satisfies descent along $f^\ast_{\Betti}$ and along every fiber product of $f^\ast_{\Betti}$. By construction $f_{\Betti}$ is an effective epimorphism, so it naturally satisfies \textasteriskcentered-descent. \textit{Universal} \textasteriskcentered-descent follows from the fact that effective epimorphisms are stable under pullbacks(\cite[Proposition 6.2.3.15]{htt}).
\end{proof}
\begin{lemma}\label{prop: !-cover betti stacks}
    Let $X$ be a finite dimensional, metrizable and compact topological space. Then there is a universal !-cover of $f_{\Betti} \colon K_{\Betti} \to X_{\Betti}$ where $K$ is a light pro-finite sets and the \v{C}ech nerve of $f_{\Betti}$ is made of affine analytic stacks. 
\end{lemma}
\begin{proof}
    Arguing as in the above proof we can find a light pro-finite set $K$ and a surjection $K \to X$ such that the associated \v{C}ech nerve is made of light pro-finite sets. In virtue of \cite[Proposition II.1.1]{geometrizationreallanglands} the object $f_{\Betti}(\boldone)$ is descendable, so if $\boldone_{K_{\Betti}}$ was also $f_{\Betti}$-proper we could apply \cite[Lemma 4.7.4]{heyer6ff2024} to deduce that $f_{\Betti}$ yields a universal !-cover. 
    
    Using \cite[Lemma 4.4.8]{heyer6ff2024} and \Cref{lem: pullback of light pro-finite} we can assume that $f$ is a map of light pro-finite sets. In this case, $f_{\Betti}$ is actually weakly cohomologicaly proper in the sense of \cref{def:weakly_cohom_proper}. Indeed, the analytic ring structure on $\mathcal{C}^0(T,\Z)$ is always induced from the one of $\Z^{\triv}$ for all light pro-finite sets $T$, and so we simply apply \Cref{rmk: cohomologically proper maps affine} to get our desired assertion. 
    
    The fact that the \v{C}ech nerve of $f_{\Betti}$ is made of affine analytic stacks follows directly from \Cref{prop: betti esteso} and \Cref{lem: pullback of light pro-finite}. 
\end{proof}

\begin{theorem}\label{thm: Betti is 1-affine}
Let $X$ be a finite dimensional, metrizable, compact Hausdorff space. Then the analytic stack $X_{\Betti}$ is 1-affine. Morover the functor
\begin{equation*}
    \Mod_{\scrD(-)}\colon \AnStacks \longrightarrow \catrex
\end{equation*}
satisfies descent along the universal !-cover described in \Cref{prop: !-cover betti stacks}.
\end{theorem}
\begin{proof}
Applying \Cref{prop: !-cover betti stacks}, we see that the hypotheses of \cref{thm: main theorem} are evidently satisfied, so $X_{\Betti}$ is indeed $1$-affine. 
      The last part follows form \Cref{cor: 2 descent}.
\end{proof}
\begin{remark}
In derived algebraic geometry, one can define another Betti stack of a topological space $X$. This is the sheafification of the constant prestack with values in the underlying anima $|X|$. 

In \cite{PPS_2}, it is proved that this "homotopical" Betti stack is never $1$-affine when the second homotopy group of $X$ is non-trivial. However, there are plenty metrizable and compact finite dimensional topological manifolds which admit non-trivial $\pi_2$: one example is provided by the real projective space $\RR\PP^2$. Thus, \cref{thm: Betti is 1-affine} shows a curious discrepancy between the analytic and the homotopical Betti stacks. 
\end{remark}
We collect now some consequences of \cref{thm: main theorem}. First, suppose that $X$ is a compact complex manifold. Using the analytic Riemman--Hilbert correspondence as stated in \cite[Theorem II.3.1]{geometrizationreallanglands} and our \Cref{thm: Betti is 1-affine} we can immediately prove the following result, that can be understood as the "analytic” version of a theorem of Gaitsgory (see \cite[Theorem 2.6.3]{1affineness}).
\begin{corollary}
\label{cor:de_rham_1_affine}
    Let $X$ be a compact complex manifold, Then the analytic de Rham stack, $X^{\mathrm{an}}_{\mathrm{dR}}$ as defined in \cite{geometrizationreallanglands} is $1$-affine.
\end{corollary}
Another immediate consequence of our result is $2$-descent for derived categories of sheaves over well-behaved topological spaces.
\begin{corollary}\label{cor: two descent condensed coef}
    Let $X$ be a finite dimensional, metrizable, compact Hausdorff space. Fix a surjection $f\colon T \to X$ from a light pro-finite set $T$. Then the functor $\mathrm{Sh}(-,\scrD(\Z^{\triv}))$ sending a topological space $X$ to its derived $\infinity$-category of sheaves with values in $\scrD(\Z^{\triv})\simeq\Mod_{\underline{Z}}$ satisfies $2$-descent along $f$. In other words, the functor
    \[
    \LPr_{\mathrm{Sh}(-,\scrD(\Z^{\triv}))}\colon\mathrm{Top}\longrightarrow\catrex
    \]
    satisfies descent along $f$. 
\end{corollary}
\begin{proof}
    This follows from \Cref{thm: Betti is 1-affine}, once we observe that for any topological space $X$ the $\infinity$-category $\scrD(X_{\Betti})$ agrees with the derived $\infinity$-category $\mathrm{Sh}(X,\scrD(\Z^{\triv}))$ (\cite[Corollary II.1.2]{geometrizationreallanglands}).
\end{proof}
\begin{warning}
    Note that it is not reasonable to extend the above corollary to every compact Hausdorff space. Indeed, in \cite{haine2022} it was proved that for a general presentable $\infinity$-category $\scrC$ the functor
    \begin{equation*}
        \mathrm{Sh(-, \scrC)}\colon \mathrm{CHaus} \longrightarrow \widehat{\mathrm{Cat}}_{\scriptstyle\infty}
    \end{equation*}
    from compact Hausdorff spaces does not satisfy descent along a surjection $f\colon T\to X$ from a light pro-finite set. It does, however, satisfy descent \textit{ up to Postnikov completion} (\cite[Corollary 2.8]{haine2022}). In particular, for $\mathrm{Sh}(-,\scrC)$ to satisfy descent along $f$ it is necessary that $\mathrm{Sh}(X,\scrC)$ is itself Postnikov complete. This is the case for example when $X$ is a pro-finite set or is a CW-complex, see \cite[Examples 1.28 and 1.29]{haine2022}.
\end{warning}


\subsection{Nuclear categories in rigid geometry}
\label{sec:nuclear}
We want to provide an answer to \cref{question:main} in the setting of rigid analytic geometry. In order to do so, we need to turn to the setting of $\infinity$-categories of nuclear modules.

In this section we recall the relevant definitions, and we prove that this problem fits in the general framework of categorified sheaves described in \cref{sec:modular_case} (see \Cref{prop: nuclear categories and module categories}).
\begin{assumption}
For the remainder of this section, we will work with \textit{classical} affinoid algebras and rigid analytic varieties defined over a complete non-archimedean base field $\Bbbk$ equipped with the $I$-adic topology with respect to an ideal of definition $I$, with a pseudo-uniformizer $\varpi$ and ring of integers $\scrO_{\Bbbk}$. We always abuse our notations and identify an affinoid $\Bbbk$-algebra $A$ with its associated Huber pair $(A,A^0)$, where $A^0$ is the whole subring spanned by topologically nilpotent elements of $A$, and also with its associated analytic ring (\cref{exmp:rigid}). In particular, we shall still denote as $\Bbbk$ the associated condensed ring of the base field, with topology induced by its non-archimedean norm.\\ 

We will later consider rigid analytic varieties defined $\Bbbk$, in the sense of \cite{huber2013etale}. More explicitly, rigid analytic varieties are obtained by gluing Tate--Huber pairs $(A,A^0)$ of finite type over $(\Bbbk, \mathcal{O}_{\Bbbk})$ in the analytic topology. When a rigid analytic variety is the affinoid variety associated to an affinoid $\Bbbk$-algebra $A$, then we shall denote it as $\Spa(A)$. 
\end{assumption}


\begin{definition}\label{def:nuclear}
Let $\cC$ be a stable closed symmetric monoidal $\infty$-category. 
An object $X \in \cC$ is called \emph{nuclear} if, for every compact object $P \in \cC$, 
the canonical map
\[
 \Map_{\cC}\!\bigl(\boldone_{\cC},\, \Homin_{\scrC}(P,\boldone_{\cC}) \otimes_{\cC} X \bigr) 
  \;\longrightarrow\;
 \Map_{\cC}(P,X)
\]
is an equivalence of anima. 
The full subcategory of $\cC$ spanned by nuclear objects will be denoted as $\Nuc(\cC) \subseteq \cC$.
\end{definition}
\begin{notation}
When $\scrC$ is the derived $\infinity$-category $\scrD(A)$ of an analytic ring $A$, we shall simply write $\Nuc(A)$ instead of $\Nuc(\scrD(A)).$
\end{notation}
Given an analytic ring $A$, the $\infinity$-category $\Nuc(A)$ is a presentable full sub-$\infinity$-category of $\scrD(A)$ which is closed under colimits and tensor product (\cite[Theorem 8.6]{clausenlectures}). This equips $\Nuc(A)$ of a symmetric monoidal tensor product that we denote as $\widehat\otimes_A$. Moreover, the assignment $A\mapsto \Nuc(A)$ defines a sheaf for the analytic topology on affinoid $\Bbbk$-algebras (\cite[Theorem 5.42]{andreychev2021pseudocoherent}). Given a map of affinoid $\Bbbk$-algebras $f\colon A\to B$, the functor $\Nuc(A)\to\Nuc(B)$ is the restriction of the functor $\scrD(A)\to\scrD(B)$ to the full sub-$\infinity$-categories of nuclear objects (\cite[Theorem 2.9 (c)]{meyer2024derived}).

In this subsection we want to study a categorified version of the sheaf of nuclear modules. More precisely, we will study the functor
\begin{align}
\label{eq: categorifed nuc}
\begin{split}
\Afd_{\Bbbk}&\longrightarrow\catrex\\
A&\mapsto \LPr_{\Nuc(A)}.
\end{split}
\end{align}
We will prove that such functor is a sheaf for the analytic topology, and more generally for the fppf topology (\Cref{cor: Nuc satisy descent}). 

We start with the following.
\begin{proposition}[{\cite[Remark 2.13]{mikami2023fppf}}]\label{prop: affinoids are nuclear}
    Let $A$ be an affinoid $\Bbbk$-algebra. Then $A$ is nuclear over $\Bbbk$.
\end{proposition}
\cref{prop: affinoids are nuclear} tells us that given any affinoid $\Bbbk$-algebra $A$, we can view it as a $\Einf$-algebra object inside $\Nuc(\Bbbk)$. So we can consider two \textit{a priori} distinct symmetric monoidal $\infinity$-categories associated to an affinoid $\Bbbk$-algebra.
\begin{enumerate}
    \item The $\infinity$-category $\Nuc(A)$, with tensor product given by the restriction of the tensor product of $\scrD(A)$.
    \item The $\infinity$-category $\Mod_A(\Nuc(\Bbbk))$, with tensor product given by the usual relative tensor product over $A$.
\end{enumerate}
We now show that these two $\infinity$-categories are actually symmetric monoidally equivalent. 
\begin{proposition}\label{prop: nuclear categories and module categories}
    Let $A$ be an affinoid algebra over $\Bbbk$. Then there is an equivalence of $\infinity$- categories between $\Mod_{A}(\Nuc(\Bbbk))$ and $\Nuc(A)$. In particular $\Nuc(A)$ is a rigid $\infinity$-category.
\end{proposition}
\begin{proof}
    As explained in \cite[Corollary 2.7]{mikami2023fppf}, if $A$ is an  affinoid $\Bbbk$-algebra the functor $\scrD(A) \to \Mod_{A}(\scrD(\Bbbk))$ induces an equivalence between the associated nuclear sub-$\infinity$-categories. In particular, since $A$ is already nuclear in $\scrD(\Bbbk)$, we deduce that $\Nuc(A)$ is equivalent to $\Mod_{A}(\Nuc(\Bbbk))$. The last assertion follows from \cite[Remark 3.34]{anschutz2024}.
\end{proof}
\begin{remark}
Combining \cref{prop: affinoids are nuclear} with \cref{prop: nuclear categories and module categories}, we can see that for any map of affinoid $\Bbbk$-algebras $A\to B$ we can consider $B$ as an object in $\Nuc(A)$. Indeed, $A\to B$ is in particular a map of commutative algebras inside $\Nuc(\Bbbk)$, so $B$ inherits a natural $A$-module structure and so corresponds to an object in $\Mod_A(\Nuc(\Bbbk))$. This $\infinity$-category is the same as the $\infinity$-category $\Nuc(A)$, so our claim follows.
\end{remark}
\begin{proposition}\label{prop: pullback of nuclear categories as base change}
    Let $A \to B$ a map of affinoid algebras over $\Bbbk$. Under the identification $\Mod_A(\Nuc(\Bbbk))\simeq\Nuc(A)$ of \Cref{prop: nuclear categories and module categories}, the pullback functor $\Nuc(A)\to\Nuc(B)$ can be identified with the usual base-change functor $-\widehat\otimes_{A}B$.
\end{proposition}
\begin{proof}
Using the categorical Eilenberg--Watts theorem (\cite[Theorem 4.8.4.1]{HA}), we can identify the $\infinity$-category of $\Nuc(\Bbbk)$-linear, colimit-preserving functors from $\Nuc(A)\simeq\Mod_{A}(\Nuc(\Bbbk))$ to $\Nuc(B)\simeq\Mod_{B}(\Nuc(\Bbbk))$ with the $\infinity$-category of $(A,B)$-bimodules in $\Nuc(\Bbbk)$. More explicitly, this equivalence is given by
\begin{align*}
\FunL_{\Nuc(\Bbbk)}(\Nuc(A),\Nuc(B))&\overset{\simeq}{\longrightarrow}{_A}\mathrm{BMod}_B(\Nuc(\Bbbk))\\
F&\mapsto F(A).
\end{align*}
In particular the diagram
    \begin{equation*}
    \begin{tikzpicture}[scale=0.75]
        \node (a) at (-3,2){$\Nuc(A)$};
        \node (b) at (3,2){$\Nuc(B)$};
        \node (c) at (-3,0){$\Mod_A(\Nuc(\Bbbk))$};
        \node (d) at (3,0){$\Mod_B(\Nuc(\Bbbk))$};
        \draw[->,font=\scriptsize] (a) to node[above]{$f^\ast$}(b);
        \draw[->,font=\scriptsize] (c) to node[rotate=90,above]{$\simeq$}(a);
        \draw[->,font=\scriptsize] (d) to node[rotate=90,above]{$\simeq$}(b);
        \draw[->,font=\scriptsize] (c) to node[below]{$-\widehat\otimes_AB$}(d);
    \end{tikzpicture}
    \end{equation*}
    yields two functors $\Mod_A(\Nuc(\Bbbk))\to\Mod_B(\Nuc(\Bbbk))$, both mapping the object $A$ to the object $B$ with the $(A,B)$-bimodule structure induced by the map $A\to B$. It follows that the two functors must be equivalent; in particular, the above diagram is commutative up to natural equivalence.
\end{proof}
\begin{remark}\label{rmk: comparison}
Applying \cref{prop: pullback of nuclear categories as base change} to the case $A=B$, we see that the endofunctor $\Nuc(A)\to\Nuc(A)$ given by the tensor product with a nuclear $A$-module $M$ corresponds to the functor $-\widehat\otimes_AM\colon\Mod_A(\Nuc(\Bbbk))\to\Mod_A(\Nuc(\Bbbk))$. In particular, the equivalence of \cref{prop: nuclear categories and module categories} is symmetric monoidal.

In the following, we will blur the distinction between the functor $\Nuc(-)$ and the functor $\Mod_{(-)}(\Nuc(\Bbbk))$, since the previous statements imply that they are naturally equivalent.
\end{remark}
\begin{proposition}[\cite{mikami2023fppf}]\label{prop:descendable}
    Let $A \to B$ an fppf map of affinoid $\Bbbk$-algebras. Then the object $B$ is descendable in $\scrD(A)$, in the sense of \cref{def: descendable}. In particular $B$ is descendable as an algebra in $\Nuc(A)$.
\end{proposition}
\begin{proof}
    In \cite[Theorem 4.15]{mikami2023fppf}, the statement is proved using the definition of descendability given in \cite[Definition 2.6.1]{mann2022p}. So it is sufficient to observe that that definition implies the one given in \cref{def: descendable}.
\end{proof}
We now endow the $\infinity$-category of affinoid $\Bbbk$-algebra with the Grothendieck topology generated by fppf covers. 
\begin{corollary}\label{cor: Nuc satisy descent}
    Let $\varphi\colon A \to B$ an fppf map of affinoid $\Bbbk$-algebras. Then the functor
    \begin{equation*}
    \Nuc(-)\colon \Afd_{\Bbbk} \longrightarrow \LPr_{\Nuc(\Bbbk)}
    \end{equation*}
    satisfies descent along $\varphi$.
\end{corollary}
\begin{proof}
    This is an immediate application of \cite[Proposition 3.22]{mathew2016galois}, using \cref{prop: nuclear categories and module categories,prop:descendable}.
\end{proof}
The previous corollary was the last brick of the theory that we needed in order to apply the general machinery developed in \cref{sec:modular_case}. In particular, we obtain the following.
\begin{theorem}\label{thm: local case}
    Let $A \to B$ be an fppf map of affinoid $\Bbbk$-algebras. Then the adjunction
    \begin{equation*}
    \adjunction{-\otimes_{\Nuc(A)}\Nuc(B)}{\LPr_{\Nuc(A)}}{\LPr_{\Nuc(B)}}{\mathrm{forget}}
    \end{equation*}
    is comonadic. In particular the functor \eqref{eq: categorifed nuc} is a sheaf for the fppf topology.
\end{theorem}
\begin{proof}
    This follows from \Cref{prop:descendable,prop:modular_case}, after the identification $\Nuc(-)\simeq\Mod_{\Nuc(-)}(\Nuc(\Bbbk))$ of \Cref{rmk: comparison}.
\end{proof}
\subsection{\texorpdfstring{$1$}{1}-affineness for categories of nuclear modules over rigid varieties}
\label{sec:1_affineness_nuclear}
After proving that the functor $\LPr_{\Nuc(-)}$ is a sheaf for the fppf topology on affinoid $\Bbbk$-algebras, we are ready to tackle \cref{question:main} for more general rigid analytic varieties. We first set some piece of terminology.
\begin{notation}
Let $\mathrm{Rig}_{\Bbbk}$ denote the $\infinity$-category of \textit{discrete} rigid analytic varieties over our complete non-archimedean base field $\Bbbk$. Since $\Nuc(-)$ is a sheaf for the fppf topology on affinoid $\Bbbk$-algebras, it can be extended to a sheaf for the fppf topology on all $\mathrm{Rig}_{\Bbbk}$, i.e., we have
\[
\Nuc(-)\colon\mathrm{Rig}_{\Bbbk}\longrightarrow\CAlg(\LPr_{\Nuc(\Bbbk}).
\]
We shall denote with $\Nuccat$ the categorification of the above functor, as in \cref{def:categorification}. More precisely, this is the right Kan extension along the inclusion  $\Afd_{\Bbbk}\subseteq\mathrm{Rig}_{\Bbbk}$ of the functor
\begin{align*}
\Nuccat\colon\Afd_{\Bbbk}&\longrightarrow\catrex\\
\Spa(A)&\mapsto\LPr_{\Nuc(A)}.
\end{align*}
In virtue of \cref{thm: local case} and \cref{prop:modular_case}, $\Nuccat$ is a sheaf for the fppf topology on rigid analytic varieties. 
\end{notation}
\begin{warning}
\label{warning:derived}
    \cref{assumption:modular_final} \textit{almost} holds for the case of rigid analytic varieties and nuclear modules, but not quite. The deal-breaker of the theory is that the functor $\Nuc(-)$ will not send finite limits of affinoid spaces to colimits of $\Nuc(\Bbbbk)$-linear presentably symmetric monoidal $\infinity$-categories: indeed, in general the pullback $X\times_YZ$ of affinoid spaces needs to be derived. The theory then can work in the context of derived rigid geometry as defined in \cite{soor2024sixfunctorformalismquasicoherentsheaves}, or using derived Tate adic spaces in the sense of \cite{camargo2024analytic}.\\
    Still, pullbacks are sent to tensor products when we restrict to \textit{flat} morphisms of affinoid spaces (indeed, in this case we do not have to derive the tensor product of their coordinate affinoid $\Bbbk$-algebras); in particular, pullbacks along analytic open immersions are sent to tensor products of $\infinity$-categories of nuclear modules. So, the results in \cref{sec:modular_case} still hold under these stricter assumptions.
\end{warning}
We can now state our main result.
\begin{theorem}\label{thm: 1-affine Nuc}
    Let $X$ be a quasi-compact and quasi-separated rigid analytic variety over $\Bbbk$. Then $X$ is $1$-affine with respect to the sheaf of $\infinity$-categories $\Nuccat(-)$. More precisely, there is an equivalence of $\infinity$-categories
    \begin{equation*}
        \LPr_{\Nuc(X)} \simeq \lim_{\Spa(A) \subset X}\LPr_{\Nuc(A)},
    \end{equation*}
    where the limit is taken over an affinoid cover of $X$. 
\end{theorem}
\cref{thm: 1-affine Nuc} is a particular instance of \cref{thm:modular_case}. Indeed, the proof of \cref{thm: 1-affine Nuc} boils down to proving that the assumptions of \cref{thm:modular_case} are satisfied. 

Recall that the proof of \cref{thm:modular_case} relies on various auxiliary results, among which some variants of categorical K\"unneth formulas. Since they might be of independent interest, we explicitly write the statements down when applying the formalism of \cref{sec:modular_case} to the particular case of nuclear modules and rigid analytic varieties.
\begin{proposition}\label{lem: right adjoint general case}
    Let $\iota\colon U \to X$ and $j\colon\Spa(B)\to X$ be morphisms of rigid analytic varieties, and assume that both $U$ and $X$ are quasi-compact and separated. Assume that either $\iota$ or $j$ is flat.
    \begin{propenum}
        \item\label{lemma:base change 1-nuclear} (\emph{Base change}, \cref{lemma:base change 1}) Consider the diagram
               \[
        \begin{tikzpicture}[scale=0.75]
        \node (a) at (-3,2){$\Spa(A)\times_{X}U$};
        \node (b) at (3,2){$\Spa(A)$};
        \node (c) at (-3,0){$U$};
        \node (d) at (3,0){$X.$};
        \draw[->,font=\scriptsize](a) to node[above]{$\iota'$}(b);
        \draw[->,font=\scriptsize](b) to node[right]{$j$}(d);
        \draw[->,font=\scriptsize](c) to node[above]{$\iota$}(d);
        \draw[->,font=\scriptsize](a) to node[left]{$j'$}(c);
        \end{tikzpicture}
       \]
        Then for every nuclear module $M$ over $U$ the canonical map 
        \begin{equation*}
            j^{\ast}\iota_{\ast}M \longrightarrow \iota'_{\ast}j'^{\ast}M
        \end{equation*}
        is an equivalence.
        \item\label{lemma:pushforward nuclear colimits-nuclear} (\cref{lemma:pushforward nuclear colimits}) When restricted to nuclear modules, the functors $\iota_\ast$ and $j_\ast$ are conservative and commute with colimits. In particular, they are both monadic and comonadic functors. The monad $j_\ast j^\ast$ on $\Nuc(X)$ is naturally identified with the monad $-\widehat\otimes_{\Nuc(X)}j_\ast A$, and thus we have an equivalence of $\infinity$-categories
        \[
        \Mod_{j_\ast A}(\Nuc(X))\simeq\Nuc(A).
        \]
        \item If $\iota$ is an analytic open immersion, then $\iota_\ast$ is moreover fully faithful.
        \item \label{prop: kunneth for Nuc using opens_nuclear} (\emph{Categorical K\"unneth formula}, \cref{prop: kunneth for Nuc using opens}) The canonical functor
    \begin{equation*}
        \Nuc(U)\otimes_{\Nuc(X)}\Nuc(B) \longrightarrow \Nuc(U\times_{X}\Spa(B)) 
    \end{equation*}
    is an equivalence.
    \item\label{cor: projection formula of open-nuclear} (\emph{Projection formula}, \cref{cor: projection formula of open}) For every $N$ in $\Nuc(U)$ and $M$ in $\Nuc(X)$, the natural map
    \[
    \iota_\ast(N\otimes_{\Nuc(U)}\iota^\ast M)\longrightarrow\iota_\ast N\otimes_{\Nuc(X)}M
    \]
    is an equivalence.
    \item \label{thm: Nuc is rigid-nuclear}(\cref{thm: Nuc is rigid}) The $\infinity$-category $\Nuc(X)$ is rigid.
    \end{propenum}
\end{proposition}
In order to apply \cref{thm:modular_case}, the only thing that is not a formal consequence of the theory is that analytic open immersions of affinoid subsets $\iota\colon \Spa(A)\to X$ provide an adjunction $\iota_!\dashv \iota^\ast$, which is key in order to apply \cref{lem: right adjointable modular case}. This is the content of the following proposition.
\begin{proposition}
    \label{lem: left adjoint pullback opens-nuclear}
    Let $\iota \colon U \to X$ be an analytic open immersion into a separated and quasi-compact rigid analytic variety. Then the functor
    \begin{equation*}
        \iota^{\ast}\colon \Nuc(X) \longrightarrow \Nuc(U)
    \end{equation*}
    admits a left adjoint $\iota_{!}$.
\end{proposition}
\begin{proof}
    We will show that $\iota^\ast$ commutes with limits, i.e., that for any diagram $I\to\Nuc(X)$ selecting objects $(M_i)_{i\in I}$ the canonical map
    \[
    \lim_{i\in I}\iota^\ast M_i\longrightarrow\iota^\ast\lp\lim_{i\in I}M_i\rp
    \]
    is an equivalence.
    If $X$ is not affinoid, we can still consider a finite affinoid cover $\left\{V_j\to X\right\}_{j\in J}$. Using fppf descent of $\Nuc(-)$ we write
    \[
    \Nuc(X)\simeq\lim_{j\in J} \Nuc(V_j),
    \]
    so each $M_i$ can be realized as a limit of a cosimplicial diagram $M^{\bullet}_i$, associated to the \v{C}ech nerve of the cover $\left\{V_j\to X\right\}$. The limit $\lim_I\iota^\ast M_i$ is then computed as the totalization of the cosimplicial diagram $\lim_JM^{\bullet}_j$; so it is sufficient to prove the statement in the case $U\to X$ is an analytic open immersion into an affinoid variety.\\
    Considering a rational covering of $U$ and using once again the fact that $\Nuc(-)$ is a sheaf for the fppf topology, we can assume that $U$ is a rational subset. In virtue of \cite[Lemma 5.3 (b)]{meyer2024derived} the inclusion of a rational open subset $U\hookrightarrow X$ is also  an open immersion in the sense of \cref{def:open immersion}, so we obtain an adjunction at the level of derived $\infinity$-categories
    \begin{equation*}
        \adjunction{\iota_{!}}{\scrD(U)}{\scrD(X)}{\iota^{\ast}}
    \end{equation*}
    where $\iota_{!}$ is moreover fully faithful and satisfies the projection formula. In particular, for every $N$ and $M$ in $\scrD(U)$ we have that
    \begin{equation*}
        \iota_{!}(N\otimes_{\scrD(U)} M) \simeq \iota_{!}(N\otimes_{\scrD(U)}\iota^{\ast}\iota_{!}M) \simeq \iota_{!}N \otimes_{\scrD(X)} \iota_{!}M. 
    \end{equation*}
    So $\iota_!$ is actually a strongly monoidal colimit-preserving functor, and thus using \cite[Theorem 2.9 (c)]{meyer2024derived} we can conclude that $\iota_{!}$ preserves nuclear subcategories. In particular, $\iota^\ast\colon\scrD(X)\to\scrD(U)$ restricts to a right adjoint $\iota^\ast\colon\Nuc(X)\to\Nuc(U)$ and thus obviously commutes with limits.
\end{proof}
\begin{remark}
In the context of derived rigid analytic varieties developed in \cite{soor2024sixfunctorformalismquasicoherentsheaves}, we still have that an analogue of \cref{lem: left adjoint pullback opens-nuclear} holds when $\iota$ is a \textit{Zariski-open immersion} (\cite[Proposition 4.5]{soor2024sixfunctorformalismquasicoherentsheaves}). In particular, one could get a $1$-affineness result also for quasi-compact and separated \textit{derived} rigid analytic varieties with respect to the sheaf of presentably symmetric monoidal $\infinity$-categories $\Nuc(-)$ on derived rigid analytic varieties: indeed, $\Nuc(\Bbbk)$ is rigid thanks to \cite[Corollary 6.38]{kelly2025localisinginvariantsderivedbornological}.
\end{remark}
We conclude by presenting some consequences of \cref{thm: 1-affine Nuc}.

\begin{corollary}
\label{corollary: auxiliary kunneth}
    Let $X \to Z$ and $Y \to Z$ be two maps of quasi-compact and separated rigid analytic varieties, and assume that one of the two maps is flat. Then the diagram of $\infinity$-categories 
    \begin{equation}\label{eq: diagramma right adjointable}
    \begin{tikzpicture}[scale=0.75,baseline=(current  bounding  box.center)]]
    \node (a) at (-4,2.5){$\LPr_{\Nuc(Z)}$};
    \node (b) at (4,2.5){$\LPr_{\Nuc(Y)}$};
    \node (c) at (-4,0){$\LPr_{\Nuc(X)}$};
    \node (d) at (4,0){$\LPr_{\Nuc(X\times_ZY)}$};
    \draw[->,font=\scriptsize] (a) to node[above]{$-\otimes_{\Nuc(Z)}\Nuc(Y)$} (b);
    \draw[->,font=\scriptsize] (a) to node[left]{$-\otimes_{\Nuc(Z)}\Nuc(X)$} (c);
    \draw[->,font=\scriptsize] (b) to node[right]{$-\otimes_{\Nuc(Y)}\Nuc(X\times_ZY)$} (d);
    \draw[->,font=\scriptsize] (c) to node[below]{$-\otimes_{\Nuc(X)}\Nuc(X\times_ZY)$} (d);
    \end{tikzpicture}
    \end{equation}
    is vertically right adjointable.
\end{corollary}
\begin{proof}
    The proof can be obtained by adapting the argument of \cite[Proposition 3.2.14]{stefanich2023tan}.  We want to show that for every $\scrC$ in $\LPr_{\Nuc(X)}$ the canonical functor
    \begin{equation*}
         \scrC \otimes_{\Nuc(Z)}\Nuc(Y)\longrightarrow\scrC \otimes_{\Nuc(X)}\Nuc(X\times_{Z}Y)
    \end{equation*}
    is an equivalence. Consider an affinoid open cover of $X$, so as to write
    \[
    \LPr_{\Nuc(X)}\simeq\lim_{\Spa(A_i)\subseteq X}\LPr_{\Nuc(A_i)}
    \]
    using \cref{thm: 1-affine Nuc}.  Notice that all the forgetful functors and the functors $-\otimes_{\Nuc(X)}\Nuc(X\times_{Z}Y)$ and $-\otimes_{\Nuc(Z)}\Nuc(Y)$ commute with limits (the latter in virtue of \Cref{thm: Nuc is rigid-nuclear}), so we can reduce ourselves to the case when $X$ is affinoid. Now we apply (the dual of) \cite[Lemma 2.2.7]{stefanich2023tan} as in the proof of \cite[Proposition 3.2.14]{stefanich2023tan} to reduce to the case where both $Z$ and $Y$ are affinoid as well. In this case, the statement follows from the fact that the diagram
    \begin{equation*}
    \begin{tikzpicture}[scale=0.75]
    \node (a) at (-3,2.5){$\Nuc(Z)$};
    \node (b) at (3,2.5){$\Nuc(Y)$};
    \node (c) at (-3,0){$\Nuc(X)$};
    \node (d) at (3,0){$\Nuc(X\times_ZY)$};
    \draw[->,font=\scriptsize](a) to node[above]{}(b);
    \draw[->,font=\scriptsize](a) to node[above]{}(c);
    \draw[->,font=\scriptsize](b) to node[above]{}(d);
    \draw[->,font=\scriptsize](c) to node[above]{}(d);
    
    \end{tikzpicture}
    \end{equation*}
    is a pushout diagram of commutative algebras in $\LPr$ (\Cref{prop: kunneth for Nuc using opens_nuclear}). 
\end{proof}
\begin{porism}\label{cor: Kunneth nuclear rigid}
Applying \cref{corollary: auxiliary kunneth} to the categorical $\Nuc(X)$-module $\Nuc(X)$, we obtain an equivalence of categorical $\Nuc(Y)$-modules
\[
\Nuc(X\times_ZY)\simeq\scrC\otimes_{\Nuc(Z)}\Nuc(Y),
\]
which represents a stronger variant of the categorical K\"unneth formula for nuclear modules of \cref{prop: kunneth for Nuc using opens_nuclear} (indeed, in this case we do not assume that the source of one of the two maps is affinoid). We still remark that we can obtain a different proof of \Cref{cor: Kunneth nuclear rigid} without using the notion of $1$-affineness, by adapting the arguments used in the proof of \Cref{prop: kunneth for Nuc using opens}. 
\end{porism}
\printbibliography
\end{document}